%% file: bodrgr09_2.tex
\documentclass[12pt, a4paper]{amsart}

\newtheorem{theorem}{Theorem}[section]
\newtheorem{prop}[theorem]{Proposition}

\newtheorem{lemma}[theorem]{Lemma}
\newtheorem{theo}[theorem]{Theorem}

\theoremstyle{definition}
\newtheorem{defin}[theorem]{Definition}

\newtheorem{exam}[theorem]{Example}

\theoremstyle{remark}
\newtheorem{remk}[theorem]{Remark}
\numberwithin{equation}{section}

\usepackage{hyperref}
\usepackage[all]{xy}
\usepackage{xypic}
\input xy \xyoption{all}
\xyoption{arc}

\usepackage{multirow}
\usepackage{amsfonts, bbm, dsfont,mathrsfs,yfonts, amsmath,amssymb, amsopn, stmaryrd }

\DeclareMathOperator{\rank}{rank}

\DeclareMathOperator{\Hom}{Hom}
\DeclareMathOperator{\End}{End}

\DeclareMathOperator{\GL}{GL}
\DeclareMathOperator{\SL}{SL}
\DeclareMathOperator{\Mat}{Mat}

\DeclareMathOperator {\Br}{\mathsf{Br}}
\DeclareMathOperator{\VB}{\mathsf{VB}}
\DeclareMathOperator{\vb}{\mathsf{VB}}
\DeclareMathOperator{\Sch}{\mathsf{Sch}}
\DeclareMathOperator{\Coh}{\mathsf{Coh}}
\DeclareMathOperator{\Perf}{\mathsf{Perf}}
\DeclareMathOperator{\PIC}{\mathsf{Pic}}

\DeclareMathOperator{\TF}{\mathsf{TF}}

\DeclareMathOperator{\Tr} {\mathsf{Tr}}

\DeclareMathOperator{\cmod} {\mathsf{mod}}

\DeclareMathOperator{\T} {\mathsf{T}}

\DeclareMathOperator{\MP}{\mathsf{MP}}
\DeclareMathOperator{\BT}{\mathsf{BT}}
\DeclareMathOperator{\BC}{\mathsf{BC}}

\def \Id {\mathbb I}
\def \one {\mathbbm 1}
\def\F {\mathbb F}

\newcommand{\Pp}{\mathcal{P}}
\newcommand{\PP}{\mathbb{P}}
\newcommand{\XX}{ X}
\newcommand{\ZZ}{S}
\newcommand{\EE}{E}

\def\mP{\mathbb P}
\def\mZ{\mathbb Z}

\def\mN{\mathbb N}
\def\kE{\mathbb E}

\newcommand{\kI}{\mathcal I}
\newcommand{\kM}{\mathcal M}

\newcommand{\ee}{\mathcal{E}}
\newcommand{\ff}{\mathcal{F}}

\newcommand{\kl}{\mathcal{L}}

\newcommand{\jj}{\mathcal{J}}

\newcommand{\mm}{\mathcal{M}}
\newcommand{\oo}{\mathcal{O}}

\newcommand{\xn}{\widetilde{X}}
\newcommand{\sn}{\widetilde{S}}

\newcommand{\fn}{\widetilde{\mathcal{F}}}
\newcommand{\en}{\widetilde{\mathcal{E}}}
\newcommand{\on}{\widetilde{\mathcal{O}}}
\newcommand{\jn}{\widetilde{\mathcal{J}}}

\def \aa {\textswab{a}}

\def\d{\partial}

\def\kA{\mathfrak A}

\def\Mm{\mathscr{M}}

\def \mn {\tilde{\mu} }
\def \pn {\tilde{\pi} }
\def \i {\imath}
\def \iw {\tilde{\i} }
\def \e {\varepsilon}

\def \mo {\rightarrow}
\def \lar {\longrightarrow}
\def \emb{\hookrightarrow}
\def \to {\mapsto}

\newcommand{\ko}{{\mathbbm k}}
\newcommand{\dd}{{\mathbbm d}}
\def \bd{\bar{d}}
\newcommand{\sss}{{\mathbbm s}}
\newcommand{\rr}{{\mathbbm r}}
\def \Sum{\mathop{\textstyle{\sum}}}

\newcommand{\bul}{\scriptstyle \bullet}
\def \ci {\cdot}
\def\Arr {\xy*{}\ar@{=}@{*{}-*{\dir2{>}}}(10,0)\endxy}
\def\rrA {\xy*{}\ar@{=}@{*{}-*{\dir2{>}}}(-10,0)\endxy}

\newcommand{\mk}[1]{\makebox[16pt][c]{\rule{0pt}{11.5pt} ~$ #1 $ }}
\newcommand{\fmk}[1]{
\begin{array}{@{}|@{}c@{}|@{}}
 \hline
\mk{#1} \\
\hline
\end{array}
}
\newcommand{\mmk}[1]{\makebox[25pt][c]{\rule{0pt}{20pt}~$ ^{\displaystyle #1} $ }}

\newcommand{\cmk}[1]{\makebox[25pt][c]{$ {\displaystyle #1} $ }}
\newcommand{\rmk}[1]{\makebox[5pt][c]{\rule{0pt}{20pt}~$ ^{\displaystyle #1} $ }}

\newcommand{\smk}[1]{\makebox[15pt][c]{\rule{0pt}{5pt} ~$ #1 $ }}
\newcommand{\tmk}[1]{\makebox[8pt][c]{\rule{0pt}{11pt} ~$ #1 $ }}

\def\M{{M}}
\newcommand{\A}{{A}}
\newcommand{\B}{{B}}

\newcommand{\D}{{D}}

\begin{document} 

\title{Simple vector bundles on plane degenerations 
\linebreak
of an elliptic curve}
%

\author{Lesya Bodnarchuk} 
\address{Max-Planck-Institut f\"ur Mathematik, Bonn}
\email{lesyabod@mpim-bonn.mpg.de}
\

\thanks{
We express our sincere thanks to Professor Serge Ovsienko
for fruitful discussions and helpful advice.
The first named author would like to thank the  
Mathematisches Forschungsinstitut Oberwolfach,
where she stayed as a Leibniz fellow
during the period that this paper was written.}

\author{Yuriy Drozd}
\address{Institute of Mathematics, National Academy of Sciences of Ukraine}
\email{drozd@imath.kiev.ua}


\author{Gert-Martin Greuel}
\address{University of Kaiserslautern}
\email{greuel@mathematik.uni-kl.de}






\begin{abstract}
In 1957  Atiyah classified simple and indecomposable 
vector bundles on an  elliptic curve.
In this article we generalize his classification 
by describing the simple vector bundles on  all reduced plane 
cubic curves. Our main result states that a simple vector 
bundle on such a curve is completely determined  by its rank,
multidegree and determinant.  Our approach, based on the representation
theory of boxes, also
yields an explicit description of the corresponding
universal families of simple  vector bundles. 
\end{abstract}

\maketitle


\input{intro.tex}
\input{mainpart.tex}

\input{2comp.tex}
\input{3comp.tex}
\input{exams.tex}

\section{Properties of simple vector bundles} 
\label{sec_pr}
\subsection{Tensor products} 
Let  ${\ee(\lambda)\in \VB_E^s(r,\dd)}$ 
and $\kl(\lambda)\in \PIC^{(0,\dots, 0)}_E$
be respectively a simple 
vector bundle and the line bundle 
with the matrix $M=M(\lambda)\in \Br_{\kA}$ and
the parameter $\lambda \in \Lambda,$
where $\Lambda:= \ko^*$ if $E$ a Kodaira cycle and  $\ko$ 
if $E$ is a Kodaira fiber ($\Lambda\cong \PIC^{(0,\dots,0)}(E)$). 
\begin{prop}
For $\lambda_1,\lambda_2\in \Lambda$ we have
\begin{align*}
 \ee(\lambda_1)\otimes \kl(\lambda_2) 
=
\left\{
\begin{array}{ll}
\ee(\lambda_1\ci\lambda_2^r) 
& 
\hbox{if $E$ is a Kodaira cycle I$_1,$ I$_2$ or I$_3$;} 
\\
\ee(\lambda_1+r\ci \lambda_2)
&
\hbox{if $E$ is a Kodaira fiber II, III or IV.} 
\end{array}
\right .
\end{align*}

\end{prop}

\begin{proof}
Let $(\fn, V, \mn'(\lambda_1))$ 
and $(\on,\oo_S, \mn''(\lambda_2))$ 
be the triples of the 
vector bundle $\ee(\lambda_1)$  and
the line bundle $\kl(\lambda_2).$
Then the triple of the 
vector bundle $\ee(\lambda_1)\otimes_{\oo}\kl(\lambda_2)$ 
is $(\fn, V, \mn),$ where \linebreak
${\mn:=\mn'(\lambda_1)\otimes \mn''(\lambda_2))}.$

\noindent
For Kodaira fiber I:
$\oo_S=\ko$ and $\oo_{\sn}=\ko\oplus \ko,$
 $\mn'(\lambda_1))=(\Id,\, M({\lambda_1}))$ 
and
$\mn''(\lambda_2))= (1,\,  \lambda_2) .$  
\begin{align*}
\mn=\mn'(\lambda_1)\mathop\otimes\limits_{\oo_{\sn}} \mn''(\lambda_2))
&=
\big( \Id,\, M({\lambda_1} \big)\ci (1,\, \lambda_2)\\
&=
(\Id,\, \lambda_2\ci M(\lambda_1)\big)
= (\Id, \, M(\lambda_1\ci\lambda_2^r)\big).
\end{align*}
To obtain the last equality one should reduce 
$\lambda_2\ci M(\lambda_1)$ to the canonical form preserving 
the first $\Id$-matrix. We illustrate it on the case $r=2:$
\begin{align*}
\big(\mn(0),\mn(\infty) \big )
& =
\left(
\begin{array}{|@{}     c       @{}|@{}     c       @{}|}
\hline
{\mk{1}}  &    {\mk{0}}    \\
\hline
 \mk{0}        &        {\mk{1}} \\
\hline
\end{array}
,\;
\lambda_2 
\begin{array}{ |@{}     c       @{}|@{}     c       @{}|}
\hline
 \mk{0}  &    \mk{1}   \\
\hline
 \mk{\lambda_1}      &        {\mk{0}} \\
\hline
\end{array}
\right)
=
\left(
\begin{array}{|@{}     c       @{}|@{}     c       @{}|}
\hline
{\mk{1}}  &    {\mk{0}}    \\
\hline
 \mk{0}        &        {\mk{1}} \\
\hline
\end{array}
,\;
\begin{array}{ |@{}     c       @{}|@{}     c       @{}|}
\hline
 \mk{0}  &    \mk{\lambda_2}   \\
\hline
{\lambda_1 \lambda_2}      &        {\mk{0}} \\
\hline
\end{array}
\right)
\\
&=
\left(
\begin{array}{|@{}     c       @{}|@{}     c       @{}|}
\hline
{\mk{\frac{1}{\lambda_2}}}  &    {\mk{0}}    \\
\hline
 \mk{0}        &        {\mk{1}} \\
\hline
\end{array}
,\;
\begin{array}{ |@{}     c       @{}|@{}     c       @{}|}
\hline
 \mk{0}  &    \mk{1}   \\
\hline
 { \lambda_1 \lambda_2}      &        {\mk{0}} \\
\hline
\end{array}
\right)
=
\left(
\begin{array}{|@{}     c       @{}|@{}     c       @{}|}
\hline
{\mk{1}}  &    {\mk{0}}    \\
\hline
 \mk{0}        &        {\mk{1}} \\
\hline
\end{array}
,\;
\begin{array}{ |@{}     c       @{}|@{}     c       @{}|}
\hline
 \mk{0}  &    \mk{1}   \\
\hline
{\lambda_1 \lambda_2^2}      &        {\mk{0}} \\
\hline
\end{array}
\right).
\end{align*}

\noindent
For Kodaira fiber II:
$\oo_S=\ko$ and $\oo_{\sn}=\ko[\e]/{\e^2},$
 $\mn'(\lambda_1)=\Id+ \e\ci M({\lambda_1})$ 
and
$\mn''(\lambda_2)= 1 + \e\cdot \lambda_2 .$  
\begin{align*}
\mn=\mn'(\lambda_1)\mathop\otimes\limits_{\oo_{\sn}} \mn''(\lambda_2))
&=
\big(\Id+ \e\ci M({\lambda_1})\big)\ci (1 + \e\cdot \lambda_2)
\\
&=
\Id+ \e \ci\big( M({\lambda_1})+ \lambda_2 \ci\Id \big)
=
\Id+ \e \ci M(\lambda_1+ \lambda_2).
\end{align*}
The last equality follows immediately if we
rewrite $M(\lambda)$ in the ``diagonal'' form (\ref{form_diag}). 
For example, if $r=2$ and $d=1$ we have
\begin{align*}
\mn=\mn(0)+ \e\mn_{\e}(0)
& =
%
\begin{array}{|@{}     c       @{}|@{}     c       @{}|}
\hline
{\mk{1}}  &    {\mk{0}}    \\
\hline
 \mk{0}        &        {\mk{1}} \\
\hline
\end{array}
+\e\ci
\left(
\begin{array}{ |@{}     c       @{}|@{}     c       @{}|}
\hline
 \mk{\frac{\lambda_1}{2}}  &    \mk{1}   \\
\hline
 \mk{0}      &        {\mk{\frac{\lambda_1}{2}}} \\
\hline
\end{array}
+
\begin{array}{ |@{}     c       @{}|@{}     c       @{}|}
\hline
 \mk{\lambda_2}  &    \mk{0}   \\
\hline
 \mk{0}      &        \mk{\lambda_2} \\
\hline
\end{array}
\right)
=
\begin{array}{|@{}     c       @{}|@{}     c       @{}|}
\hline
{\mk{1}}  &    {\mk{0}}    \\
\hline
 \mk{0}        &        {\mk{1}} \\
\hline
\end{array}
+\e
\begin{array}{ |@{}     c       @{}|@{}     c       @{}|}
\hline
 {\frac{\lambda_1+ 2\lambda_2}{2}}  &    \mk{1}   \\
\hline
\mk{0}      &        {\frac{\lambda_1+ 2\lambda_2}{2}} \\
\hline
\end{array}
.
\end{align*}

\noindent
For Kodaira cycles I$_2,$ I$_3$ and fibers III and IV 
the calculations should be carried out on each component.
On the first component the picture is similar to the cases
of I and II. On the other components we have $\mn_k =\mn_k'.$
\end{proof}

\begin{exam}
If $r=3$ and $\dd=(1,1)$ 
for Kodaira cycle I$_2$ we have
\begin{align*}
\left(
\begin{array}{@{}c@{}}
\begin{array}{ |@{}c @{}@{} c   @{}|@{}   c   @{}|}
\hline
{\mk{1}}      &\mk{0}         &     \mk{0}    \\
 {\mk{0}}     &{\mk{1}}    &     \mk{0}      \\
\hline
 {\mk{0}}     &{\mk{0}}    &     \mk{1}    \\
\hline
\end{array}
\;\;
\begin{array}{ |@{}c @{}@{} c   @{}|@{}   c   @{}|}
\hline
{\mk{0}}      &\mk{1}         &     \mk{0}    \\
 {\mk{0}}     &{\mk{0}}    &     \mk{1}      \\
\hline
 {\mk{\lambda_1}}     &{\mk{0}}    &     \mk{0}    \\
\hline
\end{array}
\\
\rule{0pt}{33pt}
\begin{array}{@{}c@{}}
\begin{array}{ |@{}c @{} c   @{}|@{}   c   @{}|}
\hline
{\mk{1}}      &\mk{0}         &     \mk{0}    \\
{\mk{0}}     &{\mk{0}}    &     \mk{1}      \\
\hline
 {\mk{0}}     &{\mk{1}}    &     \mk{0}    \\
\hline
\end{array}
\;\;
\begin{array}{ |@{}c @{} c   @{}|@{}   c   @{}|}
\hline
{\mk{1}}      &\mk{0}         &     \mk{0}    \\
 {\mk{0}}     &{\mk{1}}    &     \mk{0}      \\
\hline
 {\mk{0}}     &{\mk{0}}    &     \mk{1}    \\
\hline
\end{array}
\end{array}
\end{array}
\right )
\cdot
\left(
\begin{array}{@{}c@{}}
\fmk{1}\; \fmk{\lambda_2}\\
\rule{0pt}{13pt}
\fmk{1}\; \fmk{1}
\end{array}
\right)
=
\left(
\begin{array}{@{}c@{}}
\begin{array}{ |@{}c @{}@{} c   @{}|@{}   c   @{}|}
\hline
{\mk{1}}      &\mk{0}         &     \mk{0}    \\
 {\mk{0}}     &{\mk{1}}    &     \mk{0}      \\
\hline
 {\mk{0}}     &{\mk{0}}    &     \mk{1}    \\
\hline
\end{array}
\;\;
\begin{array}{ |@{}c @{}@{} c   @{}|@{}   c   @{}|}
\hline
{\mk{0}}      &\mk{1}         &     \mk{0}    \\
 {\mk{0}}     &{\mk{0}}    &     \mk{1}      \\
\hline
 {\lambda_1\lambda_2^3}     &{\mk{0}}    &     \mk{0}    \\
\hline
\end{array}
\\
\rule{0pt}{33pt}
\begin{array}{@{}c@{}}
\begin{array}{ |@{}c @{} c   @{}|@{}   c   @{}|}
\hline
{\mk{1}}      &\mk{0}         &     \mk{0}    \\
 {\mk{0}}     &{\mk{0}}    &     \mk{1}      \\
\hline
 {\mk{0}}     &{\mk{1}}    &     \mk{0}    \\
\hline
\end{array}
\;\;
\begin{array}{ |@{}c @{} c   @{}|@{}   c   @{}|}
\hline
{\mk{1}}      &\mk{0}         &     \mk{0}    \\
 {\mk{0}}     &{\mk{1}}    &     \mk{0}      \\
\hline
 {\;\mk{0}\;}     &{\mk{0}}    &     \mk{1}    \\
\hline
\end{array}
\end{array}
\end{array}
\right )
\end{align*}
and for Kodaira fiber III taking 
$\lambda:={\lambda_1+3\lambda_2}$ we have
\begin{align*}
\left( 
\begin{array}{@{}c @{}}
\begin{array}{ |@{}c @{}@{} c   @{}|@{}   c   @{}|}
\hline
{\mk{1}}      &\mk{0}         &     \mk{0}    \\
 {\mk{0}}     &{\mk{1}}    &     \mk{0}      \\
\hline
 {\mk{0}}     &{\mk{0}}    &     \mk{1}    \\
\hline
\end{array}
+
\e_1
\begin{array}{ |@{}c @{}@{} c   @{}|@{}   c   @{}|}
\hline
{\mk{\lambda_1}}      &\mk{1}         &     \mk{0}    \\
 {\mk{0}}     &{\mk{0}}    &     \mk{1}      \\
\hline
 {\mk{0}}     &{\mk{0}}    &     \mk{0}    \\
\hline
\end{array}
\\
\rule{0pt}{33pt}
\begin{array}{ |@{}c @{} c   @{}|@{}   c   @{}|}
\hline
{\mk{1}}      &\mk{0}      &     \mk{0}    \\
 {\mk{0}}     &{\mk{0}}    &     \mk{1}      \\
\hline
 {\mk{0}}     &{\mk{1}}    &     \mk{0}    \\
\hline
\end{array}
+
\e_2
\begin{array}{ |@{}c @{} c   @{}|@{}   c   @{}|}
\hline
{\mk{0}}      &\mk{0}         &     \mk{0}    \\
 {\mk{0}}     &{\mk{0}}    &     \mk{0}      \\
\hline
 {\mk{0}}     &{\mk{0}}    &     \mk{0}    \\
\hline
\end{array}
\end{array}
\right)
\cdot
\left(
\begin{array}{@{}c@{}}
\fmk{1} +\e_1\fmk{\lambda_2}\\
\rule{0pt}{13pt}
\fmk{1} +\e_2\fmk{0}\\
\end{array}
\right)
=
\left( 
\begin{array}{@{}c @{}}
\begin{array}{ |@{}c @{}@{} c   @{}|@{}   c   @{}|}
\hline
{\mk{1}}      &\mk{0}         &     \mk{0}    \\
 {\mk{0}}     &{\mk{1}}    &     \mk{0}      \\
\hline
 {\mk{0}}     &{\mk{0}}    &     \mk{1}    \\
\hline
\end{array}
+
\e_1
\begin{array}{ |@{}c @{}@{} c   @{}|@{}   c   @{}|}
\hline
\mk{\lambda}      &\mk{1}         &     \mk{0}    \\
 {\mk{0}}     &{\mk{0}}    &     \mk{1}      \\
\hline
 {\mk{0}}     &{\mk{0}}    &     \mk{0}    \\
\hline
\end{array}
\\
\rule{0pt}{33pt}
\begin{array}{ |@{}c @{} c   @{}|@{}   c   @{}|}
\hline
{\mk{1}}      &\mk{0}         &     \mk{0}    \\
 {\mk{0}}     &{\mk{0}}    &     \mk{1}      \\
\hline
 {\mk{0}}     &{\mk{1}}    &     \mk{0}    \\
\hline
\end{array}
+
\e_2
\begin{array}{ |@{}c @{} c   @{}|@{}   c   @{}|}
\hline
{\mk{0}}      &\mk{0}         &     \mk{0}    \\
 {\mk{0}}     &{\mk{0}}    &     \mk{0}      \\
\hline
 {\mk{0}}     &{\mk{0}}    &     \mk{0}    \\
\hline
\end{array}
\end{array}
\right).
\end{align*}
\end{exam}

\medskip
\subsection{Morphisms }

\begin{prop}
\label{prop_morph}
Let $E$ be one of the curves from Table (\ref{list_1})
and $\ee(\lambda_1),\ee(\lambda_2) \in\VB_{E}^s(r,\dd)$
with $\lambda_1\neq \lambda_2. $
Then  
$ \Hom_E\big(\ee(\lambda_1),\ee(\lambda_2) \big)=0.$
\end{prop}

\begin{proof}
From the equivalence 
$\VB^s_E(r,\dd)\stackrel{\sim}\mo \Br_{\kA}(\sss)$
we have:
$$\Hom_E\big(\ee(\lambda_1),\ee(\lambda_2) \big)=
\Hom_{\kA}\big(M(\lambda_1),M(\lambda_2) \big).$$
  Let $(S,S')\in\Hom_{\kA}\big(M(\lambda_1),M(\lambda_2) \big).$ 
If $r=1$ and $(S,S')\neq 0$ then $S'=S\in\ko^*$ and since $M(\lambda_1)=
\lambda_1,$ and $M(\lambda_2)=\lambda_2$, we get
a contradiction:  $S\lambda_1S^{-1} =\lambda_2 .$
Recall that a path $p$ on a small reduction automaton 
gives an equivalence of the categories 
$\Br_{\kA}(\sss) \stackrel{p} \lar \Br_{\kA'}(\sss'),$
where $\sss'\leq\sss .$ Thus the statement follows
by induction on the dimension vector $\sss$ along the path $p.$
\end{proof}

\begin{remk}
By the same approach one can also describe torsion free
sheaves which are not vector bundles. We are going to 
consider this situation in further works. 
One can consult \cite{thesis} Sections 3.3, 4.5 and 7.7
about torsion free sheaves on cuspidal and tacnode curves.
\end{remk}


\end{document}

%% file: intro.tex

\section{Introduction}
\label{intro} 
The theory  of vector bundles on an elliptic curve and its
degenerations  is known to
be closely 
related with the theory of  integrable systems
(see e.g. \cite{Krichever, Manin, Mulase}). 
Another motivation  for studying  
vector bundles on elliptic fibrations 
comes from the work of Friedman, Morgan and Witten \cite{FMW},
who discovered  their  importance  
 for  heterotic string theory.
The  main motivation of  our investigation was the following 
problem. Let $\kE \rightarrow T$ be an elliptic fibration,
where $T$ is some basis such that for any point $t \in T$ the fiber
$\kE_t$ is a reduced projective curve with the trivial dualizing sheaf.
$$
{\xy 0;/r0.17pc/:
\POS(40,0);
{(25,0)\ellipse(30,10){-}},
\POS(-15,20)*++{\kE}="a";
\POS(-15,-2)*++{T}="b";
{\ar "a";"b"};
%
\POS(1,23)*++{\kE_t};
{(12,20)\ellipse(1.5,2.5)_,=:a(180){-}},
{(12,20)\ellipse(2,2.5)^,=:a(180){-}},
\POS(3,30);
@(,
\POS(3,25)@+,  \POS(8,22)@+, \POS(9,20)@+,
\POS(8,18)@+, \POS(3,15)@+, \POS(3,10)@+,
**\crvs{-}
,@i @);
%
\POS(20,30)*++{\kE_t};
{\ar@/_10pt/@{-}(22,19);(22,34)},
{\ar@/_10pt/@{-}(28,34);(28,19)},
%
\POS(40,27)*++{\kE_0};
{\ar@/_9pt/@{-}(41,16);(41,30)},
{\ar@/_9pt/@{-}(50,30);(50,16)},
%
\POS(15,8); \POS(45,0)*{\bul};
**\crv{(16,0)};
\POS(15,-8); \POS(45,0);
**\crv{(16,0)};
\POS(15,8);
{\ar@{.}(45,12); (45,0) }
\POS(25,2)*{\bul};
\POS(29,4)*++{t};
{\ar@{.}(25,2); (25,16) }
\POS(8,-2)*{\bul};
\POS(6,-2)*++{t};
{\ar@{.}(8,-2); (8,11) }
\POS(48,0)*++{0};
\endxy}
$$
In most applications, 
a generic fiber of this fibration is an elliptic curve
and for the points of the discriminant locus $\Delta \subset T$
the fibers are singular (and possibly reducible). 
Can one give   a uniform description of 
simple vector bundles both on the smooth and the singular fibers?

It is known that the category 
of \emph{all} vector bundles of a singular genus one curve $E$ essentially
depends on the singularity type of the curve. 
For example, in the case of the Weierstra\ss{} family 
$\kE \mo \mathbb{C}^2$ given by the equation 
$zy^2 = 4x^3 + g_2 xz^2 + g_3 z^3,$
the cuspidal fiber $E=\kE_{(0,0)}$ is vector-bundle-wild, whereas
all the other fibers E=$\kE_{(g_2,g_3)}$ (smooth and nodal)
are vector-bundle-tame\footnote{
In representation theory a category is called 
\emph{tame} if its indecomposable objects 
can be described by some discrete and one continuous parameters,
and \emph{wild} if they are non-classifiable. 
An algebraic variety $X$ is called vector-bundle-wild 
or  vector-bundle-tame
if the category $\VB_X$ of vector bundles on $X$ 
is wild or respectively tame (see \cite{drgr}). 
}.
This phenomenon seems to be rather strange, since
 very strong continuity results
for the  Picard functor are 
known to be true \cite{Altman2}.
It is one of the results of this paper that 
 the situation is completely different if one
restricts to the study of the \emph{simple} 
\footnote{A bundle is called simple if it admits 
no endomorphisms but  homotheties.} 
vector bundles. 
Namely we prove that the category $\VB^s_E$ of simple vector bundles
on $E$ is indeed tame.
Moreover, we provide a complete classification of simple bundles
and describe a bundle on the moduli space, 
having certain universal properties.

The starting point of our  investigation and 
the  main source of inspiration was 
the following classical  result of  Atiyah.

\begin{theorem}[\cite{Ati57}]
\label{ati}
Let $E$ be an elliptic  
curve over an algebraically closed field $\ko$. 
Then
  a simple vector bundle
$\ee$ on $E$ is 
uniquely determined by its rank $r$, degree $d,$ 
which should be coprime,
and determinant $det(\ee) \in \PIC^d(E) \cong E.$
\end{theorem}

\noindent
The main result of our article generalizes Atiyah's theorem  to all 
reduced plane degenerations of an elliptic curve.
Singular fibers of elliptic fibred  surfaces were described  by Kodaira 
and throughout  this article we make use of his classification, see for 
example 
 \cite[Table 3, p.150]{BPV}. 
In what follows the  cycles of projective lines 
(also called Kodaira cycles) are denoted by I$_N,$  
where $N$ is the number of irreducible components. 
Note that a Kodaira cycle I$_N$ is a  plane curve 
if and only if $N\leq 3.$ 
Besides them, 
there are precisely three other Kodaira fibers.
Thus, we study simple vector bundles on the following six configurations:  

\begin{table}[ht]
$$
\begin{array} {@{}| c |l  l |l  l |@{}}
\hline
N &
\multicolumn{2}{c|}{ \;\;\hbox{ Kodaira cycles} \hspace{1cm}\,}
&
\multicolumn{2}{c|}
{ \;\;\hbox{ Kodaira fibers} \hspace{1cm}\, } \\
\hline
\rule{0pt}{15pt}
 N=1 & \hbox{I}_1:& 
 y^2z=x^3+x^2z 
 &\hbox{II}:&
y^2z=x^3  
\\ 
\rule{0pt}{25pt}
&&
\hspace{0.3cm}
\begin{array}{c}
\xy /r0.15pc/:
\POS(4,18);@-
{
(6,15)@+,  (16,5)@+, (21,10)@+
,(16,15)@+, (6,5)@+, \POS(4,2)@-,
,**\qspline{}}  
\endxy
\end{array}
&&
\hspace{0.3cm}
\begin{array}{c}
{\xy /r0.15pc/:
\POS(10,18); \POS(18,10); **\crv{(11,12)};
\POS(10,2); \POS(18,10); **\crv{(11,8)};
\endxy } 
\\ 
\end{array}
\\
\hline
\rule{0pt}{15pt}
 N=2 & \hbox{I}_2:& z^3=xyz &  \hbox{III}:& y^2z=x^2y \\
&&
\rule{0pt}{25pt}
\hspace{0.5cm}
\begin{array}{c}
\xy /r0.15pc/:
{\ar@/^10pt/@{-}(42,19);(42,2)},
{\ar@/^10pt/@{-}(48,2);(48,19)},
\endxy
\end{array}
&&
\hspace{0.3cm}
\begin{array}{c}
\xy /r0.15pc/:
{\ar@/^8pt/@{-}(36,19);(36,2)},
{\ar@/^8pt/@{-}(45,2);(45,19)},
\endxy
\end{array}
\\ 
\hline
\rule{0pt}{15pt}
 N=3 & \hbox{I}_3:&  xyz=0 
 &  \hbox{IV}:& xy^2=x^2y 
 \\
\rule{0pt}{25pt}
&&
\hspace{0.25cm}
\begin{array}{c}
{
\xy /r0.15pc/:
{\ar@{-}(59.5,7);(75.5,7)},
{\ar@{-}(65,16);(73.5,3.5)},
{\ar@{-}(69,16);(61.5,3.5)},
\endxy
}
\end{array}
&&
\begin{array}{c}
{\xy /r0.15pc/:
\POS(0,0);
{\ar@{-}(60,10);(80,10)};
{\ar@{-}(64,1.5);(76,18.5)};
{\ar@{-}(76,1.5);(64,18.5)};
\endxy
}
\end{array}
\\ 
\hline
\end{array}
$$
\caption{\label{list_1}}
\end{table}

\noindent
In order to present our main theorem,  let us fix some  notations.
Throughout  this article,  let $\ko$ be an algebraically 
closed field and a \emph{curve} be a reduced projective curve.
Let $E$ be a plane degeneration of an elliptic curve,
$N=1,2,3$  the number of its irreducible components 
and $L_k $  the 
$k$-th component of $E$. 
For a vector bundle $\ee$ on $E$ we denote
\begin{itemize}
\item $d_k = d_k(\ee) = \deg(\ee|_{L_k}) \in \mathbb{Z}$ the degree of
the restriction of $\ee$ on $L_k$;
\item $\dd = \dd(\ee) = (d_1,\dots,d_N) \in \mathbb{Z}^N$ the 
\emph{multidegree}
of $\ee$;
\item $d = \deg(\ee) = d_1 + \dots + d_N$ the degree of $\ee$. 
In our cases it is equal  to the  Euler-Poincar\'{e} characteristic:
$\chi(\ee) = h^0(\ee) - h^1(\ee);$
\item $r = \rank(\ee)$ the rank of $\ee$.
\end{itemize}

\noindent 
Moreover, let $\PIC^{\dd}(E)$\footnote{
Note that $\PIC^{\dd}(E)$
is $E$ for an elliptic curve, $\ko^*$ for Kodaira cycles 
and $\ko$ for the other Kodaira fibers.}
 be the Picard group of
invertible sheaves of multidegree $\dd$ on $E .$
The following theorem generalizes Atiyah's classification and is the main 
result of this article.

\begin{theo}
\label{theo_main}
Let $E$ be a reduced plane cubic curve with  $N$ irreducible components,
$1 \le N \le 3$. 
\begin{itemize}
  \item[(i)]
 Then the  rank $r$ and the degree $d$ 
 of a simple vector bundle on $E$ are coprime.
For any 
 tuple of integers $(r, \dd) \in \mN \times
\mathbb{Z}^N$ 
such that $\gcd(r,d_1+\dots +d_N)=1$, 
let $\Mm = \VB^s_E(r, \dd)$ be the  set of simple vector bundles
of rank $r$ and multidegree $\dd$. Then the map
$\det: \Mm  \rightarrow \PIC^{\dd}(E)$ is a bijection.
    \item[(ii)]
The Jacobian   $\PIC^{(0,\dots,0)}(E)$ 
acts transitively on $\Mm$. 
The stabilizer of a point is isomorphic
to $\mathbb{Z}_r$ if $E$ is a Kodaira cycle, and is trivial 
in the remaining cases.  
\end{itemize}
\end{theo}

\noindent
Let $\Lambda := \ko^*$ if $E$ is a Kodaira cycle and $\Lambda :=\ko$ if 
$E$ is a Kodaira fiber of type II, III or IV.
By \ref{theo_main} (i) $\Lambda$ is a moduli space of 
simple vector bundles of given rank $r$ and multidegree $\dd$ 
provided $\gcd(r,d)=1.$ 
By an observation of Burban and Kreu\ss{}ler 
\cite{BurbanKreussler4},  for a given 
 tuple of integers $(r, \dd) \in \mathbb{N} \times
\mathbb{Z}^N$ such that $\gcd(r,d)=1$, our method 
yields an explicit construction of a 
 vector bundle $\Pp = \Pp(r, \dd) \in \VB_{E \times \Lambda}$ 
satisfying in the general case only the following universality properties:
\begin{itemize}
\item for any point $\lambda \in \Lambda$ the vector bundle  
$\Pp(\lambda) := \Pp|_{E \times \{\lambda\}} \in \VB(E)$ is simple
of rank $r$ and multidegree $\dd$;
\item for any vector bundle $\ee \in \VB^s_E(r, \dd)$ there exists a unique
$\lambda \in \Lambda$ such that $\ee \cong \Pp(\lambda)$;
\item for two points $\lambda \ne  \mu$   from $\Lambda$ we have
$\Pp(\lambda) \not\cong \Pp(\mu)$.
\end{itemize} 
If the curve $E$ is irreducible, the vector bundle $\Pp$ 
is the universal family of \emph{stable} vector bundles of rank 
$r$ and degree $d$.

Similarly to Atiyah's proof \cite{Ati57},  
the main ingredient of our approach  is a 
construction of various bijections 
$
\vb^s_{E}(r,\dd) \mo \vb^s_{E}(r^\prime, \dd^\prime), 
$
where $r^\prime < r$. 
However, our method is completely different from Atiyah's.
We use a reduction of our classification problem to the description
of \emph{bricks}  in the category of representations of a certain 
\emph{box} (or a differential biquiver). 
Moreover, we provide  an explicit 
algorithm (algorithm \ref{algorithm})
that for a given tuple  
$(r,\dd) \in\mN\times \mZ^N$  
constructs a canonical form of a matrix, describing 
the universal family
of simple vector bundles of rank $r$ and multidegree $\dd$.   
The core of this algorithm
is the automaton of reduction, 
which is given for each of the listed curves and operates 
on discrete parameters like Euclidean algorithm.

For a rather long time (till the middle of the 70s) 
there were no efficient methods for studying  moduli spaces 
of vector bundles of higher ranks
on singular curves. 
In order to study vector bundles on (possibly reducible)
projective curves with only nodes or cusps as singularities, 
Seshadri  introduced  the concept of the so-called parabolic bundles 
(see \cite[Section 3]{Seshadri}). This approach was later 
 developed by 
Bhosle, who introduced the notion of generalized parabolic bundles 
 \cite{Bhosle92,  Bhosle96}. 
 
%
%
%
%
%
%

Our method of studying vector bundles on genus one curves
is a certain categorification of the language of parabolic bundles
of Seshadri and Bhosle.
It was originally proposed in \cite{drgr}, see also 
\cite{BDG01} and \cite{surv} for some further elaborations. 
The idea
of this method can be explained as follows. 
Let $X$ be a 
singular reduced projective curve 
(typically rational, but with arbitrary singularities), 
$\pi:\widetilde{X} \rightarrow X$ its normalization.
Then a description of the fibers of the functor $\pi^*:  
\VB_X \rightarrow \VB_{\widetilde{X}}$ can be 
converted to some representation theory problem,
called a \emph{matrix problem}. 
The main application of this method 
concerns the case of curves of 
arithmetic genus one. 
In the case of a cycle of $N$ projective lines 
(Kodaira cycles $I_N$),
the obtained matrix problem turns out to be 
representation-tame, see \cite{Bon92} and \cite{CB89}.
As a result, it  
 allows to 
obtain a complete  classification of 
indecomposable torsion free sheaves on these genus one curves, 
see \cite{drgr}   and \cite{surv}.

However, a description of the exact combinatorics of \emph{simple}  
vector bundles on a cycle of projective lines 
  requires some extra work. This was done in \cite{BDG01}, but
the resulting answer was not very explicit. 
For the case of a nodal cubic curve $zy^2 = x^3 + x^2z$, 
in \cite{burb} Burban derived 
the statement of Theorem \ref{theo_main} using the classification 
of all indecomposable objects.
In this article we give an improved  description
of simple vector bundles 
on cycles I$_1,$ I$_2$ and I$_3$
using the
technique of the so-called  \emph{small reductions} of matrix problems.

As we have mentioned above, the representation-theoretic properties 
of the category of torsion free sheaves on  Kodaira cycles and the other
degenerations of elliptic curves are rather different. 
For example, for a cuspidal rational curve $zy^2 = x^3$ 
even the classification of indecomposable semi-stable 
vector bundles of a given slope  
is a representation-wild problem. 
However, if 
we additionally impose the 
simplicity assumption,
then the wild fragments of the  matrix problem
disappear and we can  reduce the
matrices to a canonical form (see \cite{boddr}).

The matrix  problems describing simple vector bundles
on nodal and cuspidal cubic curves are relatively easy to deal with,
since they are \emph{self-reproducing}, i.e.
after applying one step of small reduction we obtain
the same problem but with matrices of smaller sizes.
In fact, the matrix reduction operates 
on discrete parameters of vector bundles
as Euclidean algorithm. 
Carrying this out 
we obtain the statement of
Theorem \ref{theo_main}
for \emph{irreducible} 
degenerations of an elliptic curve.
Unfortunately, the matrix problems for curves with 
many components are \emph{no longer} self-reproducing.
However, they turn out to be such in  
some bigger class of matrix problems. 
To study this class in a conceptual way
we need
more sophisticated methods from representation theory.
Namely, we describe our matrix problem as the category
of representations  of a certain 
box 
(also called bocs, 
``bimodule over a category with a coalgebra structure''
or differential biquiver) see \cite{thesis, boddr2}.

The technique of boxes  
is known to be very useful  for proving tame-wild
dichotomy theorems and various semi-continuity results, see 
\cite{Dro79}, \cite{Dro01}, \cite{Dro05}, \cite{CB90}  etc.
 A new feature, 
illustrated in this article, is that the formalism of boxes
can be very efficiently 
applied for constructing canonical forms of
representations ``in general position''. 
A usual approach to a matrix problem is a consecutive 
application of a \emph{minimal edge reduction}, which is a reduction 
of a certain block to its Gau\ss{} form.
However, since we are interested in bricks 
it turns out that it is sufficient to take into account only 
\emph{small reductions}, 
which are Gau\ss{} reductions provided that 
the rank of the block is maximal.
This way for each plane singular cubic curve 
and the
matrix problem corresponding to
the family
of simple vector bundles of rank $r$ and multidegree $\dd$
we get 
an explicit algorithm constructing its canonical form. 
The course of the construction is given as a path
on some \emph{automaton}, whose  states are 
boxes and transition arrows are 
small reductions.

To put our results in a broader mathematical context 
we would like to mention that
the case of singular curves of genus one is special in many respects.
We are especially interested in the study of vector bundles
on curves having trivial dualizing bundle. This automatically implies
that they have arithmetic genus one, but not vice versa. 
In \cite{FMW} Friedman, Morgan and Witten proposed a powerful  method
of constructing vector bundles on irreducible genus one curves and elliptic
fibrations, based on the technique of the  so-called spectral covers.
Later,  it was realized that their construction can be 
alternatively described using the language of Fourier-Mukai transforms,
see e.g. \cite{BurbanKreussler1}, \cite{BBHM}, 
\cite{hispanci2}. 
Although for
irreducible cubic curves 
 Theorem \ref{theo_main} was previously known
and can be proven
using either geometric invariant theory or 
Fourier-Mukai transforms,
our approach has  one particular advantage. Namely, it yields a very 
explicit description of a universal family of simple vector bundles, which
turned out to be important in applications. In particular, 
it was used to get new  solutions of the associative and quantum 
Yang-Baxter equations,  see  \cite{Polishchuk_4} 
and  \cite[Section 9]{BurbanKreussler4}.

We should also mention that
the geometric point of view 
suggests to replace the simplicity condition by 
Simpson stability. 
Despite the fact that for reducible curves there are even line bundles 
which are not Simpson semi-stable, both notions are closely related 
for curves of arithmetic genus one.
In \cite{Lopez1} and \cite{Lopez2}
L\'opez-Martin described the geometry of the compactified 
Jacobian in case of Kodaira fibers and elliptic fibrations.

\subsection*{Organization of the material}
In Section \ref{sec_tripl} we recall
the construction of \cite{drgr}
 and replace the category of vector bundles $\VB_E$ by the 
 equivalent category
  of triples $\Tr_E.$
  Fixing bases of triples we obtain the category 
  of matrices $\MP_E.$
In Sections \ref{sec_matr_pr_cycles} and  \ref{sec_matr_pr_fibers}
this procedure is applied to all the curves from Table \ref{list_1}. 
In Section \ref{sec_simpl_cond} we study the properties 
induced by the simplicity condition and obtain some 
additional restrictions for the matrix problem $\MP_E $.
In Section \ref{sec_prim_red}
we fix discrete parameters $(r,\dd)$ 
and reduce a brick-object\footnote{A \emph{brick} or a \emph{schurian} object is 
a representation with no nonscalar endomorphisms.  }
 of $\MP_E(r,\dd)$ to its partial canonical form.
Remarkably, this new matrix problem and its dimension vector $\sss$
are completely determined by the curve $E$  the rank $r$
and the multidegree $\dd.$ 
This correspondence for the curves with many components
is given in Tables \ref{table_I_2} -- \ref{table_IV}. 
   
In Section \ref{sec_boxes} we provide 
a formal approach to the partially reduced matrix problem: 
we interpret it as the category   
of bricks $\Br_{\kA}(\sss)$ 
of some box $\kA$ and dimension vector $\sss.$
We prove that any break is a module in a general position,
thus the Gau\ss{}  reduction can be replaced by the small one.
A course of reduction can be presented as a path on some 
automaton, where the states are matrix problems and transitions are 
small reductions. 
We call a box principal if $\Br_{\kA}(\sss)\cong \VB_E^s(r,\dd).$
For fixed 
rank $r$ and multidegree $\dd,$ if 
the set $\Br_{\kA}(\sss)$ is nonempty, then 
there is a path $p:\kA \mo \kA',$ 
where $\kA'$ is principal,
reducing the  
dimension vector $\sss$ to $(1,0,\dots, 0):$
\begin{align*}
{\xymatrix{
\VB_E^s(r,\dd)\ar[d]^{\cong}
&&& \PIC^{(0,\dots,0)}_E\ar[d]^{\cong}
\\
\Br_{\kA}(\sss) \ar[rrr]^{p}_{\sim}
&&& \Br_{\kA'}(1,0,\dots,0).
}}
\end{align*}
A transition operates 
on the pair $(d\bmod r, r-d\bmod r)$ as Euclidean
algorithm and for ${\ee\in\VB_E^s(r,\dd)}$  we obtain 
$\gcd(r,d)=1.$
It turns out that this condition 
is not only necessary but also sufficient for $\VB_E^s(r,\dd)$
to be nonempty.  
The canonical form of a brick from 
$\Br_{\kA}(\sss)$ can be recovered by 
reversing the path $p.$ 
The whole procedure
is emphasized in algorithm \ref{algorithm}.

In Sections \ref{sec_small_red_node_cusp} --
\ref{sec_small_red_I_3} we construct 
automatons for each Kodaira cycle $I_N$ ($N\leq 3$) and show that
a path on it also encodes a course of reduction for the Kodaira fiber 
with $N$-components. 

Analyzing how 
a path operates on the dimension vector $\sss$ we deduce 
the first part of Theorem
\ref{theo_main}.
In Section \ref{sec_exam} we illustrate 
algorithm \ref{algorithm} on some concrete examples.   
In Section \ref{sec_pr} we describe the action of 
$\PIC^{(0,\dots,0)}(E)$ on $\VB_E^s(r,\dd)$ and 
morphisms between simple bundles, thus deduce the second part
of the Theorem \ref{theo_main}.


%

%% file: mainpart.tex
\section{General approach}
\subsection*{Category of triples}
\label{sec_tripl}
Let $\ko$ be an algebraically closed 
field\footnote{
Although the construction of triples 
and many classification results
 are valid for an arbitrary field, 
 the matrix problems that we obtain can be quite special
 and require different methods to deal with. 
 In order to get a uniform description for 
 all cases we assume from the beginning 
the ground field $\ko$ to be algebraically 
closed.
},
 $\Sch:=\Sch/_\ko$  the category 
of Noetherian  schemes over $\ko$
and
for any scheme $T\in \Sch$  by $\VB_T,$ $\TF_T$ and $\Coh_T$ 
we denote
the categories of vector bundles, torsion free coherent  
and coherent sheaves on $T$ respectively.

Let $\XX$ be a singular curve over $\ko.$
Fix the following notations:

\begin{itemize}
\item  $\pi: \widetilde{\XX} \lar \XX$  the normalization of $\XX$;
\item   $\oo:=\oo_{\XX}$  and $\on := \oo_{\xn}$ the structure sheaves of $X$ and $\xn$
respectively;
\item   $\jj= Ann_{\oo}(\pi_*\on/\oo)$ the conductor of $\oo$ in  $\pi_*\on$;
\item  $\i:S\emb X$ the subscheme of $X$ defined by the conductor $\jj$
and $\iw:\sn\emb \xn$ its scheme-theoretic pull-back to the normalization $\xn.$
\end{itemize}
Altogether they fit into a cartesian diagram:
\begin{align}
\label{diag_sch}
\begin{array}{c}
{
\xymatrix{
\widetilde{\ZZ} \ar[r]^{\iw} \ar[d]_{\tilde{\pi}}  & \widetilde\XX \ar[d]^{\pi}\\
\ZZ \ar[r]^{\i}  &  \XX.
}}
\end{array}
\end{align}

\noindent

%

\begin{remk}
\begin{itemize}  
\item[1.]
 In what follows we shall identify the structure sheaf $\oo_T$
 of an artinian scheme $T$ with the coordinate ring $\ko[T].$

\item [2.] 
The main property of the conductor is that 
for the ideal $\jn:= \kI_{\sn}$ in $\on$
we have $\jj=\pi_*\jn.$ 

\noindent
\item[3.]
Let $\ff\in\Coh_X$ and $\fn\in\Coh_{\xn}$ be coherent sheaves on 
$X$ and $\xn$ respectively.
With a little abuse of notation one can write:
 $\i^*\ff = \ff\otimes_{\oo}\oo_S=\ff/\jj\ff\in \Coh_S$ 
 and ${\iw^*\fn=\fn\otimes_{\on}\oo_{\sn}}={\fn/\jn\fn \in \Coh_{\sn}.}$ 
Since $S$ and $\sn$ are schemes of dimension zero,
$\i_*\i^*\ff$ and  $\iw_*\iw^*\fn$ are
 skyscraper sheaves on $X$ and $\xn$ respectively.
 


\end{itemize}
\end{remk}

The usual way to deal with  vector bundles on a singular curve is to lift them
 to the normalization, and work on a
smooth curve, see for example \cite{Seshadri, Bhosle92, Bhosle96}.
To describe the fibers of the map
$\VB_X \mo \VB_{\xn}$  
we recall the following construction:

\begin{defin}
\label{dfTriples }
The  \emph{category of triples}   ${\Tr}_{\XX}$  is defined as follows:
\begin{itemize}
\item
Its objects are  triples 
$(\fn,\kM, \mn)$, 
where $\fn\in \VB_{\xn},$ $\,\kM\in \VB_{S}$  and
$\mn : \pn^*\kM\mo \iw^*\fn$
is an isomorphism of $\oo_{\sn}$--modules.

\item  A morphism
 $
 \xymatrix
{
(\fn, \kM, \mn)
\ar[rr]^{(F, f)} &&
(\fn', {\kM'}, \mn')
}
 $
is given by a pair  $(F, f)$, where
$F:\fn\mo\fn'$
is a morphism in $\VB_{\xn}$ and
$f: {\kM} \mo {\kM'}$
is a morphism in $\Coh_{S}$, such that the following diagram
commutes in $\Coh_{\sn}:$
\begin{equation}
\label{diagMorphisms}
\begin{array}{c}
     \xymatrix
{
\pn^*\kM  \ar[rr]^{\mn} \ar[d]_{\pn^*f} & &
\iw^*\fn
\ar[d]^{\iw^* F}  \\
\pn^*\kM'  \ar[rr]^{\mn'}  & &
\iw^*\fn'.\\
}
\end{array}
\end{equation}
\end{itemize}
\end{defin}

\noindent
Raison d'\^etre  for the formalism of triples is the following theorem:

\begin{theorem}[\cite{drgr}]
\label{thm:triples}
The functor
$
\mathsf{\Psi}:\VB_{X}  \lar {\Tr}_{\XX}
$
taking a vector bundle $\ff$ to the triple
$(\fn, \kM, \mn)$, where
$\fn := \pi^*\ff$,
$\kM := \i^*\ff$ and $\mn$ is the canonical morphism 
$\mn:\tilde{\pi}^* \i^*\ff \mo \iw^*\pi^* \ff ,$
is
an equivalence of categories. 
\end{theorem}

\medskip
\noindent
Although the statement 
of Theorem \ref{thm:triples} holds for arbitrary reduced curves, 
the method based on it can be efficiently used mainly
for rational curves, since 
in this case
the description of vector 
bundles on the normalization is well understood. 

\subsection*{Vector bundles on a projective line}
According to the classical result 
known as Birkhoff-Grothendieck Theorem,  
a vector bundle $\fn$ 
on a projective line $\mP^1$ 
splits into a direct sum of line bundles:
\begin{equation}
\label{eqVectBundProjLine} 
\fn\cong \mathop\oplus_{n \in \mZ}\big(\oo_{\mP^1}(n)\big)^{r_n}.
\end{equation}
Let $(z_0:z_1)$ 
be homogeneous coordinates 
on $\mP^1.$ Then an endomorphism $F$ of
$\fn $ can be written in a
matrix form: 
\begin{equation}
  \label{eqFormF}
  F=
  \left(
  \begin{array} {cccccc}
    \ddots& 0&\dots& 0&0\\
    \dots& F^{nn}&\dots& 0&0\\
    & \vdots  &\ddots& \vdots&\vdots\\
    \dots& F^{mn}&\dots& F^{mm} &0\\
    & \vdots && \vdots&\ddots\\
  \end{array}
  \right),
\end{equation}
where $F^{mn}$ are blocks of sizes $r_{m} \times r_{n}$
with  coefficients in the vector space
\begin{align}
\label{eq_sq_form}
{\Hom_{\mP^1}\big(\oo_{\mP^1}(n), \oo_{\mP^1}(m)\big) \cong
\ko[z_{0}, z_{1}]_{m-n}},
\end{align}
since a morphism $\oo_{\mP^1}(n)\mo \oo_{\mP^1}(m)$
is determined by a homogeneous form
$Q(z_0,z_1)$ of degree $m-n.$
In particular, the matrix $F$
is lower-block-triangular and the diagonal $r_{n}\times r_{n}$ blocks
$F^{nn}$ are matrices over $\ko$. The morphism
$F$ is an isomorphism if and only if
all the diagonal blocks $F^{nn}$ are invertible.


\subsection*{Matrix problem $\MP_X$.}
To classify vector bundles on 
 a rational projective curve $X$ with the normalization
$\xn=\mathop\sqcup\limits_{k=1}^N L_k$ 
 one should describe
 iso-classes of objects in ${\Tr}_{\XX}.$
Note that two triples $(\fn,\kM,\mn)$ and $(\fn',\kM',\mn')$
are isomorphic only if 
$\fn\cong \fn'$ and $\kM\cong \kM'.$ 
By Birkhoff-Grothendieck theorem a bundle $\fn$ on $\xn$ can be given by 
a tuple of integers 
$\rr:= \{r(n,k)\},$ where $n\in \mZ,$ 
$1\leq k\leq N$ and $\sum_{n\in\mZ}r(n,k)=r$ for each $k.$
Let
$\MP_X:=\bigcup_{\rr}\MP_X(\rr)$ 
be the following Krull-Schmidt category:
an object of a stratum $\MP_X(\rr)$ is a matrix $\mn$
for which there exists a triple  $(\fn,\kM,\mn)\in \Tr_X$
and the vector bundle 
$\fn\in \VB_{\xn}$ splits into a direct sum of line bundles
with the tuple of multiplicities $\rr.$ 
For two objects $\mn$ and $\mn'$ with triples
$(\fn,\kM,\mn)$ and $(\fn',\kM',\mn')$ respectively,
a morphism from $\mn$ to $\mn'$ is a pair
$({\iw^* F}, \pn^*f)$ such that
$\iw^* F\ci\mn=\mn'\ci \pn^*f,$ 
where
$F\in \Hom_{\xn}(\fn,\fn')$ and  $f\in \Hom(\kM,\kM').$
The functor $ \mathsf{H}: \Tr_X\lar \MP_X$ 
is full and dense and
there is a natural projection
\begin{align}
  \label{eq_surj_hom}
\Hom_{\Tr_X} \big((\fn,\kM, \mn),(\fn,\kM, \mn')\big)\twoheadrightarrow
 \Hom_{\MP_X}(\mn,\mn').
 \end{align}

\noindent
\begin{defin}
Replacing the set of morphisms by the set of
invertible morphisms in $\MP_X(\rr)$
(also called \emph{matrix transformations})
we obtain some groupoid.
A \emph{matrix problem} is the problem of describing
orbits of indecomposable objects.
If it is possible, a solution consists in finding
 a \emph{canonical form} of $\mn.$ 
\end{defin}

The precise description of this procedure can be found in \cite{thesis}.
For convenience we choose $\ko$-bases of $\oo_S$  and $\oo_{\sn}$
and rewrite $\mn, \iw^*F$ and $\pn^*f$
as tuples of matrices over $\ko.$



\section{Matrix problem for cycles of projective lines.}
\label{sec_matr_pr_cycles}
Let $E$ be a cycle of $N$ projective lines.
The normalization $\widetilde{\EE}$ is a disjoint union 
of $N$ copies of $\mP^1.$
For example, for $N=3$  we have:
$$
\xy 
{\ar@{-}(-10,0);(-10,20)},
{\ar@{-}(0,0);(0,20)}, 
{\ar@{-}(10,0);(10,20)},
\POS(-12,0)*{_{L_1}},
\POS(-9.8,15)*{\bul};\POS(-13,15)*{_{\infty}}; 
 \POS(-9.8,5)*{\bul};\POS(-12,5)*{_{0}}; 
\POS(-2,0)*{_{L_2}},
\POS(0.2,15)*{\bul};\POS(-3,15)*{_{\infty}}; 
 \POS(0.2,5)*{\bul};\POS(-2,5)*{_{0}}; 
\POS(8,0)*{_{L_3}},
\POS(10,15)*{\bul};\POS(7,15)*{_{\infty}}; 
 \POS(10,5)*{\bul};\POS(8,5)*{_{0}}; 
 {\ar@{.}(-9.8,15);(0.2,5)},
 {\ar@{.}(0.2,15);(10,5)},
 {\ar@{.}(-9.8,5);(10,15)},
{\ar^{\pi}(13,10);(25,10)},
{\ar@{-}(29.5,7);(45.5,7)},
{\ar@{-}(35,16);(43.5,3.5)},
{\ar@{-}(39,16);(31.5,3.5)},
\POS(37,13)*{\bul}; \POS(34, 7)*{\bul}; \POS(41 ,7)*{\bul}
\POS(34,13)*{_{s_1}}; \POS(35,5)*{_{s_2}}; \POS(40, 5)*{_{s_3}};
%
\POS(36,0)*{_{E}}
\endxy
$$
Let $s_1,\dots s_N$ be the intersection points ordered in such a way
that $s_k$ and $s_{k+1}$ belong to the component $L_k$
 for $k=1,\dots , N-1$  and the points $s_N$ and $s_1$ lay on $L_N .$ 
On each component $L:=L_k$
choose the local coordinates such that the preimages 
of $s_k$ and $s_{k+1}$ on $L_k,$ for $k=1,\dots ,N-1,$
and $s_N$ and $s_1$ on $L_N$ have
coordinates $0 := (0 : 1)$ and $\infty := (1 : 0).$
Then
$$\oo_S =  \ko(s_1)\oplus\dots\oplus\ko(s_N) 
\hbox{ and } 
\oo_{\sn} = \bigoplus\limits_{k=1}^N(\ko(0_k)\oplus\ko(\infty_k)).$$

\noindent
To obtain the matrix problem $\MP_E$
 we fix:
\begin{itemize}
\item
 a splitting  
 $\fn \cong
\mathop\bigoplus\limits_{k=1}^N 
\Big(\mathop\oplus\limits_{n \in \mathbb{Z}} \oo_{L_k}(n)^{r(n,k)}\Big)$
with $\mathop\sum\limits_{n\in\mZ} r(n,k)=r$ for each component $k;$

\item an isomorphism $\kM \cong \oo_S^{r}=(\oplus_{k=1}^N\ko(s_k) )^r.$

   \item
The choice of coordinates
on each component $L$ of $\xn$
fixes two canonical sections $z_0$ and $z_1$ of $H^0\big(\oo_L(1)\big),$
and we  use the following trivializations 
\begin{align*}
\oo_L(n)\otimes \oo_{L\cap \sn}
&\stackrel{\sim}\lar \ko(0) \times  \ko(\infty)\\
\zeta\otimes 1 &\longmapsto (\zeta/z_1^n(0), \zeta/z_0^n(\infty)).
\end{align*}
 \end{itemize}
This isomorphism  only  depends on the choice 
of coordinates on $L\cong\PP^1$.
In such a way we equip the $\oo_{\sn}$--module $\iw^*\fn,$
where
$\iw^*\fn|_L = \fn|_L(0) \oplus \fn|_L(\infty),$ with a basis and get
isomorphisms
$\fn|_L(0) \cong \ko(0)^{r}$
 and 
$\fn|_L(\infty) \cong \ko(\infty)^{r}.$

\noindent
\subsection*{Matrix problem $\MP_E$ for Kodaira cycles $I_N$ }
With respect to  all  the choices  the maps 
$\mn,$   $\iw^*F$ and $\pn^*f$ can be written as matrices.
\begin{itemize}
  \item The gluing map
  $\mn :  \pn^*\kM \stackrel{\sim}\lar    \iw^*\fn$
consists of $2N$ invertible matrices over $\ko$
\begin{align}
\label{matr_cycles}
\mn&=\big(
\mu_{1}(0), \mu_{1}(\infty),\mu_{2}(0),\mu_{2}(\infty),\dots,
\mu_{N}(0),
\mu_{N}(\infty)
\big).
\end{align}

  \item
  If we have a morphism  $\oo_L(n) \mo \oo_L(m)$ 
given by a homogeneous form
$Q(z_0,z_1)$ of degree $m-n,$ then it induces a map
$\oo_L(n) \otimes \oo_{\sn}
\lar \oo_L(m) \otimes \oo_{\sn}$
given by ${(Q(0),Q(\infty)):=(Q(0:1) , Q(1:0))}.$
Hence, with  respect to the chosen trivializations
of $\oo_L(n)$ at $0$ and $\infty$ the map
\begin{align}
{\iw^* F}|_L
=\big(F_k(0), F_k(\infty)\big) 
: \ko^{r}(0) \oplus \ko^{r}(\infty) \lar \ko^{r}(0) \oplus \ko^{r}(\infty)
\end{align}
is given by a pair of lower
block triangular
matrices  $\big(F_k(0), F_k(\infty)\big)$ consisting of
blocks
 $F_k^{mn}(0), F_k^{mn}(\infty)
 \in \Mat_{\ko}(r(m,k) \times r(n,k)),$ for $ m>n$
 and
with common diagonal blocks
 $F_k^{nn}\in \Mat_{\ko}(r(n,k) \times r(n,k)).$ 
 The morphism $F$ is invertible, 
 if all the diagonal blocks $F_k^{nn}$ belong to
$\GL(\ko,r(n,k)).$ 

\item  The same holds for the induced map $\pn^*f=(f_1,\dots,f_N):$ 
if $(F,f)$ is invertible then 
${f_k\in\GL(\ko,r)}$ for each component $k .$

\item
The transformation rule  $\mn\to (\iw^*F)\ci \mn\ci (\pn^* f)^{-1}  $ 
can be rewritten for each component $k$ as 
$\mu_{k}(0)\to F_k(0)\mu_{k}(0)f_{k}^{-1}$ and
$\mu_{k}(\infty)\to F_k(\infty)\mu_{k+1}(\infty)f_{k}^{-1}$
assuming $f_{N+1}=f_1.$ 
For $N=3$ it can be sketched as follows:
\end{itemize}
$$
\xy ;/r0.17pc/:
{\ar@/_5pt/_{F_1(0)}(-2,18);(-2,2)},
{\ar@{.}_{f_1}(10,-1);(10,-39)},
%
{\ar@{-}(0,0);(0,20)}, 
{\ar@{-}(20,0);(20,20)},
{\ar@{-}(0,20);(20,20)}, 
{\ar@{-}(0,15);(20,15)},
{\ar@{-}(0,5);(20,5)},  
{\ar@{-}^{}(0,10);(20,10)},  
{\ar@{-}^{}(0,0);(20,0)}, 
{\ar@/^5pt/^{F_1(\infty)}(62,18);(62,2)},
{\ar@{-}(40,0);(40,20)}, 
{\ar@{-}(40,20);(60,20)}, {\ar@{-}(60,0);(60,20)},
{\ar@{-}(40,15);(60,15)},
{\ar@{-}(40,5);(60,5)},  
{\ar@{-}^{}(40,10);(60,10)},
{\ar@{.}(40,3);(20,3)}, {\ar@{.}(40,17);(20,17)},
 {\ar@{.}(40,13);(20,13)},  {\ar@{.}(40,7);(20,7)},
{\ar@{-}^{}(40,0);(60,0)}, 
{\ar@{.}^{f_2}(50,-1);(50,-9)},
{\ar@/_5pt/_{F_2(0)}(38,-12);(38,-28)},
%
{\ar@{-}(40,-30);(40,-10)}, 
{\ar@{-}(40,-10);(60,-10)}, 
{\ar@{-}(60,-30);(60,-10)},
{\ar@{-}^{}(40,-17);(60,-17)},{\ar@{-}^{}(40,-23);(60,-23)},
{\ar@{.}(80,-13);(60,-13)}, {\ar@{.}(80,-20);(60,-20)}, 
{\ar@{.}(80,-27);(60,-27)},
{\ar@{-}^{}(40,-30);(60,-30)}, 
{\ar@{-}(80,-30);(80,-10)}, {\ar@{-}(100,-30);(100,-10)},
{\ar@{-}(80,-10);(100,-10)}, 
{\ar@{-}^{}(80,-17);(100,-17)},{\ar@{-}^{}(80,-23);(100,-23)},
{\ar@{-}^{}(80,-30);(100,-30)}, 
{\ar@{.}^{f_3}(90,-31);(90,-39)},
{\ar@/^5pt/^{F_2(\infty)}(102,-12);(102,-28)},
{\ar@/_5pt/_{F_3(\infty)}(-2,-42);(-2,-58)},
{\ar@{-}(0,-60);(0,-40)}, {\ar@{-}(20,-60);(20,-40)},
{\ar@{-}(0,-40);(20,-40)}, 
{\ar@{.}(20,-42);(80,-42)},
{\ar@{.}(20,-49);(80,- 49)},
{\ar@{.}(20,-56);(80,-56)},
{\ar@{-}^{}(0,-53);(20,-53)},  
{\ar@{-}^{}(0,-45);(20,-45)},  
{\ar@{-}^{}(0,-60);(20,-60)}, 
{\ar@{-}(80,-60);(80,-40)}, {\ar@{-}(100,-60);(100,-40)},
{\ar@{-}(80,-40);(100,-40)}, 
{\ar@{-}^{}(80,-53);(100,-53)},  
{\ar@{-}^{}(80,-45);(100,- 45)},  
{\ar@{-}^{}(80,-60);(100,-60)}, 
{\ar@/^5pt/^{F_3(0)}(102,-42);(102,-58)},
\endxy
$$

\medskip
\medskip

\noindent
Since the matrices $F_k(0)$ and $F_k(\infty)$ have the block-triangular structure,
as described above,
thus the matrices
$\mu_k(0)$ and $\mu_k(\infty)$ 
split into horizontal blocks labeled by $n\in\mZ,$ 
as in the decomposition  of $\fn.$
Such blocks contain $r(n,k)$ rows, and 
can be transformed only together by
 $F_k^{nn}(0)=F_k^{nn}(\infty).$ 
 We call them \emph{conjugated} blocks and connect by 
 dotted lines. 

These types of matrix problems\index{matrix problem} 
are well-known in  representation theory.
They are called
 \emph{Gelfand problems}
or \emph{representations of bunches of chains} 
(see \cite{GP, Bon92}).
For an application of Gelfand problems to the classification 
of  torsion free
sheaves on cycles of projective lines we refer to \cite{drgr}
(see also\cite{surv}).

\section{Matrix problem for Kodaira fibers II, III and IV}
\label{sec_matr_pr_fibers}
In this section we formulate the matrix problem $\MP_E$ 
for the other curves from the Table \ref{list_1}. 
Let $E$ be a Kodaira fiber with $N$ ($N\leq 3$) components, 
 $s$ the unique singular point 
and $\pi: \widetilde{E} \mo E$ the normalization map. 
For example, for $N=3$ we have
$$
\xy 
{\ar@{-}(0,0);(0,20)}, {\ar@{-}(10,0);(10,20)}, {\ar@{-}(20,0);(20,20)},
\POS(-2,0)*{_{L_1}},\POS(8,0)*{_{L_2}}, \POS(18,0)*{_{L_3}},
\POS(0.2,10)*{\bul};\POS(-2,12)*{_{0}}; 
\POS(10.2,10)*{\bul}; \POS(8,12)*{_{0}};
\POS(20.2,10)*{\bul}; \POS(18,12)*{_{0}};
{\ar^{\pi}(25,10);(37,10)},
{\ar@{-}(50,0);(50,20)};
{\ar@{-}(41,4);(59,16)};
{\ar@{-}(59,4);(41,16)};
\POS(53,10)*{_{s}};\POS(50.3,10)*{\bul};
\POS(57,0)*{_{E}}
\endxy
$$
Note that $\widetilde{E}$
consists of a disjoint union of $N$ projective lines.
On each component  $L_k$
choose coordinates $(z_0: z_1)$ 
such that the preimage of
the singular point $s = (0:0:1)$ on $L_k$  is $0:=(0:1)$. 
Let $U_k = \{(z_0:z_1)| z_1 \ne 0\}$
be  affine neighborhoods 
of $0$ on $L_k$
with local coordinates  $t_k := z_0/z_1$ for $k=1,\dots , N ;$
and let $U$ be
the union $\bigcup_{k=1}^N\pi(U_k).$
Calculate 
the normalization map 
$\oo\emb\pi_*\on=
\pi_*\big(\mathop\oplus\limits_{k=1}^N \oo_{L_k}\big),$
the conductor $\jj$ and 
the structure sheaves $\oo_S,$ and $\oo_{\sn}$ for each 
Kodaira fiber:

\begin{itemize}
  \item[II.] 
  Let $\EE$ be a cuspidal cubic curve in $\mP^2$ 
given by the equation  $x^3-y^2z=0.$ Then locally 
the normalization map is
  $\ko[U]=\ko[t^2,t^3]\emb \ko[t].$
  Since on $\pi(U)$ the conductor is 
  $\jj=\langle t^2 , t^3 \rangle,$ we have 
  $\oo_S\cong\ko(s) \hbox{ and } 
\oo_{\sn}\cong \big(\ko[\varepsilon]/\varepsilon^2\big)(0).$

  \item[III.]Let $E$ be a tacnode curve given
by the equation  $y(zy-x^2)=0.$
 Then the normalization map is
  $
\ko[U]\emb \ko[t_1]\oplus \ko[t_2]
$ taking $1 \mapsto (1,1),$
$x \mapsto (t_1,t_2),$ and
$y \mapsto (0,t_2^2).$
On $\pi(U)$ for the conductor we have
$\jj=\langle(t_1^2,0), (0,t_2^2) \rangle.$
In other words, 
the ideal sheaf of the scheme-theoretic preimage of $s$ is
$\jn = \Big(\kI^2_{L_1,0}, \kI^2_{L_2,0}\Big),$ 
where $\kI_{L_k,0}$ denotes 
the ideal sheaf of the point $0$ on the component $L_k.$
Hence, $\oo_{\sn}\cong \on/\jn=\oo_{L_1}/\kI^2_{L_1,0}\oplus \oo_{L_2}/\kI^2_{L_2,0}.$
Altogether we get
$
\oo_S\cong\big(\ko[ \e]/\e^2\big)(s),
$
and 
$
\oo_{\sn} \cong\big(\ko[\e_1]/\e_1^2\big)(0)\oplus
\big(\ko[\e_2]/\e_2^2\big)(0)
$
and the induced map $\oo_S\emb \pn_*\oo_{\sn}$ 
takes $\e$ to $(\e_1,\e_2).$

  \item[IV.]Let $E$ be a curve  consisting of 
three concurrent projective lines in $\mP^2,$ 
given by the equation $xy(x-y)=0.$
Then the normalization map is
$
\ko[U]\emb \ko[t_1]\oplus \ko[t_2]\oplus \ko[t_3],
$
sending
$1\mapsto (1,1,1),$
$x\mapsto (t_1,t_2,0),$ and
$y\mapsto (t_1,0,t_3).$
Since  $\jj(U)=\langle x^2,y^2,xy\rangle,$ we have
$\oo_S =\ko[  x,y ]/\langle x^2,y^2,xy\rangle.$
%
Note that the ideal sheaf $\jn:=\pi^*\jj$ is
locally generated by
$
(t_1^2,0,0),$ $(0,t_2^2,0)$ and $(0,0,t_3^2)$
i.e.
$\jn = \Big(\kI^2_{L_1,0}, \kI^2_{L_2,0},\kI^2_{L_3,0}\Big),$ 
where $\kI_{L_k,0}$ is as above.
Hence,
$\oo_{\sn}\cong \mathop\oplus\limits_{k=1}^3
\oo_{L_k}/\kI^2_{L_k,0}.$

\end{itemize}

\subsection*{Matrix problems $\MP_E$ for Kodaira fibers II, III and IV}
 For a triple
$(\fn, \kM, \mn)$ 
 we fix:
\begin{itemize}
\item
 a splitting  
 $\fn \cong
\mathop\bigoplus\limits_{k=1}^N 
\Big(\mathop\oplus\limits_{n \in \mathbb{Z}} \oo_{L_k}(n)^{r(n,k)}\Big)$
with $\mathop\sum\limits_{n\in\mZ} r(n,k)
=r;$

\item an isomorphism $\kM \cong \oo_S^{r};$

\item  
for each component $L:=L_k$
we take the trivializations 
\begin{align*}
\oo_L(n) \otimes  \oo_L/\kI^2_{L,0}& \lar\ko[\e_k]/ \e_k^2,\\
\zeta\otimes 1 &\longmapsto pr(\frac{\zeta}{z_1^n})
\end{align*}
for a local section $\zeta$
of $\oo_{L_k}(n)$  on the open set $U_k,$ 
where the projection 
$$pr: \ko[U_k] \lar \ko[\e_k]/\e_k^2$$
is the map induced by  
$
\ko[t_k] \lar \ko[\e_k]/\e_k^2,$
mapping
$t_k \mapsto \e_k.$
\end{itemize}

\noindent
With respect to  all  these choices we have:
\begin{itemize}
 \item
The map $\mn$ can be written as a combination 
of $2N$  $r\times r$-matrices over $\ko$:  
\begin{equation}
\label{eq form of mu_trypr}
\mn = (\mu_1,\dots,\mu_N)= 
\Big(  \mu_1(0) + \e_1 \cdot \mu_{\e_1}(0),~~
   \dots,~~
   \mu_N(0) + \e_N \cdot \mu_{\e_N}(0)
   \Big).
\end{equation}
\noindent
The morphism  $\mn$ is  invertible if and only if
 all $\mu_k(0)\in\GL(\ko,r).$ 
 
\item 
If on a component $L=L_k$ we have a morphism  $\oo_L(n) \mo \oo_L(m)$ 
given by a homogeneous form
$Q(z_0,z_1)$ of degree $m-n,$ then the induced map
${\oo_L(n) \otimes \oo_{\sn}
\lar \oo_L(m) \otimes \oo_{\sn}}$
is given by the map
$$pr(Q(z_0, z_1)/z_1^{m-n}) = Q(0:1) + \varepsilon_k {\textstyle\frac{d Q}{d z_0}(0:1).}$$
Hence,
for a morphism $(F,f):(\fn, \kM, \mn)\lar (\fn', \kM', \mn')$ 
the induced map $\iw^*F: \iw^*\fn\lar \iw^*\fn'$ 
 is
$${\iw^* F}|_{L}=F_k(0)  + \varepsilon_k \textstyle{\frac{d F_k}{d z_0}(0)}
\in \Mat(\ko[\varepsilon_k]/\varepsilon_k^2,r),
$$
where, 
as usual, $F_k(0)$ denotes $F_k(0:1).$

\item  The morphism
$\pn^*f$ consists of $N$ copies of the matrix $f,$
where 
\begin{itemize}
  \item $f\in\Mat(\ko,r\times r) $ for the cuspidal cubic;
  \item $f=f(0)+f_{\e}(0)\in\Mat(\ko[\e]/\e^2,r\times r),$ 
for $\e=(\e_1,\e_2)$  for the tacnode curve (Kodaira fiber III);
  \item 
$f=f(0)+x \cdot f_x(0)+y \cdot f_y(0)
\in \Mat\Big({\ko[  x,y ]/\langle x^2,y^2,xy\rangle},r\times r\Big)$
for the three lines through a point in a plane (Kodaira fiber IV).
\end{itemize}
\end{itemize}
\noindent
A morphism $(F,f)$ is an automorphism if and only if 
all ${F_k}(0)$ for $k\in\{1,\dots, N \}$ and  $f(0)$ are  invertible 
$r\times r$ matrices over $\ko.$
For example, for the Kodaira fiber IV we get the following
matrix problem. 
There are six $r\times r$ matrices $\mu_1(0),$  $\mu_{\e_1}(0),$ 
$\mu_2(0),$  $\mu_{\e_2}(0)$ and $\mu_3(0),$  $\mu_{\e_3}(0),$
where  all $\mu_k(0)$
are invertible.
 The pairs
$\mu_k(0), \mu_{\e_k}(0)$ 
are simultaneously 
 divided into horizontal blocks 
 with respect to the splitting of $\fn|_{L_k}.$
$$
\xy ;/r0.17pc/:
{\ar@/_5pt/_{F_1(0)}(-2,58);(-2,42)},
{\ar@/_7pt/_{f_y(0)}(10,38);(50,38)},
{\ar@/^7pt/_{f_x(0)}(10,62);(50,62)},
{\ar@{-}(0,40);(0,60)}, {\ar@{-}(20,40);(20,60)},
{\ar@{-}(0,60);(20,60)}, 
{\ar@{-}^{}(0,54);(20,54)}, 
{\ar@{-}^{}(0,51);(20,51)},
{\ar@{-}^{}(0,48);(20,48)},
{\ar@{-}^{}(0,45);(20,45)},  
{\ar@{-}^{}(0,40);(20,40)}, 
{\ar@/^5pt/^{F_1(0)}(62,58);(62,42)},
{\ar@{-}(40,40);(40,60)}, 
{\ar@{-}(40,60);(60,60)}, {\ar@{-}(60,40);(60,60)},
{\ar@{-}^{}(40,54);(60,54)}, 
{\ar@{-}^{}(40,51);(60,51)},{\ar@{.}(40,52.5);(20,52.5)},
{\ar@{-}^{}(40,48);(60,48)},{\ar@{.}(40,49.5);(20,49.5)},
{\ar@{-}^{}(40,45);(60,45)}, {\ar@{.}(40,46.5);(20,46.5)}, 
{\ar@{-}^{}(40,40);(60,40)}, 
\POS(22,57);\POS(38,45)*\dir{>}**\crv{(35,56)&(25,46)};
\POS(29,41)*{\scriptstyle{\frac{d F_1}{d z_0}(0)}};
{\ar@/_5pt/_{F_2(0)}(-2,18);(-2,2)},
{\ar@/^7pt/_{f_x(0)}(10,22);(50,22)},
{\ar@{-}(0,0);(0,20)}, {\ar@{-}(20,0);(20,20)},
{\ar@{-}(0,20);(20,20)}, 
{\ar@{-}^{}(0,16);(20,16)}, 
{\ar@{-}^{}(0,13);(20,13)},
{\ar@{-}^{}(0,10);(20,10)},
{\ar@{-}^{}(0,7);(20,7)},  
{\ar@{-}^{}(0,0);(20,0)}, 
{\ar@{.}_{f(0)}(8,39);(8,21)},
{\ar@{.}^{f(0)}(52,39);(52,21)},
%
%
{\ar@/^5pt/^{F_2(0)}(62,18);(62,2)},
{\ar@{-}(40,0);(40,20)}, 
{\ar@{-}(40,20);(60,20)}, {\ar@{-}(60,0);(60,20)},
{\ar@{-}^{}(40,16);(60,16)}, 
{\ar@{-}^{}(40,13);(60,13)},{\ar@{.}(40,14.5);(20,14.5)},
{\ar@{-}^{}(40,10);(60,10)},{\ar@{.}(40,11.5);(20,11.5)},
{\ar@{-}^{}(40,7);(60,7)}, {\ar@{.}(40,8.5);(20,8.5)}, 
{\ar@{-}^{}(40,0);(60,0)}, 
%
\POS(22,17);\POS(38,5)*\dir{>}**\crv{(35,16)&(25,6)};
\POS(30,1)*{\scriptstyle{\frac{d F_2}{d z_0}(0)}};
{\ar@{.}^{f(0)}(8,-9);(8,-1)},
{\ar@{.}_{f(0)}(52,-8);(52,-1)},
{\ar@/_5pt/_{F_3(0)}(-2,-12);(-2,-28)},
{\ar@/_7pt/_{f_y(0)}(10,-32);(50,-32)},
{\ar@{-}(0,-30);(0,-10)}, {\ar@{-}(20,-30);(20,-10)},
{\ar@{-}(0,-10);(20,-10)}, 
{\ar@{-}^{}(0,-16);(20,-16)}, 
{\ar@{-}^{}(0,-19);(20,-19)},
{\ar@{-}^{}(0,-22);(20,-22)},
{\ar@{-}^{}(0,-25);(20,-25)},  
{\ar@{-}^{}(0,-30);(20,-30)}, 
{\ar@/^5pt/^{F_3(0)}(62,-12);(62,-28)},
{\ar@{-}(40,-30);(40,-10)}, 
{\ar@{-}(40,-10);(60,-10)}, {\ar@{-}(60,-30);(60,-10)},
{\ar@{-}^{}(40,-16);(60,-16)}, 
{\ar@{-}^{}(40,-19);(60,-19)},{\ar@{.}(40,-17.5);(20,-17.5)},
{\ar@{-}^{}(40,-22);(60,-22)},{\ar@{.}(40,-20.5);(20,-20.5)},
{\ar@{-}^{}(40,-25);(60,-25)}, {\ar@{.}(40,-23.5);(20,-23.5)}, 
{\ar@{-}^{}(40,-30);(60,-30)}, 
%
\POS(22,-13);\POS(38,-25)*\dir{>}**\crv{(35,-14)&(25,-24)};
\POS(30,-29)*{\scriptstyle{\frac{d F_3}{d z_0}(0)}};
\endxy
$$

\noindent
Note that $f_x(0)$ does not act on $\mu_{\e_3}$ and 
$f_y(0)$ does not act on $\mu_{\e_2},$
since as explained above, the normalization map 
 $\oo_S\emb \pn_*\oo_{\sn}$ 
sends
$x\mapsto (\e_1,\e_2,0),$ and
$y\mapsto (\e_1,0,\e_3).$

If we restrict this problem on the first two components
and assume $f_y(0)=0$  and $f_{\e}:=f_x(0)$ we obtain the matrix problem
for a tacnode curve.
If we restrict it to the first component  with 
$f_y(0)=f_x(0)=0$ we get the matrix problem for the cuspidal cubic curve.
Each of this problems is wild even for two horizontal blocks,
see \cite[Section 1]{Dro92} or \cite{boddr2}. 
However, the simplicity 
condition of a triple $(\fn,\kM,\mn )$ imposes
some additional restrictions making the problem tame.

\section{Simplicity condition}
\label{sec_simpl_cond}

  A vector bundle on a curve $X$ is called \emph{simple}
  if it admits no  endomorphisms but homotheties, i.e.
  $\End_{X}(\ff)=\ko$  and
the subcategory 
of simple vector bundles is denoted by $\VB_{X}^s.$
This notion can be obviously translated to the language of triples.
In terms of matrix problems:
an object $\mn$ of $\MP_X$ 
is called a \emph{brick}
if $\End_{\MP_X}(\mn)=\ko$.
The full subcategory of bricks is denoted by $\MP^s_X$
and $\MP^s_X(\rr)$ if the dimension vector $\rr$ is fixed.
Note that
a nonscalar morphism $(F,f)$ can have a scalar
restriction $(\iw^*F,\pn^*f).$ 

\begin{lemma}
\label{main_simple} 
Let $X$ be a rational singular curve and 
$(\fn,\kM,\mn)\in\Tr_X$ be a triple.
Then the map 
$ \End_{\Tr_X}(\fn,\kM,\mn)\mo \End_{\MP_X}(\mn)$
is bijective if and only if 
for all the components $L$ of $\xn$
 and for all summands $\oo_L(n)\oplus\oo_L(m)$ of $\fn|_L$
the canonical maps
 $\Hom(\oo_L(n),\oo_L(m))\mo\ko[\sn\cap L],$
 taking
$
 Q\to \iw^*Q \nonumber,
 $
 are bijective. 
\end{lemma}

\noindent
This obvious lemma implies certain nice properties
for a matrix problem under the simplicity condition.
For the curves under consideration, we have the following:

\begin{lemma}
\label{lemma_simple_form}
Let $E$ be a Kodaira fiber I$_N,$ (for $N\in\mN$)  II, III or IV, 
and let 
 $(\fn,\kM,\mn)\in\Tr_X$ be a simple triple, i.e.  $\End_{\Tr_X}(\fn,\kM,\mn)=\ko. $
 Then for each component $L:=L_k$ 
($1\leq k \leq N$) 
 \begin{align}
 \label{eq_form_fn}
 \fn|_L=\big(\oo_{L}(n_k)\big)^{r-\bd_{k}}
 \oplus\big(\oo_{L}(n_k+1) \big)^{\bd_{k}}
 \end{align}
for some 
$n_k\in \mZ$ and $1\leq \bd_k\leq r.$
\end{lemma}

\begin{proof}
 Assume that $\pi^*\ff|_L$ contains a summand 
$\oo_L(n)\oplus\oo_L(m)$ with $m\geq n+2.$
Let $(z_0:z_1)$ be the
local coordinates  as in Sections
\ref{sec_matr_pr_cycles} and \ref{sec_matr_pr_fibers}.
Since the degree $m-n\geq 2$
there exists a nonzero homogeneous form 
$Q\in {\Hom_{L}(\oo_{L}(n), 
\oo_{L}(m)) \cong \ko[z_{0}, z_{1}]_{m-n}}$
such that $\iw^*Q=0.$
Indeed, if $E$ is a Kodaira cycle then $\iw^*Q= (Q(0),Q(\infty))$
and if $E$ is a Kodaira fiber of type II, III or IV 
then the restriction of $\jn$ to the component $L$ is 
  $\kI_{L,0}^2\subset \oo_{L,0}$ and thus
$\iw^*Q =Q(0)+\frac{ \d Q}{\d z_0}(0).$
In both cases  
the map $Q\to \iw^*Q$ is not injective and
we get a contradiction to the condition of Lemma \ref{main_simple}.
\end{proof}

\begin{remk}
  \label{remk_shifts}
Note that the twists $n_k$ 
do not affect the matrix problem. 
Hence we can assume that the blocks have weights 
0 and 1 for each component $L_k$ and
replace  the multidegree $\dd$ by
$(\bd_1,\dots,\bd_N)$ 
and the degree
$d$ by $\bd:=\bd_1+\dots +\bd_N,$
where
$\bd_k=d_k\bmod r.$ 
Having the
twists $n_k$ we can recover the multidegree of $\dd$
 by the rule $d_k =r\cdot n_k + \bd_k.$
\end{remk}


\section{Primary reduction.}
\label{sec_prim_red}
Applying condition (\ref{eq_form_fn}) to the
matrix problem $\MP_E$ we obtain 
that 
each matrix consists of 
at most two horizontal blocks. 
Despite of this simplification the problem remains quite cumbersome.
However, it can be reduced to a 
\emph{partial canonical form,} such that all 
its matrices but one consist of identity and zero blocks.
We denote by $\M$ the remaining nonreduced matrix and formulate 
for it a new matrix problem.
It seems reasonable to introduce some simplified system of notations.
\begin{itemize}
  \item Let $\one$ denotes  the identity blocks,  $0$  the zero blocks,
    \item use the star $*$  to denote nonreduced blocks and
small Latin letters for a finer specification. 
\end{itemize}
The matrix $\M$ is divided into blocks,
the set of column-blocks coincides 
with the set of row-blocks and is denoted by $I=\{1,2,\dots , |I|\} .$ 
Then $\sss= (s_1,\dots,s_{|I|})\in \mN^{I}$ 
is the dimension vector of $M.$

\subsection{Nodal cubic curve}
\label{subsec_prim_node}
According to Section \ref{sec_matr_pr_cycles}
the matrix problem $\MP_E$
for the nodal cubic curve $E$ with two blocks is as follows:
$$
\xy ;/r0.15pc/:
{\ar@/_5pt/_{F(0)}(-2,18);(-2,2)},
{\ar@{-}(0,0);(0,20)}, {\ar@{-}(20,0);(20,20)},
{\ar@{-}(0,20);(20,20)}, 
{\ar@{-}^{}(0,10);(20,10)},
{\ar@{-}^{}(0,0);(20,0)}, 
%
{\ar@/^5pt/^{F(\infty)}(62,18);(62,2)},
{\ar@{-}(40,0);(40,20)}, 
{\ar@{-}(40,20);(60,20)}, {\ar@{-}(60,0);(60,20)},
{\ar@{-}^{}(40,10);(60,10)},{\ar@{.}(40,15);(20,15)},
{\ar@{.}(40,5);(20,5)}, 
{\ar@{-}^{}(40,0);(60,0)}, 
{\ar@/^10pt/@{.}_{f}(10,20);(50,20)}, 
\POS(10,-5)*{_{\mu(0)}};
\POS(50,-5)*{_{\mu(\infty)}};
\endxy
$$
Since the normalization consists of a unique component $L$
we skip the indices by $F,f$ and $\mu.$ 
As it was mentioned above both matrices
$\mu(0)$ and $\mu(\infty)$ are invertible. 
We reduce one of them,
say $\mu(0),$ to the identity form:
\begin{align*}
\mu(0)=
\begin{array}{|@{}     c       @{}|@{}     c       @{}|}
\hline
{\mk{\one}}  &    {\mk{0}}    \\
\hline
 \mk{0}        &        {\mk{\one}} \\
\hline
\end{array}
~~~\;\; \hbox{ and }~~~~
\M:=\mu(\infty)=
\begin{array}{ |@{}     c       @{}|@{}     c       @{}|}
\hline
 \mk{a_{1}}  &    \mk{b}   \\
\hline
 \mk{c}      &        {\mk{a_{2}}} \\
\hline
\end{array}.
\end{align*}
To preserve $\mu(0)$ unchanged we assume $f=F(0).$
Reformulate the problem for the matrix \linebreak
${\M:=\mu(\infty).}$
The transformation rule is
$ \M \to S \M  (S')^{-1},$
where 
$$
(S,S'):= (F(\infty), F(0))=
\left(
\begin{array}{
|@{}     c       @{}|@{}     c       @{}|}
\hline
 {\mk{w_1}}  &    \mk{0}        \\
\hline
 {\mk{u}}    &        {\mk{w_2}}  \\
\hline
\end{array}
\,
\begin{array}{ |@{}    c       @{}|@{}     c       @{}|}
\hline
 {\mk{w_1}}  &    \mk{0}        \\
\hline
{\mk{v}}      &        {\mk{w_2}}  \\
\hline
\end{array}
\right).
$$
Note that the sizes of blocks are determined by 
rank and degree: $(s_1,s_2)=(r-\bd,\bd),$
where $\bd:= d\bmod r .$
\subsection{Cuspidal cubic curve}
\label{subsec_prim_cusp}
Recall the problem $\MP_E$ on two blocks for the cuspidal curve: 
$$
\xy ;/r0.15pc/:
{\ar@/_5pt/_{F(0)}(-2,18);(-2,2)},
{\ar@{-}(0,0);(0,20)}, {\ar@{-}(20,0);(20,20)},
{\ar@{-}(0,20);(20,20)}, 
{\ar@{-}^{}(0,10);(20,10)},
{\ar@{-}^{}(0,0);(20,0)}, 
%
{\ar@/^5pt/^{F(0)}(62,18);(62,2)},
{\ar@{-}(40,0);(40,20)}, 
{\ar@{-}(40,20);(60,20)}, {\ar@{-}(60,0);(60,20)},
{\ar@{-}^{}(40,10);(60,10)},{\ar@{.}(40,15);(20,15)},
{\ar@{.}(40,7);(20,7)}, 
{\ar@{-}^{}(40,0);(60,0)}, 
{\ar@/^10pt/@{.}_{f}(10,20);(50,20)}, 
\POS(22,17);\POS(38,5)*\dir{>}**\crv{(35,16)&(25,6)};
\POS(29,0)*{\scriptstyle{\frac{d F}{d z_0}(0)}};
\POS(10,-5)*{_{\mu(0)}};
\POS(50,-5)*{_{\mu_{\e}(0)}};
\endxy
$$
As in the case of a nodal curve 
we skip the indices by $F,f$ and $\mu.$ 
The matrix $\mu(0)$ can be reduced to the identity form. 
To preserve this form unchanged we assume $F(0) = f.$ 
Moreover,
 using transformations $\frac{d F}{d z_0}(0)$
 we can make zero on the left lower block of $\mu_\varepsilon(0)$:
\begin{align*}
\mu(0)=
\begin{array}{|@{}     c       @{}|@{}     c       @{}|}
\hline
{\mk{\one}}  &    {\mk{0}}    \\
\hline
 \mk{0}        &        {\mk{\one}} \\
\hline
\end{array}
~~~\;\; \hbox{ and }~~~~
\M:=\mu_\varepsilon(0)=
\begin{array}{|@{}     c       @{}|@{}     c       @{}|}
\hline
 {\mk{a_{1}}}  &    {\mk{b}}      \\
\hline
 \mk{0}        &        {\mk{a_{2}}} \\
\hline
\end{array}.
\end{align*}
We obtain a new matrix problem which reads:
$ \M \to S \M  S^{-1}\bmod 
\left(
\begin{smallmatrix}
  0&0\\
  \times& 0
\end{smallmatrix}
\right),
$
where the matrix $S$ 
inherits the same lower-block-triangular structure as $F(0):$
$$S:=F(0)=f=
\begin{array}{
|@{}     c       @{}|@{}     c       @{}|}
\hline
 {\mk{w_1}}  &    \mk{0}        \\
\hline
 {\mk{u}}    &        {\mk{w_2}}  \\
\hline
\end{array}
.$$ 
As in the previous case
the sizes of blocks are determined by 
rank and degree: 
$${(s_1,s_2)=(r-\bd,\bd),} \hbox{ where }
 \bd:= d\bmod r .$$

\subsection{Cycle of two lines}
\label{subsec_prim_cycle_2}
According to Section \ref{sec_matr_pr_cycles}
the original matrix problem $\MP_E$ for a cycle of two lines
with two blocks on each component is
$$
\xy ;/r0.15pc/:
{\ar@/_5pt/_{F_1(0)}(-2,18);(-2,2)},
{\ar@{-}(0,0);(0,20)}, {\ar@{-}(20,0);(20,20)},
{\ar@{-}(0,20);(20,20)}, 
{\ar@{-}^{}(0,10);(20,10)},
{\ar@{-}^{}(0,0);(20,0)}, 
%
{\ar@/^5pt/^{F_1(\infty)}(62,18);(62,2)},
{\ar@{-}(40,0);(40,20)}, 
{\ar@{-}(40,20);(60,20)}, {\ar@{-}(60,0);(60,20)},
{\ar@{.}(40,15);(20,15)},
{\ar@{-}^{}(40,10);(60,10)},{\ar@{.}(40,5);(20,5)},
{\ar@{-}^{}(40,0);(60,0)}, 
{\ar@{.}^{f_1}(10,-1);(10,-9)},
{\ar@{.}^{f_2}(50,-1);(50,-9)},
{\ar@/_5pt/_{F_2(\infty)}(-2,-12);(-2,-28)},
{\ar@{-}(0,-30);(0,-10)}, {\ar@{-}(20,-30);(20,-10)},
{\ar@{-}(0,-10);(20,-10)}, 
{\ar@{-}^{}(0,-20);(20,-20)},  
{\ar@{-}^{}(0,-30);(20,-30)}, 
{\ar@/^5pt/^{F_2(0)}(62,-12);(62,-28)},
{\ar@{-}(40,-30);(40,-10)}, 
{\ar@{-}(40,-10);(60,-10)}, {\ar@{-}(60,-30);(60,-10)},
{\ar@{-}^{}(40,-20);(60,-20)},
{\ar@{.}(20,-15);(40,-15)},
{\ar@{.}(20,-25);(40,-25)},
{\ar@{-}^{}(40,-30);(60,-30)}, 
\endxy
$$

\noindent
All four matrices $(\mu_1(0),\mu_1(\infty),\mu_2(\infty)\mu_2(0))$
are invertible.
Two diagonal matrices, say $\mu_1(0)$ and $\mu_2(0),$
can be reduced to the identity form. Then one of the others,
say $\mu_2(\infty),$ can be reduced to the form:
\vspace{-0.5cm}
\begin{align}
\label{id_form}
\begin{array}{@{}c@{}l@{}}
\begin{array}{  @{}c @{}c @{}c @{}c @{}}
\mk{_1}&\mk{_2}&\mk{_3}&\mk{_4}
\end{array}
\\
\begin{array}{ @{}|@{}c   @{} c @{}|@{} c   @{}   c   @{}|@{}}
\hline
{\mk{\one}}      &\mk{0}          & \mk{0}&     \mk{0}    \\
{\mk{0}}      &\mk{0}          & \mk{\one}&     \mk{0}    \\
\hline
{\mk{0}}      &\mk{\one}          & \mk{0}&     \mk{0}    \\
{\mk{0}}      &\mk{0}          & \mk{0}&     \mk{\one}    \\
\hline
\end{array}
&
\begin{array}{@{}l@{}}
\tmk{\scriptstyle{1}}\\\tmk{\scriptstyle{3}}\\ 
\tmk{\scriptstyle{2}}\\\tmk{\scriptstyle{4}}
\end{array}
\\\phantom{a}
\end{array}
\end{align}
\noindent
Transformations $(F,f)$ preserving the reduced matrices
$\mu_1(0),$ $\mu_2(0)$ and  $\mu_2(\infty)$ unchanged 
satisfy the equations  
\begin{equation}
\label{eq_for_F}
f_1=F_1(0),\; f_2=F_2(0)
\,\hbox{ and }\, F_2(\infty)\mu_2(\infty)=\mu_2(\infty) f_1 .
\end{equation}
This implies the following triangular structures for  $F_1(0)$ and $F_2(\infty):$
\begin{align}
  \label{eq_F_1(0)}
F_1(0)
=
\begin{array}{c}
\begin{array}{ |@{}c   @{} c @{}|@{} c   @{}   c   @{}|}
\hline
{\mk{w_1}}      &\mk{0}          & \mk{0}&     \mk{0}    \\
{\mk{x_{21}}}      &\mk{w_2}          & \mk{0}&     \mk{0}    \\
\hline
{\mk{x_{31}}}      &\mk{0}          & \mk{w_3}&     \mk{0}    \\
{\mk{x_{41}}}      &\mk{x_{42}}          & \mk{x_{43}}&     \mk{w_4}    \\
\hline
\end{array}
\end{array}
\hbox{ and }
F_2(\infty)
=
\begin{array}{c}
\begin{array}
{ |@{}c   @{} c @{}|@{} c   @{}   c   @{}|}
\hline
\mk{w_1}      &{\mk{0}}           & \mk{0}   &     \mk{0}    \\
\mk{x_{31}}   &{\mk{w_3}}         & \mk{0}   &     \mk{0}    \\
\hline
\mk{x_{21}}   &{\mk{0}}           & \mk{w_2} &     \mk{0}      \\
\mk{x_{41}}   &{\mk{x_{43}}}      & \mk{x_{42}}  & \mk{w_4}  \\
\hline
\end{array}
\end{array}.
\end{align}
Since the diagonal blocks of $F_k(0)$ and $F_k(\infty)$
coincide (for $k=1,2$), we also have:  
$$F_1(\infty)
=
\begin{array}{c}
\begin{array}{ |@{}c   @{} c @{}|@{} c   @{}   c   @{}|}
\hline
{\mk{w_1}}      &\mk{0}          & \mk{0}&     \mk{0}    \\
{\mk{x_{21}}}      &\mk{w_2}          & \mk{0}&     \mk{0}    \\
\hline
{\mk{y_{31}}}      &\mk{y_{32}}          & \mk{w_3}&     \mk{0}    \\
{\mk{y_{41}}}      &\mk{y_{42}}          & \mk{x_{43}}&     \mk{w_4}    \\
\hline
\end{array}
\end{array}
\hbox{ and }
F_2(0)
=
\begin{array}{c}
\begin{array}
{ |@{}c   @{} c @{}|@{} c   @{}   c   @{}|}
\hline
\mk{w_1}    &{\mk{0}}             & \mk{0}&     \mk{0}    \\
\mk{x_{31}} &{\mk{w_3}}            & \mk{0}&     \mk{0}    \\
\hline
\mk{z_{21}} &{\mk{z_{23}}}              &     \mk{w_2}     & \mk{0} \\
\mk{z_{41}} &{\mk{z_{43}}}              &     \mk{x_{42}}   & \mk{w_4}  \\
\hline
\end{array}
\end{array}.
$$
\subsection*{Reduced matrix problem}
Thus we obtain a new problem for
the matrix $M:= \mu_1(\infty)$ with the transformations 
${\M \to S \M  (S')^{-1},}$
where $(S,S'):=(F_1(\infty),F_2(0)).$
Note that
if the sizes of blocks 1 and 4 are both nonzero then taking 
a nonzero entry $x_{41}$ of the matrices $F_2(\infty)$ and $F_1(0)$ 
we obtain a nonscalar endomorphism.
Hence, there are no sincere bricks and the maximal
tuples of blocks are
$I=(1,2,3)$ and its dual $I=(2,3,4).$
The dimension vector 
$\sss=(s_i)_{i\in I}$ 
and the matrix problem are
determined by  $r$ and 
$(\bd_1,\bd_2),$ where $\bd_k= d_k\bmod r$ and $\bd=\bd_1+\bd_2,$ 
as follows:
\begin{table}[ht]
$$
\begin{array}{|c|c|c |c| c |}
\hline
&
\hbox{condition}& \hbox{ set }\; I& \hbox{dimension vector }\; \sss
& \hbox{state}\\
\hline
\hline
1.&
\rule{0pt}{12pt} 
r\geq \bd & (1,2,3)& (r-\bd,\bd_2,\bd_1) 
 &A^+ \\
\hline
1'.&
\rule{0pt}{12pt} 
r< \bd & (2,3,4)& (r-\bd_1,r-\bd_2,\bd-r ) & A^-\\
\hline
\end{array}
$$
\caption{ \label{table_I_2}}
\end{table}

\medskip
\noindent
where $A^+$ denotes the problem
$\M \to S \M  (S')^{-1},$ on the set of 
blocks $I=\{i_1,i_2,i_3\}$  with
$$
M = 
\begin{array}{@{}c@{}c@{}}
\begin{array}{@{}c@{}c@{}c@{}}
\smk{_{i_1}}&\smk{_{i_3}}&\smk{_{i_2}}
\end{array}
&
\\
\begin{array}{ |@{}c @{}|@{} c  @{}|@{}    c   @{}|@{}}
\hline
{\mk{a_1}}      &\mk{*}         &     \mk{*}    \\
\hline
 {\mk{*}}     &{\mk{*}}    &     \mk{a_2}      \\
\hline
 {\mk{*}}     &{\mk{a_3}}    &     \mk{*}    \\
\hline
\end{array}
&
\begin{array}{@{}l@{}}
\tmk{\scriptstyle{i_1}}\\\tmk{\scriptstyle{i_2}}\\\tmk{\scriptstyle{i_3}}
\end{array}
\\ \phantom{A}
\end{array}
\;\hbox{ and }\;
(S,S')=
\left(
\begin{array}{@{}c@{}c@{}}
&
\begin{array}{@{}c@{}c@{}c}
\smk{_{i_1}}&\smk{_{i_2}}&\smk{_{i_3}}
\end{array}
\\
\begin{array}{@{}r@{}}
\tmk{\scriptstyle{i_1}}\\\tmk{\scriptstyle{i_2}}\\\tmk{\scriptstyle{i_3}}
\end{array}
&
\begin{array}{ @{}|@{}c @{}|@{} c  @{}|@{}    c   @{}|@{}}
\hline
{\mk{w_1}}      &\mk{0}         &     \mk{0}    \\
\hline
 {\mk{*}}     &{\mk{w_2}}    &     \mk{0}      \\
\hline
 {\mk{*}}     &{\mk{*}}    &     \mk{w_3}    \\
\hline
\end{array}
\\ \phantom{A}
\end{array}
\;
\begin{array}{@{}c@{}c@{}}
\begin{array}{@{}c@{}c@{}c@{}}
\smk{_{i_1}}&\smk{_{i_3}}&\smk{_{i_2}}
\end{array}
&
\\
\begin{array}{ @{}|@{}c @{}|@{} c  @{}|@{}    c   @{}|@{}}
\hline
{\mk{w_1}}      &\mk{0}         &     \mk{0}    \\
\hline
 {\mk{*}}     &{\mk{w_3}}    &     \mk{0}      \\
\hline
 {\mk{*}}     &{\mk{*}}    &     \mk{w_2}    \\
\hline
\end{array}
&
\begin{array}{@{}c@{}}
\tmk{\scriptstyle{i_1}}\\\tmk{\scriptstyle{i_3}}\\\tmk{\scriptstyle{i_2}}
\end{array}
\\ \phantom{A}
\end{array}
\right);
$$
in accordance with our notations, the problem $A^-:$ is
$M\to SM(S')^{-1},$ on the set 
of vertices $I=\{i_1,i_2,i_3\},$ where
$$
M = 
\begin{array}{@{}c@{}l@{}}
\begin{array}{@{}c@{}c@{}c@{}}
\smk{_{i_2}}&\smk{_{i_1}}&\smk{_{i_3}}
\end{array}
& 
\\
\begin{array}{ @{}|@{}c @{}|@{} c  @{}|@{}    c   @{}|@{}}
\hline
{\mk{*}}      &\mk{a_1}         &     \mk{*}    \\
\hline
 {\mk{a_2}}     &{\mk{*}}    &     \mk{*}      \\
\hline
 {\mk{*}}     &{\mk{*}}    &     \mk{a_3}    \\
\hline
\end{array}
&
\begin{array}{@{}r@{}}
\tmk{\scriptstyle{i_1}}\\\tmk{\scriptstyle{i_2}}\\\tmk{\scriptstyle{i_3}}
\end{array}
\\
\phantom{A}
\end{array}
\;\hbox{ and} \;
(S,S')=
\left(
\begin{array}{@{}r@{}c@{}}
&
\begin{array}{@{}c@{}c@{}c}
\smk{_{i_1}}&\smk{_{i_2}}&\smk{_{i_3}}
\end{array}
\\
\begin{array}{@{}r@{}}
\tmk{\scriptstyle{i_1}}\\\tmk{\scriptstyle{i_2}}\\\tmk{\scriptstyle{i_3}}
\end{array}
&
\begin{array}{ @{}|@{}c @{}|@{} c  @{}|@{}    c   @{}|@{}}
\hline
{\mk{w_1}}      &\mk{0}         &     \mk{0}    \\
\hline
 {\mk{*}}     &{\mk{w_2}}    &     \mk{0}      \\
\hline
 {\mk{*}}     &{\mk{*}}    &     \mk{w_3}    \\
\hline
\end{array}
\\
\phantom{A}
\end{array}
\;
\begin{array}{@{}c@{}l@{}}
\begin{array}{@{}c@{}c@{}c@{}}
\smk{_{i_2}}&\smk{_{i_1}}&\smk{_{i_3}}
\end{array}
&
\\
\begin{array}{ @{}|@{}c @{}|@{}  c  @{}|@{}    c   @{}|@{}}
\hline
{\mk{w_2}}      &\mk{0}         &     \mk{0}    \\
\hline
 {\mk{*}}     &{\mk{w_1}}    &     \mk{0}      \\
\hline
 {\mk{*}}     &{\mk{*}}    &     \mk{w_3}    \\
\hline
\end{array}
&
\begin{array}{@{}l@{}}
\tmk{\scriptstyle{i_2}}\\\tmk{\scriptstyle{i_1}}\\\tmk{\scriptstyle{i_3}}
\end{array}
\\ \phantom{A}
\end{array}
\right)
.
$$
Note that since 
matrices $S$ and $S'$ are low triangular,
both problems $A^+$ or $A^-$ can be recognized by the form  
of the matrix $M.$



\subsection{Tacnode curve}
\label{subsec_prim_tacnode}
Analogously as in the previous case, we 
reduce the matrix $\mu_1(0)$ to the identity form and 
the matrix $\mu_2(0)$ to the form
(\ref{id_form}). 
Then for the transformations we have the restrictions:
\begin{equation}
\label{eq_f_restr}
f(0)=F_1(0)\, \hbox{ and } 
F_2(0)\mu_2(0)=\mu_2(0)f(0),
\end{equation}
and consequently $F_1(0)$ is as in  (\ref{eq_F_1(0)}).
By the transformation $f_{\e}$ we
can reduce 
one of the matrices either $\mu_{\e_1}(0)$ or $\mu_{\e_2}(0),$
say $\mu_{\e_2}(0),$ to the zero form.
In the remaining matrix $\M:=\mu_{\e_1}(0):$
the blocks (31),(32), (41) and (42) can be reduced to zero
by the transformation $\frac{d F_k}{d z_0}(0)$ 
and the blocks (21), (23), (41) and (43) can be killed by $f_{\e}.$

\subsection*{Reduced matrix problem}
Thus we obtain a new problem for
the matrix $\M$ with the transformations 
$M\to SM{S}^{-1}$ modulo zero block-entries of $M:$ 
$$
\M:=\mu_{\e_1}(0)=
\begin{array}{@{}c@{}l@{}}
&
\\
\begin{array}{|@{}c @{}  c  @{}|@{}  c @{}  c @{}|@{}}
\hline
\mk{a_1}&\mk{b_{12}}&\mk{b_{13}}&\mk{b_{14}} \\
\mk{0}&\mk{a_{2}}&\mk{0}&\mk{b_{24}} \\
\hline
\mk{0}&\mk{0}&\mk{a_{3}}&\mk{b_{34}} \\
\mk{0}&\mk{0}&\mk{0}&\mk{a_{4}} \\
\hline
\end{array}
&
\\ \phantom{a}
\end{array}
\;\hbox{ and }\;
S:=F_1(0)
=
\begin{array}{c}
\begin{array}{ |@{}c   @{} c @{}|@{} c   @{}   c   @{}|}
\hline
{\mk{w_1}}      &\mk{0}          & \mk{0}&     \mk{0}    \\
{\mk{x_{21}}}      &\mk{w_2}          & \mk{0}&     \mk{0}    \\
\hline
{\mk{x_{31}}}      &\mk{0}          & \mk{w_3}&     \mk{0}    \\
{\mk{x_{41}}}      &\mk{x_{42}}          & \mk{x_{43}}&     \mk{w_4}    \\
\hline
\end{array}.
\end{array}
$$

\noindent
It is easy to see 
that if the sizes of blocks 1 and 4 are both nonzero then
there is a nontrivial endomorphism. 
As in the previous case 
there are no sincere bricks and the admissible tuples of blocks 
$I$ and sizes $\sss$  are the same as in Table \ref{table_I_2},
whereas the configurations 
$A^+$ and  $A^-$ are respectively 
the matrix problems with
$$
\M=
\begin{array}{@{}c@{}l@{}}
&
\begin{array}{@{}c@{}c@{}c@{}}
\smk{_{i_1}}&\smk{_{i_2}}&\smk{_{i_3}}
\end{array}
\\
\begin{array}{@{}r@{}}
\tmk{\scriptstyle{i_1}}\\\tmk{\scriptstyle{i_2}}\\\tmk{\scriptstyle{i_3}}
\end{array}
&
\begin{array}{|@{}c@{}|@{}  c  @{}|@{}  c @{}|@{}}
\hline
\mk{*}&\mk{*}&\mk{*} \\
\hline
\mk{}&\mk{*}&\mk{} \\
\hline
\mk{}&\mk{}&\mk{*} \\
\hline
\end{array}
\\ \phantom{a}
\end{array}
\;\;
S=
\begin{array}{@{}c@{}c@{}}
\begin{array}{@{}c@{}c@{}c@{}}
\smk{_{i_1}}&\smk{_{i_2}}&\smk{_{i_3}}
\end{array} 
&
\\
\begin{array}{ |@{}c   @{}|@{} c @{}|@{} c   @{}|}
\hline
{\mk{w_1}}      &\mk{0}          & \mk{0}        \\
\hline
{\mk{*}}      &\mk{w_2}          & \mk{0}        \\
\hline
{\mk{*}}      &\mk{0}          & \mk{w_3}      \\
\hline
\end{array}
&
\begin{array}{@{}r@{}}
\tmk{\scriptstyle{i_1}}\\\tmk{\scriptstyle{i_2}}\\\tmk{\scriptstyle{i_3}}
\end{array}
\\ \phantom{a}
\end{array}
\;\hbox{ and }\;
\M=
\begin{array}{@{}c@{}l@{}}
&
\begin{array}{@{}c@{}c@{}c@{}}
\smk{_{i_1}}&\smk{_{i_2}}&\smk{_{i_3}}
\end{array}
\\
\begin{array}{@{}r@{}}
\tmk{\scriptstyle{i_1}}\\\tmk{\scriptstyle{i_2}}\\\tmk{\scriptstyle{i_3}}
\end{array}
&
\begin{array}{|@{}c @{}|@{}  c  @{}|@{}  c @{}|@{}}
\hline
\mk{*}&\mk{}&\mk{*} \\
\hline
\mk{}&\mk{*}&\mk{*} \\
\hline
\mk{}&\mk{}&\mk{*} \\
\hline
\end{array}
\\ \phantom{a}
\end{array}
\;\;
S=
\begin{array}{@{}c@{}l@{}}
\begin{array}{@{}c@{}c@{}c@{}}
\smk{_{i_1}}&\smk{_{i_2}}&\smk{_{i_3}}
\end{array}
&
\\
\begin{array}{ |@{}c   @{}|@{} c @{}|@{} c   @{}|}
\hline
{\mk{w_2}}      &\mk{0}          & \mk{0}        \\
\hline
{\mk{0}}      &\mk{w_3}          & \mk{0}        \\
\hline
{\mk{*}}      &\mk{*}          & \mk{w_4}      \\
\hline
\end{array}
&
\begin{array}{@{}r@{}}
\tmk{\scriptstyle{i_1}}\\\tmk{\scriptstyle{i_2}}\\\tmk{\scriptstyle{i_3}}
\end{array}
\\ \phantom{a}
\end{array}
$$
In the matrix $M,$  
we replaced the zero-blocks 
by the empty spaces, since they do not play 
any role in calculations,
and can be always restored for some outer reasons.

\begin{exam} Let $E$ be a Kodaira fiber I$_2$ or III
and $(\fn,\mm,\mn)$ be a triple corresponding 
to a simple vector bundle. 
If $r\geq \bd$ then
the matrix $\mn$ can be respectively transformed to the form
\begin{align*}
\left(
\begin{array}{@{}c@{}}
\begin{array}{ |@{}c @{}@{} c   @{}|@{}   c   @{}|}
\hline
{\mk{\one}}      &\mk{0}         &     \mk{0}    \\
 {\mk{0}}     &{\mk{\one}}    &     \mk{0}      \\
\hline
 {\mk{0}}     &{\mk{0}}    &     \mk{\one}    \\
\hline
\end{array}
\;\;
\begin{array}{ |@{}c @{}@{} c   @{}|@{}   c   @{}|}
\hline
{\mk{*}}      &\mk{*}         &     \mk{*}    \\
 {\mk{*}}     &{\mk{*}}    &     \mk{*}      \\
\hline
 {\mk{*}}     &{\mk{*}}    &     \mk{*}    \\
\hline
\end{array}
\\
\rule{0pt}{33pt}
\begin{array}{@{}c@{}}
\begin{array}{ |@{}c @{} c   @{}|@{}   c   @{}|}
\hline
{\mk{\one}}      &\mk{0}         &     \mk{0}    \\
 {\mk{0}}     &{\mk{0}}    &     \mk{\one}      \\
\hline
 {\mk{0}}     &{\mk{\one}}    &     \mk{0}    \\
\hline
\end{array}
\;\;
\begin{array}{ |@{}c @{} c   @{}|@{}   c   @{}|}
\hline
{\mk{\one}}      &\mk{0}         &     \mk{0}    \\
 {\mk{0}}     &{\mk{\one}}    &     \mk{0}      \\
\hline
 {\mk{0}}     &{\mk{0}}    &     \mk{\one}    \\
\hline
\end{array}
\end{array}
\end{array}
\right )
& 
\;\hbox{ or }\;
\left( 
\begin{array}{@{}c @{}}
\begin{array}{ |@{}c @{}@{} c   @{}|@{}   c   @{}|}
\hline
{\mk{\one}}      &\mk{0}         &     \mk{0}    \\
 {\mk{0}}     &{\mk{\one}}    &     \mk{0}      \\
\hline
 {\mk{0}}     &{\mk{0}}    &     \mk{\one}    \\
\hline
\end{array}
+
\e_1
\begin{array}{ |@{}c @{}@{} c   @{}|@{}   c   @{}|}
\hline
{\mk{*}}      &\mk{*}         &     \mk{*}    \\
 {\mk{0}}     &{\mk{*}}    &     \mk{0}      \\
\hline
 {\mk{0}}     &{\mk{0}}    &     \mk{*}    \\
\hline
\end{array}
\\
\rule{0pt}{33pt}
\begin{array}{ |@{}c @{} c   @{}|@{}   c   @{}|}
\hline
{\mk{\one}}      &\mk{0}         &     \mk{0}    \\
 {\mk{0}}     &{\mk{0}}    &     \mk{\one}      \\
\hline
 {\mk{0}}     &{\mk{\one}}    &     \mk{0}    \\
\hline
\end{array}
+
\e_2
\begin{array}{ |@{}c @{} c   @{}|@{}   c   @{}|}
\hline
{\mk{0}}      &\mk{0}         &     \mk{0}    \\
 {\mk{0}}     &{\mk{0}}    &     \mk{0}      \\
\hline
 {\mk{0}}     &{\mk{0}}    &     \mk{0}    \\
\hline
\end{array}
\end{array}
\right).
\end{align*}
\end{exam}

\medskip
\subsection{Cycle of tree lines}
\label{subsec_prim_cycle3}
According to Section \ref{sec_matr_pr_cycles} 
the original matrix problem $\MP_E$ with two blocks on each 
component is 
$$
\xy ;/r0.15pc/:
{\ar@/_5pt/_{F_1(0)}(-2,18);(-2,2)},
{\ar@{.}_{f_1}(10,-1);(10,-39)},
%
{\ar@{-}(0,0);(0,20)}, 
{\ar@{-}(20,0);(20,20)},
{\ar@{-}(0,20);(20,20)}, 
{\ar@{-}^{}(0,10);(20,10)},  
{\ar@{-}^{}(0,0);(20,0)}, 
{\ar@/^5pt/^{F_1(\infty)}(62,18);(62,2)},
{\ar@{-}(40,0);(40,20)}, 
{\ar@{-}(40,20);(60,20)}, {\ar@{-}(60,0);(60,20)},
{\ar@{-}^{}(40,10);(60,10)},
{\ar@{.}(40,15);(20,15)},
 {\ar@{.}(40,5);(20,5)}, 
{\ar@{-}^{}(40,0);(60,0)}, 
{\ar@{.}^{f_2}(50,-1);(50,-9)},
{\ar@/_5pt/_{F_2(0)}(38,-12);(38,-28)},
%
{\ar@{-}(40,-30);(40,-10)}, 
{\ar@{-}(40,-10);(60,-10)}, {\ar@{-}(60,-30);(60,-10)},
{\ar@{-}^{}(40,-20);(60,-20)},
{\ar@{.}(80,-15);(60,-15)}, {\ar@{.}(80,-25);(60,-25)}, 
{\ar@{-}^{}(40,-30);(60,-30)}, 
{\ar@{-}(80,-30);(80,-10)}, {\ar@{-}(100,-30);(100,-10)},
{\ar@{-}(80,-10);(100,-10)}, 
{\ar@{-}^{}(80,-20);(100,-20)},  
{\ar@{-}^{}(80,-30);(100,-30)}, 
{\ar@{.}^{f_3}(90,-31);(90,-39)},
{\ar@/^5pt/^{F_2(\infty)}(102,-12);(102,-28)},
{\ar@/_5pt/_{F_3(\infty)}(-2,-42);(-2,-58)},
{\ar@{-}(0,-60);(0,-40)}, {\ar@{-}(20,-60);(20,-40)},
{\ar@{-}(0,-40);(20,-40)}, 
{\ar@{.}(20,-45);(80,-45)},
{\ar@{.}(20,-55);(80,-55)},
{\ar@{-}^{}(0,-50);(20,-50)},  
{\ar@{-}^{}(0,-60);(20,-60)}, 
{\ar@{-}(80,-60);(80,-40)}, {\ar@{-}(100,-60);(100,-40)},
{\ar@{-}(80,-40);(100,-40)}, 
{\ar@{-}^{}(80,-50);(100,-50)},  
{\ar@{-}^{}(80,-60);(100,-60)}, 
{\ar@/^5pt/^{F_3(0)}(102,-42);(102,-58)},
\endxy
$$
Matrices $\mu_1(0),$ $\mu_2(0)$ and $\mu_3(0)$ can 
be reduced to the identity form. The matrix 
$\mu_3(\infty)$ can be reduced to the form (\ref{id_form}).
For the morphisms we have
\begin{equation}
\label{eq_for_F_I3}
 f_1=F_1(0),\, f_2=F_2(0), \, f_3=F_3(0)
\hbox{ and } F_3(\infty)\mu_3(\infty)=\mu_3(\infty)f_1.
\end{equation}
Then the matrix $f_3$ becomes a special block-triangular structure.
In other words, the matrix $\mu_2(\infty)$ is subdivided into 
four column-blocks:
a column can be added to any other column 
from a block on the left and it cannot be added 
to a column from another block on the right. 
Thus $\mu_2(\infty)$ can be reduced to the form
\begin{align}
\label{eq_redfrm}
\mu_2(\infty)=
\begin{array}{@{}c@{}l@{}}
\begin{array}{  @{}c @{}c @{}c @{}c @{}c @{}c @{}c @{}c @{} }
\smk{_{1}}&\smk{_{2}}&\smk{_{5}}&\smk{_{6}}
\smk{_{3}}&\smk{_{4}}&\smk{_{7}}&\smk{_{8}}
\end{array}
&
\\
\begin{array}{|@{}c@{}c@{}|@{}c@{}c@{}|@{}c@{}c@{}|@{}c@{}c@{}|}
\hline
    \mk{\one}&\mk{0}  &\mk{0}&\mk{0}   &\mk{0}&\mk{0}        &\mk{0} &\mk{0}\\
    \mk{0}&\mk{0}   &\mk{\one}&\mk{0}   &\mk{0}&\mk{0}        &\mk{0} &\mk{0}\\
    \mk{0}&\mk{0}   &\mk{0}&\mk{0}     &\mk{\one}&\mk{0}     &\mk{0} &\mk{0}\\
    \mk{0}&\mk{0}   &\mk{0}&\mk{0}     &\mk{0}&\mk{0}       &\mk{\one}&\mk{0}\\
\hline
    \mk{0}&\mk{\one}   &\mk{0}&\mk{0}   &\mk{0}&\mk{0}        &\mk{0} &\mk{0}\\
    \mk{0}&\mk{0}   &\mk{0}&\mk{\one}   &\mk{0}&\mk{0}        &\mk{0} &\mk{0}\\
    \mk{0}&\mk{0}   &\mk{0}&\mk{0}   &\mk{0}&\mk{\one}     &\mk{0} &\mk{0}\\
    \mk{0}&\mk{0}   &\mk{0}&\mk{0}   &\mk{0}&\mk{0}       &\mk{0} &\mk{\one}\\
\hline
    \end{array}
&\begin{array}{@{}l@{}}
\tmk{\scriptstyle{{1}}}\\\tmk{\scriptstyle{{5}}}\\ 
\tmk{\scriptstyle{{3}}}\\\tmk{\scriptstyle{{7}}}\\
\tmk{\scriptstyle{{2}}}\\\tmk{\scriptstyle{{6}}}\\ 
\tmk{\scriptstyle{{4}}}\\\tmk{\scriptstyle{{8}}}
\end{array}
\\ \phantom{a}
\end{array}
\end{align}
%
%

\subsection*{Reduced matrix problem}
The remaining nonreduced matrix is $\M:=\mu_1(\infty).$
For it we obtain the problem $M\to SM(S')^{-1},$ where 
the transformations are $(S,S')=(F_1(\infty),f_2).$
Equations (\ref{eq_for_F_I3}) together with 
${F_2(\infty)\mu_2(\infty)=\mu_2(\infty)f_3 }$ imply 
the triangular forms for the matrices $F_k(0),F_k(\infty)$ and $f_k$
(for $k=1,2,3$); in particular:
\begin{align*}
(S,S')=
\left(
\begin{array}{@{}r@{}c@{}}
&
\begin{array}{@{}c@{}c@{}c@{}c@{}c@{}c@{}c@{}c@{}}
\mk{_1}&\mk{_2}&\mk{_3}&\mk{_4}\mk{_5}&\mk{_6}&\mk{_7}&\mk{_8}
\end{array}
\\
\begin{array}{@{}l@{}}
        \tmk{\scriptstyle{1}}\\\tmk{\scriptstyle{2}}\\\tmk{\scriptstyle{3}}\\
        \tmk{\scriptstyle{4}}\\
        \tmk{\scriptstyle{5}}\\
        \tmk{\scriptstyle{6}}\\\tmk{\scriptstyle{7}}\\\tmk{\scriptstyle{8}}
        \end{array}
        \,
&
\begin{array}{|@{}c@{}c@{}|@{}c@{}c@{}|@{}c@{}c@{}|@{}c@{}c@{}|@{}}
\hline
    \mk{w_1}& \mk{0}   &0&0     &0&0       &\mk{0}&\mk{0}\\
    \mk{*}& \mk{w_2}   &0&0     &0&0       &\mk{0}&\mk{0}\\
\hline
    \mk{*}&\mk{*}   &\mk{w_{3}}&\mk{0}     &\mk{0}&\mk{0}     &\mk{0} &\mk{0}\\
    \mk{*}&\mk{*}   &\mk{*}&\mk{w_{4}}     &\mk{0}&\mk{0}     &\mk{0} &\mk{0}\\
\hline
    \mk{*}&\mk{*}   &\mk{*}&\mk{*}     &\mk{w_5}&\mk{0}     &\mk{0} &\mk{0}\\
    \mk{*}&\mk{*}   &\mk{*}&\mk{*}     &\mk{*}&\mk{w_{6}}  &\mk{0} &\mk{0}\\
\hline
    \mk{*}&\mk{*}   &\mk{*}&\mk{*}     &\mk{*}&\mk{*}     &\mk{w_{7}} &\mk{0}\\
    \mk{*}&\mk{*}   &\mk{*}&\mk{*}     &\mk{*}&\mk{*}     &\mk{*} &\mk{w_{8}}\\
\hline
    \end{array}
\\ \phantom{a}
\end{array}
\;\;
\begin{array}{@{}c@{}l@{}}
\begin{array}{@{}c@{}c@{}c@{}c @{}c@{}c@{}c@{}c@{}}
\mk{_1}&\mk{_5}&\mk{_3}&\mk{_7}\mk{_2}&\mk{_6}&\mk{_4}&\mk{_8}\\
\end{array}
&
\\
\begin{array}{@{}|@{}c@{}c@{}|@{}c@{}c@{}|@{}c@{}c@{}|@{}c@{}c@{}|@{}}
\hline
    \mk{w_1}& \mk{0}   &\mk{0}&\mk{0}     &0&0       &\mk{0}&\mk{0}\\
    \mk{*}& \mk{w_{5}}   &\mk{0}&\mk{0}     &0&0       &\mk{0}&\mk{0}\\
\hline
    \mk{*}&\mk{*}   &\mk{w_{3}}&\mk{0}     &\mk{0}&\mk{0}     &\mk{0} &\mk{0}\\
    \mk{*}&\mk{*}   &\mk{*}&\mk{w_{7}}     &\mk{0}&\mk{0}     &\mk{0} &\mk{0}\\
\hline
    \mk{*}&\mk{*}   &\mk{*}&\mk{*}     &\mk{w_{2}}&\mk{0}     &\mk{0} &\mk{0}\\
    \mk{*}&\mk{*}   &\mk{*}&\mk{*}     &\mk{*}&\mk{w_{6}}     &\mk{0} &\mk{0}\\
\hline
    \mk{*}&\mk{*}   &\mk{*}&\mk{*}     &\mk{*}&\mk{*}     &\mk{w_4} &\mk{0}\\
    \mk{*}&\mk{*}   &\mk{*}&\mk{*}     &\mk{*}&\mk{*}     &\mk{*} &\mk{w_{8}}\\
\hline
    \end{array}
&
\begin{array}{@{}l@{}}
\tmk{_1}\\\tmk{_5}\\\tmk{_3}\\\tmk{_7} \\ \tmk{_2}
\\\tmk{_6}\\\tmk{_4}\\\tmk{_8}
\end{array}
\\ \phantom{a}
\end{array}
\right)
\end{align*}
The stars $*$
denote arbitrary 
blocks and $w_i$ for $i\in\{1,\dots, 8\}$ are the common diagonal blocks.
The transformations of row and column-blocks of $\M$ are clear:
a row can be added to any other one from a block 
below and it can not be added to a row from a block above it;
and a  column can be added to any other column 
from a block on the left and it can not be added 
to a column from a block on the right.

\subsection*{Nontrivial endomorphisms}
Analogously as in the case of a cycle of two lines 
there are some pairs $(ij) \in I\times I$ such that 
if $s_i\cdot s_j>0 $ then there exists a nontrivial endomorphism.
Such blocks are  called \emph{mutually excluding} and 
denoted by $i\cap j.$  
\begin{itemize}
  \item 
 If the matrices $F_3(\infty)$ and $F_1(0)$
contain at least one of the following entries:
$(71),$  $(81),$ $(72)$ or $(82)$
then there is a nontrivial endomorphism.
In our short notations we have intersections $1,2 \cap  7,8.$ 

\item
Analogously we have $1,5 \cap 4,8$ 
coming from the matrices $F_3(0)$ and $F_2(\infty).$

 \item The blocks $1$ and  $6 $  are mutually excluding;
the endomorphism
is induced by the entry
(61) of the matrices  $F_3(0),$ $F_3(\infty)$ and $F_2(\infty).$
Similarly, 
there is an endomorphism 
for the pair $(38)$
induced by the matrices  $F_1(0),$ $F_3(\infty)$ and $F_3(0).$
\end{itemize}

\noindent
All the mutually excluding blocks can be indicated on the 
\emph{intersection diagram}:
\begin{equation}
\label{eq_dia}
\begin{smallmatrix}
 && 1\\
2&3&    &-&5\\
&&\bigcap&&\\
7&-&    &6&4\\
&& 8
\end{smallmatrix}
\end{equation}
The diagram reads as follows: 
a matrix $\M$ is a brick
if it contains no pair of blocks $(ij)$
such that
$i$ and $j$ in the diagram are 
separated by $\cap$ and either 
in the same column or 
one of them is $1$ or $8. $ 

In the following table we
present the  maximal tuples of blocks $I=(i_1,i_2,i_3,i_4)$
for $M$ being a brick, 
express the dimension vector
$\sss=(s_{i_1},s_{i_2},s_{i_3},s_{i_4})\in\mN^4$ in terms of 
rank and multidegree and 
moreover, 
 answer the question when such tuple of 
blocks appears and specialize the matrix problems in each case. 
\begin{table}[ht]
$$
\begin{array}{@{}|c@{} |c|c |c|@{}c@{}|@{}}
\hline
&
\hbox{condition}& \hbox{ set }\; I& \hbox{dimension vector }\; \sss
& \hbox{ state }\\
\hline
\hline
1.&
\rule{0pt}{12pt} 
r\geq \bd & (1,2,3,5)& (r-\bd,\bd_2,\bd_3,\bd_1) &A^+ \\
\hline
1'.&
\rule{0pt}{12pt} 
\bd \geq 2r& (4,6,7,8)& (r-\bd_1,r-\bd_3,r-\bd_2, \bd-2r ) & A^-\\
\hline\hline
2.&
\rule{0pt}{12pt} 
\bd>r> (\bd_2+\bd_3), (\bd_1+\bd_3)&
(2,3,5,6)& (r-(\bd_1+\bd_3),\bd_3,r-(\bd_2+\bd_3),\bd-r) & A^- \\
\hline
2'.&
\rule{0pt}{12pt} 
(\bd_2+\bd_3), (\bd_1+\bd_3)>r&
(3,4,6,7)&
(2r-\bd, 
(\bd_2+\bd_3)-r,r-\bd_3, (\bd_1+\bd_3)-r) & A^+ \\
\hline\hline
3.&
\rule{0pt}{12pt} 
(\bd_2+\bd_3)\geq r \geq  (\bd_1+\bd_3) &
(2,3,4,6)& (r-(\bd_1+\bd_3),r-\bd_2,(\bd_2+\bd_3) -r,\bd_1)& C\\
\hline
3'. &
\rule{0pt}{12pt} 
(\bd_1+\bd_3)\geq r \geq  (\bd_2+\bd_3)&
(3,5,6,7) &
(r-\bd_1, r-(\bd_2+\bd_3),\bd_2,(\bd_1+\bd_3) -r)& C \\
\hline
\end{array}
$$
\caption{ \label{table_I_3}}
\end{table}

\noindent
The configurations $A^+,$ $A^-$  and $C$ 
on the set of blocks $I=\{i_1,i_2,i_3,i_4\}$
encode matrix problems $M\to SM(S')^{-1},$
where $S$ and $S'$ are block-triangular and the matrix 
$\M$ is defined as follows:
\begin{align}
A^+=
\begin{array}{@{}c@{}c@{}}
\begin{array}{  @{}c @{}c @{}c @{}c @{}}
\mk{_{i_1}}&\mk{_{i_4}}&\mk{_{i_3}}&\mk{_{i_2}}
\end{array}
\\
\begin{array}{ |@{}c @{}|@{} c  @{}|@{}    c   @{}|@{} c   @{}|@{}}
\hline
{\mk{a_{i_1}}}      &\mk{*}         &     \mk{*}    &     \mk{*}\\
\hline
 {\mk{*}}     &{\mk{*}}    &     \mk{*}  &     \mk{a_{i_2}}     \\
\hline
{\mk{*}}     &{\mk{*}}    &     \mk{a_{i_3}}  &     \mk{*}     \\
\hline
{\mk{*}}     &{\mk{a_{i_4}}}    &     \mk{*}  &     \mk{*}     \\
\hline
\end{array}
&
\begin{array}{@{}r@{}}
\tmk{\scriptstyle{{i_1}}}\\\tmk{\scriptstyle{{i_2}}}\\ 
\tmk{\scriptstyle{{i_3}}}\\\tmk{\scriptstyle{{i_4}}}
\end{array}
\\\phantom{a}
\end{array}
\;\;\;
\A^-= 
\begin{array}{@{}c@{}l@{}}
\begin{array}{  @{}c @{}c @{}c @{}c @{}}
\mk{_{i_3}}&\mk{_{i_2}}&\mk{_{i_1}}&\mk{_{i_4}}
\end{array}
\\
\begin{array}{ |@{}c @{}|@{} c  @{}|@{}    c   @{}|@{} c   @{}|@{}}
\hline
{\mk{*}}     &{\mk{*}}    &     \mk{a_{i_1}}  &     \mk{*}     \\
\hline
{\mk{*}}     &{\mk{a_{i_2}}}    &     \mk{*}  &     \mk{*}     \\
\hline
{\mk{a_{i_3}}}      &\mk{*}         &     \mk{*}    &    \mk{*} \\
\hline
 {\mk{*}}     &{\mk{*}}    &     \mk{*}  &     \mk{a_{i_4}}     \\
\hline
\end{array}
&
\begin{array}{@{}l@{}}
\tmk{\scriptstyle{{i_1}}}\\\tmk{\scriptstyle{{i_2}}}\\ 
\tmk{\scriptstyle{{i_3}}}\\\tmk{\scriptstyle{{i_4}}}
\end{array}
\\\phantom{a}
\end{array}
\;\;\hbox{ and }
C=\begin{array}{@{}c@{}l@{}}
\begin{array}{  @{}c @{}c @{}c @{}c @{}}
\mk{_{i_2}}&\mk{_{i_1}}&\mk{_{i_4}}&\mk{_{i_3}}
\end{array}
\\
\begin{array}{ |@{}c @{}|@{} c  @{}|@{}    c   @{}|@{} c   @{}|@{}}
\hline
{\mk{*}}     &{\mk{a_{i_1}}}    &     \mk{*}  &     \mk{*}     \\
\hline
{\mk{a_{i_2}}}      &\mk{*}         &     \mk{*}    &    \mk{*} \\
\hline
 {\mk{*}}     &{\mk{*}}    &     \mk{*}  &     \mk{a_{i_3}}     \\
\hline
{\mk{*}}     &{\mk{*}}    &     \mk{a_{i_4}}  &     \mk{*}     \\
\hline
\end{array}
&
\begin{array}{@{}l@{}}
\tmk{\scriptstyle{{i_1}}}\\\tmk{\scriptstyle{{i_2}}}\\ 
\tmk{\scriptstyle{{i_3}}}\\\tmk{\scriptstyle{{i_4}}}
\end{array}
\\\phantom{a}
\end{array}.
\end{align}

\subsection{Thee concurrent lines in a plane}
\label{subsec_prim_IV}
Let $E$ be the Kodaira fiber IV
and
$\MP_E$ 
the matrix problem  
formulated in Section \ref{sec_matr_pr_fibers}
with two blocks for each component. 
Reduce matrices $\mu_1(0)$
and $\mu_2(0)$ to the identity and (\ref{id_form}) forms
respectively,
as in the case of a tacnode curve.  
Then the transformations satisfy 
equations (\ref{eq_f_restr}).
Let us find a canonical form of $\mu_3(0)$
with respect to the transformations 
$$\mu_3(0)\to F_3(0)\mu_3(0)f(0)^{-1}.$$
The splitting of $F_3(0)$ and $f(0)$ into blocks 
induces the  same  column block
structure 
for $\mu_3(0)$
as in the case of a cycle of three lines.
However, on the contrary to that case,
there is no addition from the third column-block to the second one.
Thus proceeding as before
instead of the form (\ref{eq_redfrm})
we obtain only the following:
\begin{align*}
\mu_3(0)=    
\begin{array}{@{}c@{}l@{}}
\begin{array}{  @{}c @{}c @{}c @{}c @{}c @{}c @{}c @{}c @{} }
\smk{_{1}}&\smk{_{2}}&\smk{_{3}}&\smk{_{4}}&\smk{_{5}}&\smk{_{6}}&\smk{_{7}}&\smk{_{8}}
\end{array}
&
\\
\begin{array}{|@{}c@{}c@{}|@{}c@{}c@{}|@{}c@{}c@{}|@{}c@{}c@{}|}
\hline
    \mk{\one}&\mk{0}  &\mk{0}&\mk{0}   &\mk{0}&\mk{0}        &\mk{0} &\mk{0}\\
    \mk{0}&\mk{0}   &\mk{\one}&\mk{0}   &\mk{0}&\mk{0}        &\mk{0} &\mk{0}\\
    \mk{0}&\mk{0}   &\mk{0}&\mk{*}     &\mk{\one}&\mk{0}     &\mk{0} &\mk{0}\\
    \mk{0}&\mk{0}   &\mk{0}&\mk{0}     &\mk{0}&\mk{0}       &\mk{\one}&\mk{0}\\
    \hline
    \mk{0}&\mk{\one}   &\mk{0}&\mk{0}   &\mk{0}&\mk{0}        &\mk{0} &\mk{0}\\
    \mk{0}&\mk{0}   &\mk{0}&\mk{\one}   &\mk{0}&\mk{0}        &\mk{0} &\mk{0}\\
    \mk{0}&\mk{0}   &\mk{0}&\mk{0}   &\mk{0}&\mk{\one}     &\mk{0} &\mk{0}\\
    \mk{0}&\mk{0}   &\mk{0}&\mk{0}   &\mk{0}&\mk{0}       &\mk{0} &\mk{\one}\\
\hline
    \end{array}
&\begin{array}{@{}l@{}}
\tmk{\scriptstyle{{1}}}\\\tmk{\scriptstyle{{3}}}\\ 
\tmk{\scriptstyle{{5}}}\\\tmk{\scriptstyle{{7}}}\\
\tmk{\scriptstyle{{2}}}\\\tmk{\scriptstyle{{4}}}\\
\tmk{\scriptstyle{{6}}}\\\tmk{\scriptstyle{{8}}}
\end{array}
\\ \phantom{a}
\end{array}
\end{align*}
It turns out that the remaining  block $*$
can be  reduced 
to the form $(\begin{smallmatrix}0&0\\ \one &0\end{smallmatrix})$ as well.
That implies subdivisions for the reduced blocks marked by 4 and 5
and change of notations is required: 
$4\to (0,4) $ and $5\to (5,0).$ 
The equation $F_3(0)\mu_3(0)=\mu_3(0)f(0)$ 
implies that the matrix $F_1(0)$ 
preserving $\mu_3(0)$ is as follows: 
\vspace{-0.1cm}
  \begin{align*}
F_1(0)=
\begin{array}{c}
    \begin{array}{r@{}|@{}c@{}c@{}|@{}c@{}c@{}c@{}|@{}c@{}c@{}c@{}|@{}c@{}c@{}|@{}}
    \multicolumn{1}{c}{}
    &\multicolumn{1}{c}{_1}
    &\multicolumn{1}{c}{_2}
    &\multicolumn{1}{c}{_3}
    &\multicolumn{1}{c}{_0}
    &\multicolumn{1}{c}{_4}
    &\multicolumn{1}{c}{_5}
    &\multicolumn{1}{c}{_0}
    &\multicolumn{1}{c}{_6}
    &\multicolumn{1}{c}{_7}
    &\multicolumn{1}{c}{_8}
    \\
    \cline{2-11}
    _1  &w_0&0    &0&0&0 &0&0&0 &0&0\\
    _2  &*&w_1    &0&0&0 &0&0&0 &0&0\\
    \cline{2-11}
    _3  &*&0    &w_2&0&0  &0&0&0 &0&0\\
    _0  &*&x    &*&z&0  &0&0&0 &0&0\\
    _4  &*&*    &*&*&w_4  &0&0&0 &0&0\\
    \cline{2-11}
    _5  &*&0       &0&0&0 &w_5&0&0  &0&0\\
    _0  &*&x       &0&0&0 &*&z&0  &0&0\\
    _6  &*&*       &0&0&0 &*&*&w_7  &0&0\\
    \cline{2-11}
    _7  &*&0 &*&y&0  &*&y&0   &w_8&0\\
    _8  &*&* &*&*&*  &*&*&*   &*&w_9\\
    \cline{2-11}
    \end{array}
\\ \phantom{a} 
    \end{array}
\end{align*}

\vspace{-0.3cm}
\noindent
As usually the stars $*$
denote different blocks appearing only  one time and
$x,$ $y$ and $z$  are some blocks appearing twice.
By proper $f_x(0)$ and $f_y(0)$ 
the matrices $\mu_{\e_2}(0)$ and $\mu_{\e_3}(0)$ can be reduced
to zero.

\subsection*{Reduced matrix problem}
As usually take
$\M:=\mu_{\e_1}(0)$ and transformations 
$M\to S MS^{-1}$ modulo zero blocks of $M,$ where
 $S:=F_1(0).$
By proper $F_1(0),$ $f_x(0)$ and $f_y(0)$ it
can be reduced to the form 
\begin{align*}
M=
    \begin{array}{c}
    \begin{array}{r@{}|@{}c@{}c@{}|@{}c@{}c@{}c@{}|@{}c@{}c@{}c@{}|@{}c@{}c@{}|@{}}
    \multicolumn{1}{c}{}
    &\multicolumn{1}{c}{_1}
    &\multicolumn{1}{c}{_2}
    &\multicolumn{1}{c}{_3}
    &\multicolumn{1}{c}{_0}
    &\multicolumn{1}{c}{_4}
    &\multicolumn{1}{c}{_5}
    &\multicolumn{1}{c}{_0}
    &\multicolumn{1}{c}{_6}
    &\multicolumn{1}{c}{_7}
    &\multicolumn{1}{c}{_8}
    \\
    \cline{2-11}
    _1  &*  &*    &*&*&*      &* &* &*       &*&*\\
    _2  &0 &*    &0 &0_{y}&* &0&0_{y} &* &0 &*\\
    \cline{2-11}
    _3  &0  &0    &*&*&*      &0 &0 &0       &*&*\\
    _0  &0  &0    &0 &0_{z}&* &0 &0 &0       &0_{x} &*\\
    _4  &0  &0    &0 &0 &*    &0 &0 &0       &0 &*\\
    \cline{2-11}
    _5  &0 &0     &0 &0 &0       &*&* &* &* &*\\
    _0  &0 &0     &0 &0 &0       &0 &0_{z}&*     &0_{x}&*\\
    _6  &0 &0     &0 &0 &0       &0 &0 &*      &0 &*\\
    \cline{2-11}
    _7  &0  &0    &0 &0 &0       &0 &0 &0       &*&*\\
    _8  &0  &0    &0 &0 &0       &0 &0 &0       &0 &*\\
    \cline{2-11}
    \end{array}
\\ \phantom{a}
    \end{array}
\end{align*}

\vspace{-0.3cm}
\noindent
The blocks denoted by $0_{x}$ (respectively $0_{y}$ or $0_{z}$) 
are the so called {\it adjoint blocks},
that means
there is a unique block
$x$ (respectively $y$ or $z$) 
operating on both of them, and thus only one block from
an adjoint pair can be reduced to zero.

\subsection*{Nontrivial endomorphisms}
Let us analyze matrices 
$\frac{d F_k}{d z_0}(0),$ $f_x(0)$ and $f_y(0)$ 
looking for an endomorphism.
Taking into account 
equations $F_k(0)\mu_k(0)=\mu_k(0)f(0)$ for $k=2,3$
we see that there are nonzero matrices $f_x(0)$ and $f_{y}(0)$ 
leaving the matrices $\mu_{\e_2}(0)$ and $\mu_{\e_3}(0)$
in the zero form.
Hence, as in the case of a tacnode curve, 
there are places $(ij),$
where zero can be obtained in two or more different ways
(that is if $s_i\cdot s_j> 0$  then 
there exists a nonscalar endomorphism).
The diagram of mutually excluding blocks
is almost the same as the diagram (\ref{eq_dia}): 
\vspace{-0.3cm}
\begin{align}
\label{eq_mutexclud_IV}
\begin{smallmatrix}
 &&& 1\strut\\
2&0&&3&-&5\\
&&&\bigcap&&\\
7&-&&6&0&4\\
&&&8\strut
\end{smallmatrix}
\end{align}
where intersections 
\begin{itemize}
  \item 
$1,2\cap 7,8$ are induced  by   
$\frac{d F_1}{d z_0}(0)$ and  $\frac{d F_2}{d z_0}(0);$ 
  \item 
$1,3\cap 6,8$  by   $\frac{d F_1}{d z_0}(0)$  and  $\frac{d F_3}{d z_0}(0);$

\item
$1,5\cap 4,8$  by $\frac{d F_2}{d z_0}(0)$ and  $\frac{d F_3}{d z_0}(0);$
\item
and entries $(1,0)$ and $(0,8)$ link all three matrices
$\frac{d F_1}{d z_0}(0),$ $\frac{d F_2}{d z_0}(0)$ and  $\frac{d F_3}{d z_0}(0) ; $
\end{itemize}

\noindent
In Table \ref{table_IV} we present the maximal tuples $I=\{i_1,i_2,i_3,i_4\},$ 
interpret the dimension vector $\sss$ in terms of 
rank and multidegree and 
specialize 
matrices that we get in each case. 

\begin{table}[t]
{
$$ 
\begin{array}{@{}|c@{} |c|c |c|@{}c@{}|@{}}
\hline
&
\hbox{condition}& \hbox{ set }\; I& \hbox{dimension vector }\; \sss
& \hbox{state}\\
\hline
\hline
1.&
\rule{0pt}{12pt} 
r\geq \bd & (1,2,3,5)& \big(r-\bd,\bd_3,\bd_2,\bd_1 \big) &A^+ \\
\hline
1'.&
\rule{0pt}{12pt} 
\bd > 2r& (4,6,7,8)& \big(r-\bd_1,r-\bd_2,r-\bd_3,\bd-2r\big) & A^-\\
\hline\hline
2.&
\rule{0pt}{12pt} 
\begin{array}{@{}c@{}}
\bd>r> \bd_i+\bd_j  \rule{0pt}{12pt} \\
\hbox{ for all  } i,j\in\{1,2,3\};
\end{array} 
&
(2,3,5,0)&\Big(r-(\bd_1+\bd_2), r-(\bd_1+\bd_3),r-(\bd_2+\bd_3),\bd-r \Big) 
& A^- \\
\hline
2'.&
\rule{0pt}{12pt} 
\begin{array}{@{}c@{}}
\bd_i+\bd_j>r  \hbox{ and } 2r>\bd \rule{0pt}{12pt} \\
\hbox{ for all  } i,j\in\{1,2,3\};
\end{array} 
& (0,4,6,7)&
\Big(2r-\bd, (\bd_2+\bd_3)-r, (\bd_1+\bd_3)-r,(\bd_1+\bd_2)-r\Big)
 & A^+ \\
\hline\hline
3.&
\begin{array}{@{}c@{}}
(\bd_2+\bd_3)>r \hbox{ and } \rule{0pt}{12pt} \\
 r>(\bd_1+\bd_2),(\bd_1+\bd_3) 
\end{array}
&
(2,3,0,4)& 
\Big(r-(\bd_1+\bd_2), r-(\bd_1+\bd_3),\bd_1, (\bd_2+\bd_3)-r \Big)
& B^-(0) \\
\hline
3'.&
\begin{array}{@{}c@{}}
(\bd_1+\bd_2),(\bd_1+\bd_3)>r \rule{0pt}{12pt} \\
r> (\bd_2+\bd_3)
\end{array}
&
(5,0,6,7) &
\Big(r-(\bd_2+\bd_3), r-\bd_1,(\bd_1+\bd_3)-r, (\bd_1+\bd_2)-r\Big)
& \B^+(0)\\
\hline\hline
4.&
\begin{array}{@{}c@{}}
(\bd_1+\bd_3),(\bd_2+\bd_3)>r  \rule{0pt}{12pt} \\
\hbox{ and } r>(\bd_1+\bd_2), 
\end{array}
&
(2,0,4,6)& 
\Big(r-(\bd_1+\bd_2), r-\bd_3,(\bd_2+\bd_3)-r, (\bd_1+\bd_3)-r \Big)
& B^+(0) \\
\hline
4'.&
\begin{array}{@{}c@{}}
(\bd_1+\bd_2)>r \hbox{ and } \rule{0pt}{12pt} \\
 r>(\bd_1+\bd_3),(\bd_2+\bd_3) 
\end{array}
&
(3,5,0,7) &
\Big(r-(\bd_1+\bd_3), r-(\bd_2+\bd_3),\bd_3, (\bd_1+\bd_2)-r\Big)
& \B^-(0)
\\
\hline\hline
5.&
\begin{array}{@{}c@{}}
(\bd_1+\bd_3)>r \hbox{ and } \rule{0pt}{12pt} \\
 r>(\bd_1+\bd_2),(\bd_2+\bd_3) 
\end{array}
&
(2,5,0,6)& 
\Big(r-(\bd_1+\bd_2), r-(\bd_2+\bd_3),\bd_2, (\bd_1+\bd_3)-r \Big)
& B^-(0) \\
\hline
5'.&
\begin{array}{@{}c@{}}
(\bd_1+\bd_2),(\bd_2+\bd_3)>r  \rule{0pt}{12pt} \\
\hbox{ and } r>(\bd_1+\bd_3)
\end{array}
&
(2,0,4,7) &
\Big(r-(\bd_1+\bd_3), r-\bd_2, (\bd_2+\bd_3)-r, (\bd_1+\bd_2)-r\Big)& 
B^+(0)\\
\hline
\end{array}
$$
}
\caption{ \label{table_IV}}
\end{table}

\noindent
By $A^\sigma $ and
$B^\sigma(j)$ we denote the matrix problems
given by the following coincidence matrices $\M$:
\begin{align}
\label{form_A}
\begin{array}{@{}c@{}l@{}}
\begin{array}{  @{}c @{}c @{}c @{}c @{}}
\mk{_{i_1}}&\mk{_{i_2}}&\mk{_{i_3}}&\mk{_{i_4}}
\end{array}
\\
\begin{array}{ @{}|@{}c   @{}|@{} c @{}|@{} c   @{}|@{}   c   @{}|@{}}
\hline
{\mk{*}}      &\mk{*}          & \mk{*}&     \mk{*}    \\
\hline
{\mk{ }}      &\mk{*}          & \mk{ }&     \mk{ }    \\
\hline
{\mk{ }}      &\mk{ }          & \mk{*}&     \mk{ }    \\
\hline
{\mk{ }}      &\mk{ }          & \mk{ }&     \mk{*}    \\
\hline
\end{array}
&
\begin{array}{@{}l@{}}
\tmk{\scriptstyle{{i_1}}}\\\tmk{\scriptstyle{{i_2}}}\\ 
\tmk{\scriptstyle{{i_3}}}\\\tmk{\scriptstyle{{i_4}}}
\end{array}
\\
_{A^{+}}\phantom{a}
\end{array}
\;\;
\begin{array}{@{}c@{}l@{}}
\begin{array}{  @{}c @{}c @{}c @{}c @{}}
\mk{_{i_1}}&\mk{_{i_2}}&\mk{_{i_3}}&\mk{_{i_4}}
\end{array}
\\
\begin{array}{ @{}|@{}c   @{}|@{} c @{}|@{} c   @{}|@{}   c   @{}|@{}}
\hline
{\mk{*}}      &\mk{}          & \mk{}&     \mk{*}    \\
\hline
{\mk{}}      &\mk{*}          & \mk{}&     \mk{*}    \\
\hline
{\mk{}}      &\mk{}          & \mk{*}&     \mk{*}    \\
\hline
{\mk{}}      &\mk{}          & \mk{}&     \mk{*}    \\
\hline
\end{array}
&
\begin{array}{@{}l@{}}
\tmk{\scriptstyle{{i_1}}}\\\tmk{\scriptstyle{{i_2}}}\\ 
\tmk{\scriptstyle{{i_3}}}\\\tmk{\scriptstyle{{i_4}}}
\end{array} \\
_{A^{-}}\phantom{a}
\end{array}
\;\;
\begin{array}{@{}c@{}l@{}}
\begin{array}{  @{}c @{}c @{}c @{}c @{}}
\mk{_{i_1}}&\mk{_{i_2}}&\mk{_{i_3}}&\mk{_{i_4}}
\end{array}
\\
\begin{array}{ @{}|@{}c   @{}|@{} c @{}|@{} c   @{}|@{}   c   @{}|@{}}
\hline
{\mk{*}}      &\mk{*}          & \mk{*}&     \mk{*}    \\
\hline
{\mk{ }}      &\mk{*}          & \mk{*}&     \mk{*}    \\
\hline
{\mk{ }}      &\mk{ }          & \mk{*}&     \mk{ }    \\
\hline
{\mk{ }}      &\mk{ }          & \mk{ }&     \mk{*}    \\
\hline
\end{array}
&
\begin{array}{@{}c@{}}
\tmk{\scriptstyle{{i_1}}}\\\tmk{\scriptstyle{{i_2}}}\\ 
\tmk{\scriptstyle{{i_3}}}\\\tmk{\scriptstyle{{i_4}}}
\end{array}
\\
_{B^+(i_2)}\phantom{a}
\end{array}
\;\;
\begin{array}{@{}c@{}l@{}}
\begin{array}{  @{}c @{}c @{}c @{}c @{}}
\mk{_{i_1}}&\mk{_{i_2}}&\mk{_{i_3}}&\mk{_{i_4}}
\end{array}
\\
\begin{array}{ @{}|@{}c   @{}|@{} c @{}|@{} c   @{}|@{}   c   @{}|@{}}
\hline
{\mk{*}}      &\mk{ }          & \mk{*}&     \mk{*}    \\
\hline
{\mk{ }}      &\mk{*}          & \mk{*}&     \mk{*}    \\
\hline
{\mk{ }}      &\mk{ }          & \mk{*}&     \mk{*}    \\
\hline
{\mk{ }}      &\mk{ }          & \mk{ }&     \mk{*}    \\
\hline
\end{array}
&
\begin{array}{@{}l@{}}
\tmk{\scriptstyle{{i_1}}}\\\tmk{\scriptstyle{{i_2}}}\\ 
\tmk{\scriptstyle{{i_3}}}\\\tmk{\scriptstyle{{i_4}}}
\end{array}
\\
_{B^-(i_3).}\phantom{a}
\end{array}
 \end{align}
As usually, the matrix problems are $M\to SMS^{-1}$ 
modulo empty spaces and the transformation $S$ has
the form transposed to $M.$


\section{Matrix problems}
\label{sec_boxes}
\noindent



\noindent
In this section we use the technique of boxes 
and follow the notations of \cite{boddr2}. 
From now on let $\kA$ be a Roiter box 
and $(Q,\d)$ its differential biquiver, 
where  $Q=(I,Q_0,Q_1)$ with  the set of vertices  $I$ and 
the sets of solid and dotted arrows 
respectively $Q_0$ and $Q_1$. 
Let $\kA$-$\cmod$ be the category of finite dimensional 
$\kA$-modules and $\Br_{\kA}$ its full subcategory of
 bricks.
For details 
concerning boxes we also refer to \cite{Dro01} and \cite{thesis}.
Summarizing previous sections we  conclude that  
our approach provides a full and dense functor
$ \VB_E\stackrel{\sim}\lar \Tr_E \lar \MP_E$
and the primary reduction is an equivalence of categories
$ {\MP^s_E(\rr)\stackrel{\sim}\lar \Br_{\kA}(\sss),}$
for some special box $\kA$ and 
dimension vector $\sss.$ 
The composition of these functors
yields  an equivalence
%
$\VB^s_E(r,\dd)\stackrel{\sim}\lar \Br_{\kA}(\sss),$
where both the box $\kA$ and the tuple $\sss$ 
are uniquely defined by the curve 
$E,$ the rank $r$ and the multidegree $\dd.$

In most situations it is useful to present a representation $\M$
as a block-matrix with the block $\M(x)$ 
on the place $(ij)$ for $x\in Q_0(j,i).$
As in the previous sections, with a little abuse of 
notations, we write the matrices $\M$ and $S$ in a form of a table
with $x$ on the $ij$-entry instead of $\M(x).$  
In accordance with  Section \ref{sec_prim_red}
we denote an identity block and a zero block by ``$\one$'' and ``0'' 
respectively. Thus adjust our former notations to that of the theory of boxes.

\subsection*{Class of $\BC$-boxes}
A box $\kA$ with the differential biquiver $(Q,\d)$ is
  of \emph{$\BC$-type }
if its solid arrows form an  $I\times I$  matrix. 
There are two total orders on the set $I:$ 
a row order denoted by $<_r$ and  
a column order denoted by $<_c .$ 
The set of dotted arrows $Q_1$ consists of two subsets:
$\{u\in Q_1(k,j) | j>_r k \}$
and 
$\{v\in Q_1(i,l) | l >_c i \}.$
For each $x\in Q_0(i,j),$ the differential is
$$
\d(x) = \Sum_{l<_c i}x'v
-\Sum_{ j<_r k }u x'',
$$ 
where $x'\in Q_0(l,j) $ and $x''\in Q_0( i,k)$
are uniquely defined as the entries $(jl)$ and $(ki)$ of the matrix $I\times I.$
Such boxes can be presented via matrices $\M$  and $(S,S')$
and matrix multiplications:
$M\to SM(S')^{-1}, $ where
\begin{align*}
M =
\begin{array}{c@{}}
\begin{array}{@{}c @{}  c  @{}  c@{}}
\cmk{{_{c_1}}}& \cmk{\dots} & \cmk{{_{c_n}}}
\end{array}
\\
\begin{array}{|@{}c @{}| @{}  c  |@{}  c @{}|@{}}
\hline
\mmk{x_{r_1c_1}}&\mmk{\dots}& \mmk{x_{r_1c_n}}\\
\hline
\mmk{\vdots}&\mmk{\ddots} & \mmk{\vdots}\\
\hline
\mmk{x_{r_n c_1}}&\mmk{\dots} & \mmk{x_{r_nc_n}}\\
\hline
\end{array}
\\ \phantom{A}
\end{array}
\begin{array}{@{}c@{}}
\cmk{}\\
 \rmk{_{r_1}}\\ \rmk{\vdots} \\ \rmk{_{r_n}}
\\ \phantom{A}
\end{array}
\;\;
~~(S,S')=
\left(
\begin{array}{|@{}c @{}| @{}  c  |@{}  c@{}|}
\hline
\mmk{w_{r_1}}&\mmk{0}& \mmk{0}\\
\hline
\mmk{\vdots }&\mmk{\ddots} & \mmk{0}\\
\hline
\mmk{u_{r_nr_1}}&\mmk{\dots} & \mmk{w_{r_n}}\\
\hline
\end{array}
\,\,
\begin{array}{|@{}c @{}|@{}  c  |@{}  c@{}|}
\hline
\mmk{w_{c_1}}&\mmk{0}& \mmk{0}\\
\hline
\mmk{\vdots}&\mmk{\ddots} & \mmk{0}\\
\hline
\mmk{v_{c_nc_1}}&\mmk{\dots} & \mmk{w_{c_n}}\\
\hline
\end{array}
\right)
\end{align*}
and $(r_1\dots r_n)$ and  $(c_1,\dots c_n)$  are
orders $<_r$ and $<_c$ on $I,$ i.e.
$r_1<_r r_2<_r \dots <_r r_n$ and 
$c_1<_c c_2<_c \dots <_c c_n.$
 The reduced matrix problem for a nodal curve 
from Subsection \ref{subsec_prim_node}
as well 
as all the problems  $A^+,$ $A^-$ and $C$ 
from Subsections \ref{subsec_prim_cycle_2}
and \ref{subsec_prim_cycle3}
are of $\BC$-type.
Note that $\BC$-matrix problems are examples of \emph{bunches of chains}.
\medskip

\subsection*{Class of $\BT$-boxes }
A box $\kA$ with the differential biquiver  $(Q,\d)$ 
is of  \emph{$\BT$-type }
if there exists a set of \emph{distinguished} loops:
$\aa:=\{a_i\in Q_0(i,i)|i\in I\},$
an injective map:
 $v:Q_0\setminus\aa \hookrightarrow Q_1,$
mapping a solid arrow $a:i\mo j$ 
to an opposite directed  dotted arrow  
${\xymatrix @ -1pc { v_a:=v(a):j\ar@{..2>}[r]& i }},$
and for each distinguished loop $a_i\in \aa$ we have

\begin{align}
\label{eq_dist}
 \d a_i =
 \Sum_{c:\,\cdot \mo i}c\cdot v_c
-\Sum_{d:\,i\mo \cdot }v_d\cdot d.
\end{align}

\noindent 
The class of $\BT$-boxes was studied 
in details in \cite{boddr2}. The main property is that
a connected $\BT$-box with more than one vertex is wild but brick-tame. 
However, here we do not use any theoretical results. Our arguments 
are based on the concrete calculations for $\BT$-boxes 
with at most four vertices.
For a box $\kA$ of $\BT$-type 
its biquiver  $Q$ can be encoded as follows:
a vertex $i\in I$ is denoted by a bullet $\bullet;$
on the set of vertices we draw the  graph 
with arrows $Q_0\setminus \aa .$
Such system of notations becomes quite useful
since in most of our cases it is clear
how to recover the differential.

\medskip
  The $\BT$-box $\kA$ obtained in Subsection \ref{subsec_prim_cusp} 
for a cuspidal cubic curve is  
   ${\xymatrix  @ -0.1pc {^1\bullet  & \bullet^2 \ar@<+1pt>[l]}}.$
The problems on three vertices $A^+$ and $A^-$
from Subsection \ref{subsec_prim_tacnode}
  and the problems 
on four vertices:
$A^+,$  $A^-,$ $B^+(j)$ and $B^-(i)$
from Subsection \ref{subsec_prim_IV}
are also of $\BT$-type:
\begin{align*}
\begin{array}{cccccc}
\xy ;/r0.2pc/:
\POS(10,0)*+{\bullet}="v1",
\POS(25,0)*+{\bullet}="v3",
\POS(17,11)*+{\bullet}="v2",
{\ar"v3";"v1" },
{\ar"v2";"v1" },
\POS(7,0)*{_i}
\endxy
&
\xy ;/r0.2pc/:
\POS(10,0)*+{\bullet}="v1",
\POS(25,0)*+{\bullet}="v3",
\POS(17,11)*+{\bullet}="v2",
{\ar"v3";"v1" },
{\ar"v3";"v2" },
\POS(27,0)*{_j}
\endxy
&
    \xy ;/r0.2pc/:
\POS(0,0)*+{\bullet}="v1"; 
\POS(0,15)*+{\bullet}="v2";
\POS(15,15)*+{\bullet}="v3";
\POS(15,0)*+{\bullet}="v4";
\POS(-3,0)*{_i};
{\ar"v2";"v1" },
{\ar"v4";"v1" },
{\ar"v3";"v1" },
 \endxy
 &
\xy ;/r0.2pc/:
\POS(0,0)*+{\bullet}="v1"; 
\POS(0,15)*+{\bullet}="v2";
\POS(15,15)*+{\bullet}="v3";
\POS(15,0)*+{\bullet}="v4";
\POS(17,17)*{_j};
{\ar"v3";"v1" },
{\ar"v3";"v2" },
{\ar"v3";"v4" },
 \endxy
 &
\xy ;/r0.2pc/:
\POS(0,0)*+{\bullet}="v1"; 
\POS(0,15)*+{\bullet}="v2";
\POS(15,15)*+{\bullet}="v3";
\POS(15,0)*+{\bullet}="v4";
\POS(-3,0)*{_i};
\POS(17,17)*{_j};
{\ar"v2";"v1" },
{\ar"v2";"v3" },
{\ar"v4";"v3" },
{\ar"v4";"v1" },
{\ar"v3";"v1" },
 \endxy
 &
\xy ;/r0.2pc/:
\POS(0,0)*+{\bullet}="v1"; 
\POS(0,15)*+{\bullet}="v2";
\POS(15,15)*+{\bullet}="v3";
\POS(15,0)*+{\bullet}="v4";
\POS(-3,0)*{_i};
\POS(17,17)*{_j};
{\ar"v1";"v2" },
{\ar"v3";"v2" },
{\ar"v3";"v4" },
{\ar"v1";"v4" },
{\ar"v1";"v3" },
 \endxy
 \\
 {\scriptstyle A^+}&
{\scriptstyle A^-}&
{\scriptstyle A^+}&
{\scriptstyle A^-}&
{\scriptstyle B^+(j)}&
{\scriptstyle B^-(i)}
 \end{array}
\end{align*}

\begin{remk}
\label{remk_poset}
The listed $\BT$-boxes and that
which appear in the following sections
determine  partially ordered sets $(I,\prec),$ 
by the rule $i\prec j$ if there exists $x\in Q_0(j,i).$
In most of our cases a poset defines a box,
however in general, it does not provide enough information
to recover the differential. 
On the other hand, a pair of linear orders $<_r$ and $<_c$ 
in the definition of a $\BC$-box determine
a partial order 
$\prec$ by the rule $i\prec j$  if $i<_r j$ and $i<_c j .$
Posets obtained in such a way relay 
$\BC$ and $\BT$-boxes.  Moreover, 
for the $\BT$-box
they determine  the canonical minimal edge $(ij),$
where  $i$ is the minimal with respect to the total order $<_r$ and
$j$ is the maximal with respect to $<_c.$ 
Thus having a fixed dimension vector $\sss,$
not only for a $\BC$-box but also for the corresponding $\BT$-box
we have the \emph{canonical} course of reduction.
\end{remk}


\subsection*{Bricks and small reduction}
Boxes of $\BC$ and $\BT$-types possess 
a common property.
The following proposition allows to 
replace the 
usual matrix reduction by the small one.

\begin{prop}
\label{prop_sm}
Let $\kA$ be a box of $\BC$  or $\BT$ type,
$b:i\mo j$ its minimal edge and $\M$ a brick.   Then 
$M(b)$ has maximal rank.
\end{prop}

\begin{proof}
  Let $\kA$  be a box of $\BC$-type.
Since $\kA$ is an example of bunches of chains, we can   
assume that $\M$ is reduced to its canonical form.
Also assume that
$
M(b)=
\left(
\begin{smallmatrix}
  0&0\\
  \Id&0
\end{smallmatrix}
\right).
$
Let rows and columns of $\M$ be ordered $1,\dots, R.$ 
For a place $t\in\{1,\dots, R\}$ by $r(t)$ and $c(t)$ we denote
the row-block and the column-block containing $t$.
For example, since rows and columns are ordered, 
we have $r(1)=j$ and $c(R)=i.$
If $\M$ is invertible then there exist places $m$  
and $n $  such that $M_{1m}=M_{nR}=1$ and 
all the other entries in the first row and the last ($R$-th) column are 
zero. 
A nonscalar endomorphism $(S,S')$ of $\M$ can be constructed
by taking nonzero $S_{n1}=-S'_{Rm}$, diagonal entries 
to be, for example, 1   
and all the other non-diagonal entries to be zero.
Since $c(m)<_c i$ and $r(n)>_r j$ 
the block $S_{r(n)r(1)}$ containing the entry  $S_{n1}$ 
and the block $S'_{c(R)c(m)}$ containing the entry $S'_{Rm}$ are nonempty.

If $\kA$  is a box of $\BT$-type 
then after a step of minimal edge reduction there is a dotted arrow
which is not involved in any differential and hence
there is a nonscalar endomorphism
(for details see \cite[Lemma 3.1]{boddr2}).
\end{proof}

\medskip
\subsection{Small reduction automaton}
\label{subsec_autom}
Recall that an automaton is
an oriented graph on the set of vertices called \emph{states}, 
whose arrows 
are \emph{transitions} from a state to a state.
In our case the states are the matrix problems 
and the transitions encode canonical
steps of reduction.
 
 \begin{defin}
\label{def_aut}
A \emph{small-reduction automaton}
is an oriented graph $\Gamma,$ where
\begin{itemize}
  \item the set of states $\Gamma_0$ is finite and consists
  of boxes $\kA,$ whose differential biquivers have   
  the same finite set of vertices $I.$
  
\item 
The set of transitions $\Gamma_1$ is a subset of $ I\times I.$

\item  
For a minimal solid arrow either $j\mo i $ or $i\mo j$
the transition $(ij): \kA \mo \kA'$
acts on the space of sizes $\mN^{|I|}$ as :
$\sss\to \sss',$ where  $s_k'=s_k$ for $k\neq i$ and
$s_i\to s_i-s_j,$ provided $s_i>s_j.$
\end{itemize}
\end{defin}

\noindent
  A sequence
 $p:=(i_nj_n)\dots(i_2j_2)(i_1j_1)$ 
of transitions is called a \emph{path} 
if the target of $(i_kj_k)$ coincides with the source
of $(i_{k+1} j_{k+1}).$ 
A path operates on the set of sizes:
 $p:\sss \to \sss',$
 where $\sss \geq \sss'$ i.e. $s_i\geq s_i'$  for all $i\in I.$ 
 Two paths $p_1$ and $p_2$ with a common source 
 and a common target are called 
 \emph{equivalent } if for any tuple of
 sizes $\sss\in\mN^{I}$ we have
 $p_1(\sss)=p_2(\sss).$
The semigroup of paths modulo the equivalence relation 
is called the \emph{semigroup of the automaton}.

\subsection*{Principal states}
Let $\Gamma$ be an automaton of small reduction starting
from one of the boxes from Tables \ref{table_I_2}--\ref{table_IV}.     
A state $\kA\in\Gamma_0$ is called \emph{principal}
if it can be interpreted 
in terms of vector bundles 
$\Br_{\kA}(\sss)\cong \VB^s_E(r,\dd).$ 
We show that for a rank $r$ and a
multidegree $\dd$ such that $\gcd(r,d)=1$
there exists a path $p$ on $\Gamma$ 
connecting principal states $\kA$ and $\kA'$ such that 
\begin{align*}
{\xymatrix{
\VB_E^s(r,\dd)\ar[d]^{\cong}
&&& \PIC^{(0,\dots,0)}(E)\ar[d]^{\cong}
\\
\Br_{\kA}(\sss) \ar[rrr]^{p}_{\sim}
&&& \Br_{\kA'}(1,0,\dots,0).
}}
\end{align*}

\noindent
Then a canonical form of a simple vector bundle 
 can be constructed as follows.
\subsection{Algorithm.}
\label{algorithm}
Let $E$ be a reduced plane degeneration of an elliptic curve with $N$
components,
 $(r,\dd)\in\mN\times\mZ^N$ be a tuple of  integers, 
such that $\gcd(r,d)=1;$
where $d= \sum_{k=1}^N d_k$ 
and let $\lambda\in \ko$ be a continuous parameter.
\begin{enumerate}

\item  Use one of the Tables
\ref{table_I_2}, \ref{table_I_3} or \ref{table_IV} 
(with respect to the type of $E$)
to recover
the matrix problem $\Br_{\kA}$ 
and the dimension vector $\sss\in\mN^{N+1}$
from $(r,\dd).$

\item Take $\Br_{\kA} (\sss)$
 as the input data for the corresponding small-reduction
automaton. 
Choose a path $p$ on it such that $p(\sss)=(1,0,\dots,0).$

\item 
Starting with the one-dimensional matrix 
$\lambda \in \Br_{\ko[t]}(1)$ 
reverse the course of reduction 
along the path $p.$
This way, step-by-step recover the canonical form \linebreak
${M(\lambda)=p^{-1}(\lambda)\in\Br_{\kA}(\sss)\cong \VB_E^s(r,\dd).}$ 

\end{enumerate}

\section{Small reduction for nodal and cuspidal cubic curves}
\label{sec_small_red_node_cusp}
The categories 
obtained  in Subsections \ref{subsec_prim_node} and
\ref{subsec_prim_cusp} 
can be interpreted 
as the categories $\kA$-$\cmod(s_1,s_2),$
where $\kA$ are boxes of either $\BC$ and $\BT$-types.
In order to illustrate the language of boxes we present $\kA$ for a nodal curve 
as a differential biquiver, 
despite the agreement to present $\BC$-boxes by tables:
\begin{equation*}
\begin{array}{c}
 {\xymatrix @ -1pc
    {
    &&&\\
1 \ar@(ul,dl)_{a_{1}}\ar@/^20pt/[rrr]^{ c}
\ar@{..>}@/_5pt/[rrr]_{u}\ar@{..>}@/_20pt/[rrr]_{v}
&&& 2 \ar@/_5pt/[lll]_{b}\ar@(ur,dr)^{a_2} 
\\
}}
\\
\begin{array}{l@{}l}
 \d(b)~&=0, \\
 \d(a_1)~&= bu ,\\
 \d(a_2)~&= - vb, \\
 \d(c)~&= - va_1 + a_2u. 
\end{array}
\end{array}
\,\hbox{ and }\,
\begin{array}{c}
 {\xymatrix @ -1pc
    {
    &&&\\
1 \ar@(ul,dl)_{a_{1}}
\ar@{..>}@/_7pt/[rrr]_{v}
&&& 2 \ar@/_7pt/[lll]_{b}\ar@(ur,dr)^{a_2} 
\\
}}
\\
\begin{array}{l@{}l} 
\\
\d(b)~&=0,\\
 \d(a_1)~&= bv, \\
 \d(a_2)~&= - vb. 
\end{array}
\\
\end{array}
\end{equation*}
In both cases the steps of small reduction 
are 
$
\kA\stackrel{(12),(21)}\Arr \kA.
$
In other words,  both problems are \emph{self-reproducing}, and the
small-reduction automaton is
\begin{align}
{
\xy 
\POS(0,0)*\cir<5pt>{}="a"
+\ar@(lu,ld)_{(21)},
+\ar@(ru,rd)^{(12)}
\endxy .
}
\end{align}  
The transitions act on sizes as 
$(12):(s_1,s_2)\to (s_1-s_2,s_2)$
if $s_1\geq s_2$
and 
$(21):(s_1,s_2)\to (s_1,s_2-s_1)$
if $s_1<s_2.$
In terms or rank and degree we get
\begin{align}
\label{rkdg}
(12):\VB_E^s(r,\bd)\mo \VB_E^s(r-\bd,\bd)\, \hbox{ and }
(21):&\VB_E^s(r,\bd)\mo \VB_E^s(\bd,2\bd-r).
\end{align}
That implies the statement of Theorem \ref{theo_main} 
for irreducible cubic curves.

\begin{remk}
The semigroup of paths $\langle (21),(12)\rangle$,  
   generates an  $\SL(2,\mZ)$-action on the set of discrete parameters 
   $(s_1,s_2).$
\end{remk}

\begin{remk}
\label{remk_FM}
   The group generated by Seidel-Thomas spherical twists
   $\langle \T_{\oo}, \T_{\ko(p_0)} \rangle$  
   acts as $\SL(2,\mZ)$ on the $K$-group, or what is
   equivalent, on rank and degree $(r,d)\in \Hom_{\mZ}(K_0(E),\mZ).$ 
   For singular Weierstra\ss{} curves
   the action of the reduction automaton 
   on discrete parameters $(r,\bd)$  can be interpreted in terms of
   Fourier-Mukai transforms, 
   namely: $(12)$ acts as $\T_{\oo}$ 
   and $(12)$ acts as $\T_{\ko(p_0)}^2 \F,$     
   where $$\F= \T_{\ko(p_0)}\T_{\oo}\T_{\ko(p_0)}[-1].$$
   For a singular Weierstra\ss{} curve $E$ in \cite[Corollary 5.11]{BurbanDrozd2}
    Burban and Drozd constructed a fully-faithful functor 
   $\Perf(E)\emb \D^b(B),$ where $B$ is the so called 
   \emph{Butler-Burt algebra} associated to the matrix problem $\MP_E$  
   as introduced in \cite{BuBu}.
   The combinatorics of bricks over $B$ is the same as  that over $\kA.$
   There is a strong evidence that small-reductions (12) and (21) are shadows 
   of some derived autoequivalences of the derived category $\D^b(B).$ 
 \end{remk}

   \begin{remk}
Let us stress that for reductions (\ref{rkdg}) we take 
$\VB_E^s(r,\bd)$ with $0\leq \bd< r,$ or in other words, 
the full subcategory of vector bundles $\ee\in \VB^s_E$
 with the normalization \linebreak
${\en= \oo_{\mP^1}^{r-{\bd}}\oplus \oo_{\mP^1}^{\bd}(1).}$
If the degree $d$ is arbitrary, we identify $\VB_E^s(r,d)$ with $\VB_E^s(r,\bd),$
where $\bd=d \bmod r,$ using the Picard group action.
For instance, to proceed with the reduction 
after (12)-step, we have to replace
the category $\VB_E^s(r-\bd,\bd)$
by
${\VB_E^s(r-\bd, \bd \bmod(r-\bd))}.$ 
 \end{remk}



%% file: 2comp.tex
\section{Small reduction for Kodaira fibers I$_2$ and III.}
\label{sec_box_2cycle}

In Subsections \ref{subsec_prim_cycle_2} 
and \ref{subsec_prim_tacnode}
we obtained an equivalence
$\MP^s(r,\dd)\stackrel{\cong}\lar\Br_\kA(\sss),$
where the box $\kA$ is 
the configuration $\A^{\sigma},$
of $\BC$  or $\BT$-type, $\sigma\in\{+,-\} $
depending on whether $r>\bd$ or $r<\bd.$
Applying small reduction to the box $\A^{\sigma}$ we obtain 
another configuration on 3 blocks, 
defined by the standard numeration of blocks (1,2,3).
We denote this type of boxes by $\B.$
 In the $\BC$-case we get:
$$
M= 
\begin{array}{@{}c@{}l@{}}
\begin{array}{@{}c@{}c@{}c@{}}
\smk{_1}&\smk{_2}&\smk{_3}
\end{array}
&
\\
\begin{array}{ @{}|@{}c @{}|@{}  c  @{}|@{}    c   @{}|@{}}
\hline
{\mk{a_1}}      &\mk{*}         &     \mk{*}    \\
\hline
 {\mk{*}}     &{\mk{a_2}}    &     \mk{*}      \\
\hline
 {\mk{*}}     &{\mk{*}}    &     \mk{a_3}    \\
\hline
\end{array}
&
\begin{array}{@{}r@{}}
\tmk{\scriptstyle{1}}\\\tmk{\scriptstyle{2}}\\\tmk{\scriptstyle{3}}
\end{array}
\\ \phantom{A}
\end{array}
\;\;\;
(S,S')=
\left(
\begin{array}{@{}l@{}c@{}}
&
\begin{array}{@{}c@{}c@{}c@{}}
\smk{_1}&\smk{_2}&\smk{_3}
\end{array}
\\
\begin{array}{@{}r@{}}
\tmk{\scriptstyle{1}}\\\tmk{\scriptstyle{2}}\\\tmk{\scriptstyle{3}}
\end{array}
&
\begin{array}{ @{}|@{}c  @{}|@{} c  @{}|@{}    c   @{}|@{}}
\hline
{\mk{w_1}}      &\mk{0}         &     \mk{0}    \\
\hline
 {\mk{u_3}}     &{\mk{w_2}}    &     \mk{0}      \\
\hline
 {\mk{u_2}}     &{\mk{u_1}}    &     \mk{w_3}    \\
\hline
\end{array}
\\ \phantom{A}
\end{array}
\;
\begin{array}{@{}c@{}l@{}}
\begin{array}{@{}c@{}c@{}c@{}}
\smk{_1}&\smk{_2}&\smk{_3}
\end{array}
&
\\
\begin{array}{ @{}|@{}c @{}|@{} c  @{}|@{}    c   @{}|@{}}
\hline
{\mk{w_1}}      &\mk{0}         &     \mk{0}    \\
\hline
 {\mk{v_3}}     &{\mk{w_2}}    &     \mk{0}      \\
\hline
 {\mk{v_2}}     &{\mk{v_1}}    &     \mk{w_3}    \\
\hline
\end{array}
&
\begin{array}{@{}l@{}}
\tmk{\scriptstyle{1}}\\\tmk{\scriptstyle{2}}\\\tmk{\scriptstyle{3}}
\end{array}
\\ \phantom{A}
\end{array}
\right).
$$


\noindent
As was mentioned in Remark \ref{remk_poset}
column and row-orders define a poset.
Configurations $\A^{+}, \A^{-}$  and $\B$ determine respectively the posets
\begin{align*}
\begin{array}{c}
\xy
\POS(10,0)*+{\bullet}="v1",\POS(25,0)*+{\bullet}="v3",\POS(17,11)*+{\bullet}="v2",
{\ar"v3";"v1" <+0.1pc>},{\ar"v2";"v1" <-0.1pc>},
\POS(9,2)*{^1}
\POS(27,2)*{^3}; \POS(19,12)*{^2}.
\endxy
\\
{\scriptstyle A^+}
\end{array}
~~~
\begin{array}{c}
\xy
\POS(10,0)*+{\bullet}="v1",\POS(25,0)*+{\bullet}="v3",
\POS(17,11)*+{\bullet}="v2",
{\ar"v3";"v1" <+0.1pc>},{\ar"v3";"v2" <-0.1pc>}
\POS(9,2)*{^1}
\POS(27,2)*{^3}; \POS(19,12)*{^2}.
\endxy
\\
{\scriptstyle A^-}
\end{array}
~~\hbox{ and }~~
\begin{array}{c}
\xy
\POS(10,0)*+{\bullet}="v1",
\POS(25,0)*+{\bullet}="v3",
\POS(17,11)*+{\bullet}="v2",
{\ar"v3";"v1" <+0.1pc>},{\ar"v3";"v2" <-0.1pc>},{\ar"v2";"v1" };
\POS(9,2)*{^1}
\POS(27,2)*{^3}; \POS(19,12)*{^2}.
\endxy
\\
{\scriptstyle B}
\end{array}
\end{align*}

\noindent
Lets illustrate on an example how 
to associate a $\BT$-differential biquiver to a poset.
For $A^-$ and $B$  we have respectively:
\begin{align}
\label{diff_biq_B}
\begin{array}{@{}l@{}l@{}l@{}l}
\begin{array}{@{}l@{} }
 {\xymatrix  @ -0.7pc
    {
&&& && \\
&&2 
\ar@(ul,ur)^{a_{2}}
\ar@{..>}@/^+10pt/[dr]^{v_c}
&&\\
&1
\ar@(lu,ld)_{a_1}
\ar@{..>}@/_+10pt/[rr]_{v_{b}}
&&
3\ar@(ru,rd)^{a_3} 
\ar[ll]_{b}
\ar[ul]^{c}
&\\
}}
\end{array}
&
\begin{array}{@{}l@{}}
\d(b)=\d(c)=0, \\
\d(a_3)= -v_{b}b- v_c c\\
\d(a_1)=b v_{b}, \\ 
\d(a_2)=c v_{c},\\
\end{array}
\end{array}
\end{align}
\begin{align}
\begin{array}{cc}
\begin{array}{@{}l@{}}
 {\xymatrix  @ -0.7pc
    {
&&& && \\
&&2 
\ar@(ul,ur)^{a_{2}}
\ar@{..>}@/^+10pt/[dr]^{v_c}
\ar[dl]^{a}
&&\\
&1
\ar@(lu,ld)_{a_1}
\ar@{..>}@/^+10pt/[ur]^{v_a}
\ar@{..>}@/_+10pt/[rr]_{v_{b}}
&&
3\ar@(ru,rd)^{a_3} 
\ar[ll]_{b}
\ar[ul]^{c}
&\\
}}
\end{array}
&
\begin{array}{@{}l@{}}
  \d(b)=0, \\
  \d(a)=bv_c,  \\
  \d(c)=-v_a b, \\
\d(a_1)=b v_{b}+a v_{a}, \\ 
\d(a_2)=c v_{c}-v_{a}a,\\
\d(a_3)= -v_{b}b- v_c c.
\end{array}
\end{array}
\end{align}

\noindent
In Subsection \ref{subsec_prim_tacnode}
we obtained an equivalence
$\MP^s(r,d_1,d_2)\stackrel{\cong}\lar\Br_\kA(s_1,s_2,s_3),$
where $\kA$ was a $\BT$-box of type either $A^+$ or $A^-.$
The small reduction automaton starting at, let us say $A^+,$ 
is 
\begin{align}
\label{autom_cycle2}
{\xymatrix {
  _{\A^{+}}  \ar@/^10pt/[rr]^{(12)} 
  \ar@(ul,ld)_{(21)} 
  &&
  _{\B} \ar@/^10pt/[rr]^{(31)} \ar@/^10pt/[ll]^{(13)} &&
  _{\A^{-}}  \ar@/^10pt/[ll]^{(32)} \ar@(ur,rd)^{(23)}
  \\
  &&&&&
}}
\end{align}

\noindent
This is the small reduction automaton 
for a cycle of two lines, which is also the canonical one for a tacnode curve. 
We claim that the reduction can terminate only at the
states $\A^+$ and $\A^-,$ which are principal.
Indeed, assume that we have the box $\B$ with sizes $s_1=s_3.$ 
Then the matrix can be reduced to the canonical form:
\vspace{-0.3cm}
$$ 
\begin{array}{@{}c@{}c@{}}
\begin{array}{@{}c@{}c@{}c}
\mk{_1}&\mk{_2}&\mk{_3}
\end{array}
&
\\
\begin{array}{ |@{}c @{}|@{} c  @{}|@{}    c   @{}|@{}}
\hline
{\mk{0}}      &\mk{0}         &     \mk{\one}    \\
\hline
 {\mk{0}}     &{\mk{J_1}}    &     \mk{0}      \\
\hline
 {\mk{J_2}}     &{\mk{0}}    &     \mk{0}    \\
\hline
\end{array}
&
\begin{array}{@{}l@{}}
\tmk{\scriptstyle{1}}\\\tmk{\scriptstyle{2}}\\\tmk{\scriptstyle{3}}
\end{array}
\end{array}
$$
where $J_1$ and $J_2$ are Jordan cells with nonzero eigenvalues. 
It is quite obvious that this matrix is decomposable.
Analogously in the case of Kodaira fiber III:
the reduction can terminate only at a state of type $\A.$ 
Indeed, if $s_1=s_3$ then the configuration $B$ produces a splitting;
and for $A^+$ 
we get the problem $\Br_{\kA}(s_1,s_2),$
where $\kA$ is the box 
as for a cuspidal cubic curve with sizes $(s_1,s_2,s_3)\to (s_1,s_2):$ 
\begin{align*}
\begin{array}{c}
\xy
\POS(10,0)*+{\bullet}="v1",\POS(25,0)*+{\bullet}="v3",\POS(17,11)*+{\bullet}="v2",
{\ar^{b}"v3";"v1" <+0.1pc>},{\ar_{c}"v3";"v2" <-0.1pc>},
{\ar_{a}"v2";"v1" }
\POS(9,2)*{^1};
\POS(27,2)*{^3}; \POS(19,12)*{^2}.
 \endxy 
 \end{array}
 \stackrel{(31),(13)}\Arr  ~~ \bullet^2 ~~\bullet^3
\;\;\hbox{ and }\;
\begin{array}{c}
\xy
\POS(10,0)*+{\bullet}="v1",\POS(25,0)*+{\bullet}="v3",\POS(17,11)*+{\bullet}="v2",
{\ar^{b}"v3";"v1" <+0.1pc>},,{\ar_{a}"v2";"v1" <-0.1pc>}
\POS(9,2)*{^1};
\POS(27,2)*{^3}; \POS(19,12)*{^2}.
\endxy
\end{array}
 \stackrel{(31)}\Arr
\begin{array}{c}
\xy
\POS(10,0)*+{\bullet}="v1",
\POS(25,0)*+{\bullet}="v3",
{\ar_{a}"v3";"v1" },
\POS(9,2)*{^1};
\POS(27,2)*{^2}; 
\endxy
\end{array}
\end{align*} 

\noindent
By gluing paths we can construct the automaton on principal states:
\begin{equation}
\label{autom_I2}
{\xymatrix{
  _{\A^{+}}  \ar@/^10pt/[rr]^{(31)(12)} \ar@(ul,l)_{(21)}
        \ar@(dl,l)^{(13)(12)}
  &&
  _{\A^{-}}  \ar@/^10pt/[ll]^{(13)(32)} 
    \ar@(ur,r)^{(23)}
      \ar@(dr,r)_{(31)(32)}
  \\
  &&&&&
}}
\end{equation}
\noindent
For a principal configuration $A^{\sigma}$ we introduce its 
new discrete parameters $(\alpha,\beta).$
For  $A^+$ let ${(\alpha, \beta):= (s_1,s_2+s_3)}$ 
and  $(\alpha,\beta):= (s_1+s_2,s_3)$  for $A^-.$

\begin{lemma}
\label{lemma_autom_cycle2}
  Let $p: \A^{\sigma}\mo \A^{\sigma'}$ 
be a path on the principal automaton (\ref{autom_I2})
taking $\sss\to \sss'$ and respectively
$(\alpha,\beta)\mo(\alpha',\beta').$
 Then
$
\gcd(\alpha,\beta)=\gcd(\alpha',\beta').
$
\end{lemma} 

\begin{proof}
It is sufficient to prove the statement on the following transitions:
\linebreak
$
(23), ~(32)(31): \A^-\mo \A^- $
and 
  $~(13)(32): \A^-\mo \A^+.$
Indeed, we have 
\begin{itemize}
  \item[] $(23): (s_1,s_2,s_3)\to (s_1,s_2-s_3,s_3)$ and hence
  $(\alpha,\beta)\to (\alpha-\beta,\beta);$
    \item[] $(32)(31): (s_1,s_2,s_3)\to (s_1,s_2,s_3-(s_1+s_2))$ and
  $(\alpha,\beta)\to (\alpha,\beta-\alpha);$
    \item[] $(13)(32): (s_1,s_2,s_3)\to (s_1+s_2-s_3,s_2,s_3-s_2)$ and
  $(\alpha,\beta)\to (\alpha-\beta,\beta).$
\end{itemize}
\end{proof}

\noindent
Let
$\VB^s_E(r,\dd)\stackrel{\cong}\lar\VB^s_E(r',\dd')$
be a functorial bijection obtained by the course of small reductions 
along  the path $p.$ 
Replacing 
the dimension vector $\sss$ by the tuple $(r,\dd)$ 
using  Table \ref{table_I_2}
we obtain 
$(\alpha,\beta)=(r-d\bmod r,d\bmod r).$
If $\gcd(r,d)=1,$ 
at the end of reduction we get 
${\VB^s_E(r,\dd)\stackrel{\cong}\lar\PIC^{(0,0)}(E)},$
and there are no bricks otherwise.
Hence, Lemma \ref{lemma_autom_cycle2}
implies Theorem \ref{theo_main} for curves I$_2$ and III.

%% file: 3comp.tex
\section{Small reduction for Kodaira fibers I$_3$ and IV.}
\label{sec_small_red_I_3}
In Subsection \ref{subsec_prim_cycle3}
we obtained some equivalences
$\MP^s(r,\dd)\stackrel{\cong}\lar\Br_\kA(\sss),$ 
where $\sss\in \mN^4$ and $\kA$ is a $\BC$-box of type $A^+,$ $A^-$ or $C.$
To fix the notations we rewrite the configurations for 
the set of vertices $I=\{1,2,3,4\}.$
Then a small reduction automaton starting from 
the configuration $\A^+$ is as follows:
\begin{align}
\label{autom_I3}
\begin{array}{c}
{\xymatrix @ +1pc{
  _{\A^{+}}  \ar[r]^{(12)} 
  \ar@(ul,ld)_{(21)} 
  &
  _{\B^+(2)} \ar@<+3pt>[r]^{(31)} 
          \ar[d]^{(13)} 
  &
  _{C(2,3)} \ar@<+3pt>[l]^{(32)}\ar@<+3pt>[r]^{(23)} 
  &
  _{\B^-(3)} \ar@<+3pt>[l]^{(24)} 
          \ar[d]_{(42)} 
  &
  _{\A^{-}}  \ar[l]^{(43)} 
  \ar@(ur,rd)^{(34)}
  \\ 
 & _D \ar[ul]^{(14)}\ar[dl]_{(41)} & 
 &  _{D_*} \ar[ur]_{(41)}\ar[dr]^{(14)}
\\ 
  _{\A_*^{-}}  \ar[r]^{(42)} 
  \ar@(ul,ld)_{(24)} 
  &
  _{\B^-(2)} \ar@<+3pt>[r]^{(34)} 
          \ar[u]_{(43)} 
  &
  _{C(3,2)} \ar@<+3pt>[l]^{(32)}\ar@<+3pt>[r]^{(23)} 
  &
  _{\B^+(3)} \ar@<+3pt>[l]^{(21)} 
          \ar[u]^{(12)} 
  &
  _{\A_*^{+}}  \ar[l]^{(13)} 
  \ar@(ur,rd)^{(31)}
}}
\end{array}
\end{align}
Let us explain the notations:
configurations of type $A$ are
$$
\begin{array}{@{}c@{}l@{}}
\begin{array}{@{}c@{}c@{}c@{}c@{}}
\smk{_1}&\smk{_4}&\smk{_3}&\smk{_2}
\end{array}
&
\\
\begin{array}{ |@{}c @{}|@{} c  @{}|@{}    c   @{}|@{} c   @{}|@{}}
\hline
{\mk{a_1}}      &\mk{*}         &     \mk{*}    &     \mk{*}\\
\hline
 {\mk{*}}     &{\mk{*}}    &     \mk{*}  &     \mk{a_2}     \\
\hline
{\mk{*}}     &{\mk{*}}    &     \mk{a_3}  &     \mk{*}     \\
\hline
{\mk{*}}     &{\mk{a_4}}    &     \mk{*}  &     \mk{*}     \\
\hline
\end{array}
&
\begin{array}{@{}c@{}}
\tmk{\scriptstyle{1}}\\\tmk{\scriptstyle{2}}\\\tmk{\scriptstyle{3}}
\\\tmk{\scriptstyle{4}}
\end{array}
\\ _{\A^+}\phantom{a}
\end{array}
\;\;\;
\begin{array}{@{}c@{}l@{}}
\begin{array}{@{}c@{}c@{}c@{}c@{}}
\mk{_1}&\mk{_2}&\mk{_3}&\mk{_4}
\end{array}
&
\\
\begin{array}{ |@{}c @{}|@{} c  @{}|@{}    c   @{}|@{} c   @{}|@{}}
\hline
{\mk{*}}      &\mk{*}     &     \mk{a_3}    &     \mk{*}\\
\hline
 {\mk{*}}     &{\mk{a_2}}     &     \mk{*}  &     \mk{*}     \\
\hline
{\mk{a_1}}     &{\mk{*}}    &     \mk{*}&     \mk{*}       \\
\hline
{\mk{*}}     &     \mk{*}   &{\mk{*}}    &     \mk{a_4}   \\
\hline
\end{array}
&
\begin{array}{@{}c@{}}
\tmk{\scriptstyle{3}}\\\tmk{\scriptstyle{2}}\\\tmk{\scriptstyle{1}}
\\\tmk{\scriptstyle{4}}
\end{array}
\\ _{\A^-}\phantom{a}
\end{array}
\;\;\;
\begin{array}{@{}c@{}l@{}}
\begin{array}{@{}c@{}c@{}c@{}c@{}}
\smk{_1}&\smk{_4}&\smk{_2}&\smk{_3}
\end{array}
&
\\
\begin{array}{ |@{}c @{}|@{} c  @{}|@{}    c   @{}|@{} c   @{}|@{}}
\hline
{\mk{a_1}}      &\mk{*}         &     \mk{*}    &     \mk{*}\\
\hline
 {\mk{*}}     &{\mk{*}}    &     \mk{*}  &     \mk{a_3}     \\
\hline
{\mk{*}}     &{\mk{*}}    &     \mk{a_2}  &     \mk{*}     \\
\hline
{\mk{*}}     &{\mk{a_4}}    &     \mk{*}  &     \mk{*}     \\
\hline
\end{array}
&
\begin{array}{@{}c@{}}
\tmk{\scriptstyle{1}}\\\tmk{\scriptstyle{3}}\\\tmk{\scriptstyle{2}}
\\\tmk{\scriptstyle{4}}
\end{array}
\\ _{\A^+_*}\phantom{a}
\end{array}
\;\;\;
\begin{array}{@{}c@{}l@{}}
\begin{array}{@{}c@{}c@{}c@{}c@{}}
\mk{_1}&\mk{_3}&\mk{_2}&\mk{_4}
\end{array}
&
\\
\begin{array}{ |@{}c @{}|@{} c  @{}|@{}    c   @{}|@{} c   @{}|@{}}
\hline
{\mk{*}}      &\mk{*}     &     \mk{a_2}    &     \mk{*}\\
\hline
 {\mk{*}}     &{\mk{a_3}}     &     \mk{*}  &     \mk{*}     \\
\hline
{\mk{a_1}}     &{\mk{*}}    &     \mk{*}&     \mk{*}       \\
\hline
{\mk{*}}     &     \mk{*}   &{\mk{*}}    &     \mk{a_4}   \\
\hline
\end{array}
&
\begin{array}{@{}c@{}}
\tmk{\scriptstyle{2}}\\\tmk{\scriptstyle{3}}\\\tmk{\scriptstyle{1}}
\\\tmk{\scriptstyle{4}}
\end{array};
\\ _{\A^-_*} \phantom{a}
\end{array}
$$
the configurations of type $B$ are
$$
\begin{array}{@{}c@{}l@{}}
\begin{array}{@{}c@{}c@{}c@{}c@{}}
\smk{_1}&\smk{_2}&\smk{_4}&\smk{_3}
\end{array}
&
\\
\begin{array}{ |@{}c @{}|@{} c  @{}|@{}    c   @{}|@{} c   @{}|@{}}
\hline
{\mk{a_1}}      &\mk{*}         &     \mk{*}    &     \mk{*}\\
\hline
 {\mk{*}}     &{\mk{a_2}}    &     \mk{*}  &     \mk{*}     \\
\hline
{\mk{*}}     &{\mk{*}}    &     \mk{*}&     \mk{a_3}       \\
\hline
{\mk{*}}     &     \mk{*}&{\mk{a_4}}    &     \mk{*}   \\
\hline
\end{array}
&
\begin{array}{@{}l@{}}
\tmk{\scriptstyle{1}}\\\tmk{\scriptstyle{2}}\\\tmk{\scriptstyle{3}}
\\\tmk{\scriptstyle{4}}
\end{array}
\\ _{\B^+(2)} 
\end{array}
\;\;\;
\begin{array}{@{}c@{}l@{}}
\begin{array}{@{}c@{}c@{}c@{}c@{}}
\smk{_1}&\smk{_2}&\smk{_3}&\smk{_4}
\end{array}
&
\\
\begin{array}{ |@{}c @{}|@{} c  @{}|@{}    c   @{}|@{} c   @{}|@{}}
\hline
{\mk{*}}      &\mk{a_2}     &     \mk{*}    &     \mk{*}\\
\hline
 {\mk{a_1}}     &{\mk{*}}     &     \mk{*}  &     \mk{*}     \\
\hline
{\mk{*}}     &{\mk{*}}    &     \mk{a_3}&     \mk{*}       \\
\hline
{\mk{*}}     &     \mk{*}   &{\mk{*}}    &     \mk{a_4}   \\
\hline
\end{array}
&
\begin{array}{@{}l@{}}
\tmk{\scriptstyle{2}}\\\tmk{\scriptstyle{1}}\\\tmk{\scriptstyle{3}}
\\\tmk{\scriptstyle{4}}
\end{array}
\\ _{\B^-(3)}
\end{array}
\;\;\;
\begin{array}{@{}c@{}l@{}}
\begin{array}{@{}c@{}c@{}c@{}c@{}}
\smk{_1}&\smk{_3}&\smk{_4}&\smk{_2}
\end{array}
&
\\
\begin{array}{ |@{}c @{}|@{} c  @{}|@{}    c   @{}|@{} c   @{}|@{}}
\hline
{\mk{a_1}}      &\mk{*}         &     \mk{*}    &     \mk{*}\\
\hline
 {\mk{*}}     &{\mk{a_3}}    &     \mk{*}  &     \mk{*}     \\
\hline
{\mk{*}}     &{\mk{*}}    &     \mk{*}&     \mk{a_2}       \\
\hline
{\mk{*}}     &     \mk{*}&{\mk{a_4}}    &     \mk{*}   \\
\hline
\end{array}
&
\begin{array}{@{}l@{}}
\tmk{\scriptstyle{1}}\\\tmk{\scriptstyle{3}}\\\tmk{\scriptstyle{2}}
\\\tmk{\scriptstyle{4}}
\end{array}
\\ _{\B^+(3)}
\end{array}
\;\;\;
\begin{array}{@{}c@{}l@{}}
\begin{array}{@{}c@{}c@{}c@{}c@{}}
\smk{_1}&\smk{_3}&\smk{_2}&\smk{_4}
\end{array}
&
\\
\begin{array}{ |@{}c @{}|@{} c  @{}|@{}    c   @{}|@{} c   @{}|@{}}
\hline
{\mk{*}}      &\mk{a_3}     &     \mk{*}    &     \mk{*}\\
\hline
 {\mk{a_1}}     &{\mk{*}}     &     \mk{*}  &     \mk{*}     \\
\hline
{\mk{*}}     &{\mk{*}}    &     \mk{a_2}&     \mk{*}       \\
\hline
{\mk{*}}     &     \mk{*}   &{\mk{*}}    &     \mk{a_4}   \\
\hline
\end{array}
&
\begin{array}{@{}c@{}}
\tmk{\scriptstyle{3}}\\\tmk{\scriptstyle{1}}\\\tmk{\scriptstyle{2}}
\\\tmk{\scriptstyle{4}}
\end{array}
\\ _{\B^-(2)} 
\end{array}
$$
and configurations of types $C$ and $D$ are
$$
\begin{array}{@{}c@{}l@{}}
\begin{array}{@{}c@{}c@{}c@{}c@{}}
\smk{_1}&\smk{_2}&\smk{_4}&\smk{_3}
\end{array}
&
\\
\begin{array}{ |@{}c @{}|@{} c  @{}|@{}    c   @{}|@{} c   @{}|@{}}
\hline
\mk{*}    &{\mk{a_2}}    &     \mk{*}    &     \mk{*}\\
\hline
{\mk{a_1}} &{\mk{*}}      &     \mk{*}  &     \mk{*}     \\
\hline
{\mk{*}}     &{\mk{*}}    &     \mk{*}&     \mk{a_3}       \\
\hline
{\mk{*}}     &     \mk{*}&{\mk{a_4}}    &     \mk{*}   \\
\hline
\end{array}
&
\begin{array}{@{}c@{}}
\tmk{\scriptstyle{2}}\\\tmk{\scriptstyle{1}}\\\tmk{\scriptstyle{3}}
\\\tmk{\scriptstyle{4}}
\end{array}
\\ _{C(2,3)} 
\end{array}
\;\;\;
\begin{array}{@{}c@{}l@{}}
\begin{array}{@{}c@{}c@{}c@{}c@{}}
\smk{_1}&\smk{_3}&\smk{_4}&\smk{_2}
\end{array}
&
\\
\begin{array}{ |@{}c @{}|@{} c  @{}|@{}    c   @{}|@{} c   @{}|@{}}
\hline
\mk{*}    &{\mk{a_3}}    &     \mk{*}    &     \mk{*}\\
\hline
{\mk{a_1}} &{\mk{*}}      &     \mk{*}  &     \mk{*}     \\
\hline
{\mk{*}}     &{\mk{*}}    &     \mk{*}&     \mk{a_2}       \\
\hline
{\mk{*}}     &     \mk{*}&{\mk{a_4}}    &     \mk{*}   \\
\hline
\end{array}
&
\begin{array}{@{}c@{}}
\tmk{\scriptstyle{3}}\\\tmk{\scriptstyle{1}}\\\tmk{\scriptstyle{2}}
\\\tmk{\scriptstyle{4}}
\end{array}
\\ _{C(3,2)} 
\end{array}
\;\;\;
\begin{array}{@{}c@{}l@{}}
\begin{array}{@{}c@{}c@{}c@{}c@{}}
\smk{_1}&\smk{_3}&\smk{_2}&\smk{_4}
\end{array}
&
\\
\begin{array}{ |@{}c @{}|@{} c  @{}|@{}    c   @{}|@{} c   @{}|@{}}
\hline
{\mk{a_1}}      &\mk{*}     &     \mk{*}    &     \mk{*}\\
\hline
 {\mk{*}}     &{\mk{*}}     &     \mk{a_2}  &     \mk{*}     \\
\hline
{\mk{*}}     &{\mk{a_3}}    &     \mk{*}&     \mk{*}       \\
\hline
{\mk{*}}     &     \mk{*}   &{\mk{*}}    &     \mk{a_4}   \\
\hline
\end{array}
&
\begin{array}{@{}l@{}}
\tmk{\scriptstyle{1}}\\\tmk{\scriptstyle{2}}\\\tmk{\scriptstyle{3}}
\\\tmk{\scriptstyle{4}}
\end{array}
\\ _{\D}  
\end{array}
\;\;\;
\begin{array}{@{}c@{}l@{}}
\begin{array}{@{}c@{}c@{}c@{}c@{}}
\smk{_1}&\smk{_2}&\smk{_3}&\smk{_4}
\end{array}
&
\\
\begin{array}{ |@{}c @{}|@{} c  @{}|@{}    c   @{}|@{} c   @{}|@{}}
\hline
{\mk{a_1}}      &\mk{*}     &     \mk{*}    &     \mk{*}\\
\hline
 {\mk{*}}     &{\mk{*}}     &     \mk{a_3}  &     \mk{*}     \\
\hline
{\mk{*}}     &{\mk{a_2}}    &     \mk{*}&     \mk{*}       \\
\hline
{\mk{*}}     &     \mk{*}   &{\mk{*}}    &     \mk{a_4}   \\
\hline
\end{array}
&
\begin{array}{@{}l@{}}
\tmk{\scriptstyle{1}}\\\tmk{\scriptstyle{3}}\\\tmk{\scriptstyle{2}}
\\\tmk{\scriptstyle{4 .}}
\end{array}
\\ _{D_*} 
\end{array}
$$

\medskip

In Subsection \ref{subsec_prim_IV}
we obtained equivalences
${\MP^s(r,\dd)\stackrel{\cong}\lar\Br_\kA(\sss)},$ 
$\sss\in \mN^4,$
where $\kA$ is a $\BC$-box of type $A^+,$ $A^-$ or $B.$
As explained in the Remark \ref{remk_poset} the boxes of $\BC$ 
and $\BT$ types are related. 
Therefore the canonical small reduction automaton for Kodaira fiber IV
can be obtained 
from the automaton (\ref{autom_I3}) by gluing
states $A^{\sigma}$ with $A^{\sigma}_*$ and $D$ with $D_*:$
\begin{align}
\label{autom_IV}
\begin{array}{c}
{\xymatrix @ +1pc{
  &
  _{\B^+(2)} \ar@<+3pt>[r]^{(31)} 
          \ar[dr]_{(13)} 
  &
  _{C(2,3)} \ar@<+3pt>[l]^{(32)}\ar@<+3pt>[r]^{(23)} 
  &
  _{\B^-(3)} \ar@<+3pt>[l]^{(24)} 
          \ar[dl]^{(42)} 
  &
\\
_{\A^{+}}  \ar[rd]_{(13)} \ar[ru]^{(12)} 
  \ar@(l,u)^{(31)} \ar@(l,d)_{(21)}
&&
_D \ar[ll]_{(14)}\ar[rr]^{(41)}
&&
_{\A^{-}}  \ar[lu]_{(43)} \ar[ld]^{(42)} 
  \ar@(r,u)_{(34)}
  \ar@(r,d)^{(24)} 
\\
  &
  _{\B^+(3)} \ar@<+3pt>[r]^{(21)} 
          \ar[ur]^{(12)} 
  &
  _{C(3,2)} \ar@<+3pt>[l]^{(32)}\ar@<+3pt>[r]^{(23)} 
  &
  _{\B^-(2)} \ar@<+3pt>[l]^{(34)} 
          \ar[ul]_{(43)} 
  &
}}
\end{array}
\end{align}

\noindent
For the $\BT$-boxes we have
\begin{align*}
  \begin{array}{ccccc}
    \xy ;/r0.2pc/:
\POS(0,0)*+{\bullet}="v1"; 
\POS(0,15)*+{\bullet}="v2";
\POS(15,15)*+{\bullet}="v3";
\POS(15,0)*+{\bullet}="v4";
\POS(-3,0)*{_1};\POS(-3,17)*{_2};\POS(17,17)*{_3};\POS(18,0)*{_4};
{\ar_{}"v2";"v1" },
{\ar^{}"v4";"v1" },
{\ar^{}"v3";"v1" },
 \endxy
  ~~  ~~
&
    \xy ;/r0.2pc/:
\POS(0,0)*+{\bullet}="v1"; 
\POS(0,15)*+{\bullet}="v2";
\POS(15,15)*+{\bullet}="v3";
\POS(15,0)*+{\bullet}="v4";
\POS(-3,0)*{_1};\POS(-3,17)*{_2};\POS(17,17)*{_3};\POS(18,0)*{_4};
{\ar_{}"v4";"v1" },
{\ar^{}"v4";"v2" },
{\ar^{}"v4";"v3" },
 \endxy
  ~~  ~~
&
 \xy ;/r0.2pc/:
\POS(0,0)*+{\bullet}="v1"; 
\POS(0,15)*+{\bullet}="v2";
\POS(15,15)*+{\bullet}="v3";
\POS(15,0)*+{\bullet}="v4";
\POS(-3,0)*{_1};\POS(-3,17)*{_2};\POS(17,17)*{_3};\POS(18,0)*{_4};
{\ar_{}"v3";"v2" },
{\ar_{}"v4";"v2" },
{\ar^{}"v3";"v1" },
{\ar^{}"v4";"v1" },
{\ar^{}"v4";"v1" },
 \endxy
 ~~  ~~
&
 \xy ;/r0.2pc/:
\POS(0,0)*+{\bullet}="v1"; 
\POS(0,15)*+{\bullet}="v2";
\POS(15,15)*+{\bullet}="v3";
\POS(15,0)*+{\bullet}="v4";
\POS(-3,0)*{_1};\POS(-3,17)*{_2};\POS(17,17)*{_3};\POS(18,0)*{_4};
{\ar_{}"v2";"v3" },
{\ar_{}"v2";"v1" },
{\ar^{}"v4";"v1" },
{\ar^{}"v4";"v3" },
 \endxy
 ~~  ~~
 &
 \xy ;/r0.2pc/:
\POS(0,0)*+{\bullet}="v1"; 
\POS(0,15)*+{\bullet}="v2";
\POS(15,15)*+{\bullet}="v3";
\POS(15,0)*+{\bullet}="v4";
\POS(-3,0)*{_1};\POS(-3,17)*{_2};\POS(17,17)*{_3};\POS(18,0)*{_4};
%
{\ar_{}"v4";"v3" },
{\ar_{}"v2";"v1" },
{\ar_{}"v4";"v2" },
{\ar^{}"v3";"v1" },
{\ar^{}"v4";"v1" },
{\ar^{}"v4";"v1" },
 \endxy
 ~~  ~~
\\
\scriptstyle{A^+=A^+_*} 
&\scriptstyle{A^-=A^-_*} 
& \scriptstyle{C(2,3)} 
& \scriptstyle{C(3,2)} 
& \scriptstyle{ D=D_*}
  \end{array}
\end{align*}
and four configurations of type $B:$
\begin{align*}
  \begin{array}{cccc}
 \xy ;/r0.2pc/:
\POS(0,0)*+{\bullet}="v1"; 
\POS(0,15)*+{\bullet}="v2";
\POS(15,15)*+{\bullet}="v3";
\POS(15,0)*+{\bullet}="v4";
\POS(-3,0)*{_1};\POS(-3,17)*{_2};\POS(17,17)*{_3};\POS(18,0)*{_4};
{\ar_{}"v3";"v2" },
{\ar_{}"v2";"v1" },
{\ar_{}"v4";"v2" },
{\ar^{}"v3";"v1" },
{\ar^{}"v4";"v1" },
 \endxy
  ~~  ~~
&
 \xy ;/r0.2pc/:
\POS(0,0)*+{\bullet}="v1"; 
\POS(0,15)*+{\bullet}="v2";
\POS(15,15)*+{\bullet}="v3";
\POS(15,0)*+{\bullet}="v4";
\POS(-3,0)*{_1};\POS(-3,17)*{_2};\POS(17,17)*{_3};\POS(18,0)*{_4};
{\ar_{}"v2";"v3" },
{\ar_{}"v2";"v1" },
{\ar_{}"v4";"v2" },
{\ar_{}"v4";"v3" },
{\ar^{}"v4";"v1" },
 \endxy
  ~~  ~~
&
 \xy ;/r0.2pc/:
\POS(0,0)*+{\bullet}="v1"; 
\POS(0,15)*+{\bullet}="v2";
\POS(15,15)*+{\bullet}="v3";
\POS(15,0)*+{\bullet}="v4";
\POS(-3,0)*{_1};\POS(-3,17)*{_2};\POS(17,17)*{_3};\POS(18,0)*{_4};
{\ar_{}"v2";"v3" },
{\ar_{}"v2";"v1" },
{\ar_{}"v4";"v3" },
{\ar^{}"v3";"v1" },
{\ar^{}"v4";"v1" },
 \endxy
  ~~  ~~
&
 \xy ;/r0.2pc/:
\POS(0,0)*+{\bullet}="v1"; 
\POS(0,15)*+{\bullet}="v2";
\POS(15,15)*+{\bullet}="v3";
\POS(15,0)*+{\bullet}="v4";
\POS(-3,0)*{_1};\POS(-3,17)*{_2};\POS(17,17)*{_3};\POS(18,0)*{_4};
{\ar_{}"v3";"v2" },
{\ar_{}"v4";"v2" },
{\ar_{}"v4";"v3" },
{\ar^{}"v3";"v1" },
{\ar^{}"v4";"v1" },
 \endxy ~~  ~~
\\
\scriptstyle{B^+(2)}
&\scriptstyle{B^-(2)}
&\scriptstyle{B^+(3)}
&\scriptstyle{B^-(3)}
 \end{array}
\end{align*}

\noindent
All edge arrows of the poset $A$ are minimal.
The posets $B$ and $D$ are of height 2. 
Their differential biquivers
are uniquely defined
by the rule as follows: for any triangle $i\prec j\prec k$  
with arrows $a:j\mo i,$ $b:k\mo i,$ and $c:k\mo j$
the arrow $b$ is minimal  
$\d(a)=\sum_{k}b v_c+ \phi$ and $\d(c)=-\sum_{i}v_a b+\psi,$
where $\phi$ and $\psi$ are summands coming from the other triangles.
That is, in a triangle the differential is as in (\ref{diff_biq_B}) of type $B.$
For the poset $C$ 
one should additionally give a pair of minimal edges:
for $C(2,3)$ they are $3\mo 2$ and $4\mo 1,$
for $C(3,2)$ they are $2\mo 3$ and $4\mo 1.$
The differentials of the other arrows consist of the paths of 
length 3 and degree one.

\subsection*{Rank and degree}
For configurations of types $A,$ $C$ and $B$ 
let $I_{min}\subset I$ 
and $I_{max}\subset I$
be the subsets of minimal and maximal vertices 
with respect to the partial order $\prec .$
For a dimension vector $\sss\in\mN^4$ let us introduce 
new discrete parameters $(\alpha,\beta)$:  
\begin{itemize}
  \item[--] for a box of type either $A$ or $C$
  define $\alpha:= \mathop\sum\limits_{i\in I_{min}} s_i$
  and $\beta:= \mathop\sum\limits_{k\in I_{max}} s_k;$

  \item[--] for a box  $B^{\sigma}(j)$  define
  define $\alpha:= s_j+\mathop\sum\limits_{i\in I_{min}} s_i$
  and $\beta:= s_j+ \mathop\sum\limits_{k\in I_{max}} s_k.$
\end{itemize}

\begin{lemma}
\label{lemma_rd}
Let $\Gamma$ be the automaton either (\ref{autom_I3}) or (\ref{autom_IV})
and $p:\kA\mo \kA'$ a path on it connecting principal states 
$\kA$ and $\kA'$  and taking $\sss\to \sss'$ and $(\alpha,\beta)\to(\alpha',\beta').$
Then $\gcd(\alpha,\beta)=\gcd(\alpha',\beta').$
\end{lemma}
\begin{proof}
It is sufficient to  
check the statement on the shortest paths.
For the transitions 
  $A\mo B$  
   or  $C\mo B$ 
   we have $(\alpha',\beta')=(\alpha,\beta).$
For the transitions 
      $A\mo A,$ $B\mo C$ or a path of length two $B\mo A,$ we have
   $$(\alpha',\beta')=\left\{%
\begin{array}{ll}
    (\alpha-\beta,\beta), & \hbox{if }\; \alpha \geq \beta; \\
    (\alpha,\beta-\alpha), & \hbox{if otherwise.} \\
\end{array}%
\right. 
$$
That completes the proof.
\end{proof}

\noindent
To obtain the statement of the 
Theorem \ref{theo_main} for
Kodaira fibers I$_3$ and IV 
we should replace the pair $(\alpha,\beta)$ by the rank and degree $(r,\bd)$ using 
Tables \ref{table_I_3} and \ref{table_IV}.  
In each case 
but cases 2 and $2'$ of Table \ref{table_IV} we have 
$(\alpha,\beta)= (r-d\bmod r,d\bmod r).$
In the cases 2 and $2'$ we have 
respectively 
$(\alpha+\beta,\beta)= (r-d\bmod r,d\bmod r)$
and $(\alpha,\alpha+\beta)= (r-d\bmod r,d\bmod r).$

\begin{remk}
As in Remark \ref{remk_FM} for singular Weierstra\ss{} curves,
for reduced curves we also can interpret the action of
small reductions on discrete parameters in terms 
of the action of Fourier-Mukai transforms.   
For example, under the condition $\bd\leq r,$ the action of the path
$(13)(12)\colon A^+\mo A^+$  of the automaton (\ref{autom_I2})
and paths $(14)(13)(12)\colon A^+\mo A^+$ of the automatons (\ref{autom_I3})
and (\ref{autom_IV}) correspond to the action of the
Seidel-Thomas twist $\T_{\oo}: (r,\bd_1,\bd_2)\to (r-\bd, \bd_1, \bd_2)$ 
if $E$ has two components and 
$\T_{\oo}: (r,\bd_1,\bd_2,\bd_3)\to (r-\bd, \bd_1, \bd_2,\bd_3),$
for three components.
In general a comparison of both actions on rank and multidegree
is bulky, and it would not help much for understanding, as far as 
we have no theorems relating small reductions 
and Fourier-Mukai transforms on genus one curves. 
We hope that small reductions $p: \kA\mo \kA' $ can be lifted 
to the level of derived equivalences 
$p\colon \D^b(\kA)\stackrel{\sim}\longrightarrow \D^b(\kA'),$
where $\D^b(\kA)$ is the derived category associated to the box $\kA,$
as introduced by Ovsienko in \cite{ovs97}.
There is a certain analogy with Atiyah's bijections \cite{Ati57}
for vector bundles on elliptic curves, which obtained a conceptual explanation 
in terms of derived categories in \cite{LM} few decades later 
after their discovery.
\end{remk}


%% file: exams.tex
\section{Examples and remarks}
\label{sec_exam}
\begin{exam}
  Let $E$ be a curve from the list with 2 components 
i.e. the Kodaira cycle I$_2$ or the fiber III.
Let us describe vector bundles on $E$ of rank $r=9$ and multidegree
$(d_1,d_2)=(3,2)$ using algorithm \ref{algorithm}.
The normalization bundle $\fn$ is 
$$ 
 {\fn|_{L_1}=\oo_{L_1}^6
 \oplus\big(\oo_{L_1}(1) \big)^{3}}
 \hbox{ and  }\;
 {\fn|_{L_2}=\oo_{L_2}^{7}
 \oplus\big(\oo_{L_2}(1) \big)^{2}.}
$$

\begin{enumerate}
\item
Since $\bd=d= 5 < 9=r$ thus 
according to table \ref{table_I_2}
the input state
for the automaton is $A^+$
and the dimension vector is
$\sss=(s_1,s_2,s_3)=(4,2,3).$

\item  Taking on automaton (\ref{autom_cycle2}) the path
$
p: A^+\stackrel{(12)}\lar B \stackrel{(31)}\lar A^-
\stackrel{(23)}\lar A^-\stackrel{(23)}\lar A^- 
\stackrel{(32)}\lar B \stackrel{(13)}\lar A^+
\stackrel{(12)}\lar B\stackrel{(31)}\lar A^-
$
we get the reduction of sizes:
\begin{align*}
(4,2,3)
\stackrel{(12)}\longmapsto (2,2,3)
\stackrel{(31)}\longmapsto (2,2,1)
\stackrel{(23)}\longmapsto (2,1,1)
\stackrel{(23)}\longmapsto (2,0,1)
\\
\stackrel{(32)}\longmapsto (2,0,1)
\stackrel{(13)}\longmapsto (1,0,1)
\stackrel{(12)}\longmapsto (1,0,1)
\stackrel{(31)}\longmapsto (1,0,0).
\end{align*}
Since $s_2=0$ 
reductions $(1,0,1) \stackrel{(12)}\longmapsto (1,0,1)$
and $(2,0,1) \stackrel{(32)}\longmapsto (2,0,1)$ are degenerated, 
that is we formally change states, but the matrices remain the same. 

\item Reversing the path $p$ we construct a canonical form of the matrix
\linebreak ${M\in \Br_{\kA}(4,2,3)}.$ 
\end{enumerate}

If $E$ is a Kodaira cycle I$_2$ then
\begin{align*}
\begin{array}{@{}c@{}l@{}}
\begin{array}{@{}c@{}}
\smk{_1}
\end{array}
&
\\
\begin{array}{ |@{}c @{}|@{}}
\hline
 {\mk{\lambda }}    \\
\hline
\end{array}
&
\begin{array}{@{}l@{}}
\tmk{\scriptstyle{1}}\\
\end{array}
\\ \phantom{a}
\end{array}
\;\;\stackrel{(31)(12)}\longmapsto
\;\;
\begin{array}{@{}c@{}l@{}}
\begin{array}{@{}c@{}c@{}}
\smk{_1}&\smk{_3}
\end{array}
&
\\
\begin{array}{ |@{}c @{}|@{} c @{}|@{}}
\hline
{\mk{0}}      &   \mk{1}    \\
\hline
 {\mk{\lambda }}   &     \mk{0}    \\
\hline
\end{array}
&
\begin{array}{@{}l@{}}
\tmk{\scriptstyle{1}}\\\tmk{\scriptstyle{3}}
\end{array}
\\ \phantom{a}
\end{array}
\;\;\stackrel{(13)(32)}\longmapsto
\;\;
\begin{array}{@{}c@{}c@{}}
\begin{array}{@{}c@{}c@{}c}
\smk{_1}&\smk{}&\smk{_3}
\end{array}
&
\\
\begin{array}{ |@{}c @{} c  @{}|@{}    c   @{}|@{}}
\hline
{\mk{0}}      &\mk{1}      &     \mk{0}    \\
 {\mk{0}}     &{\mk{0}}    &     \mk{1}      \\
\hline
 {\mk{\lambda}}       &{\mk{0}}    &     \mk{0}    \\
\hline
\end{array}
&
\begin{array}{@{}l@{}}
\tmk{\scriptstyle{1}}\\\tmk{\scriptstyle{}}\\\tmk{\scriptstyle{3}}
\end{array}
\\ \phantom{a}
\end{array}
\;\stackrel{(23)}\longmapsto
\;
\begin{array}{@{}c@{}l@{}}
\begin{array}{@{}c@{}c@{}c@{}c}
\smk{_1}&\smk{}&\smk{_2}&\smk{_3}
\end{array}
&
\\
\begin{array}{ |@{}c @{} c  @{}|@{}    c   @{}|@{}  c   @{}|@{}}
\hline
{\mk{0}}      &\mk{0}     &\mk{0}      &     \mk{1}    \\
\hline
{\mk{0}}      &\mk{1}     &\mk{0}      &     \mk{0}    \\
{\mk{0}}      &\mk{0}     &\mk{1}      &     \mk{0}    \\
\hline
{\mk{\lambda}}      &\mk{0}     &\mk{0}    &     \mk{0}    \\
\hline
\end{array}
&
\begin{array}{@{}l@{}}
\tmk{\scriptstyle{2}}\\
\tmk{\scriptstyle{1}}\\\tmk{}
\\\tmk{\scriptstyle{3}}
\end{array}
\\ \phantom{a}
\end{array}
\\
\;\stackrel{(23)}\longmapsto
\;
\begin{array}{@{}c@{}l@{}}
\begin{array}{@{}c@{}c@{}c@{}c@{}c}
\mk{_1}&\mk{}&\mk{_2}&\mk{}&\mk{_3}
\end{array}
&
\\
\begin{array}{ |@{}c @{} c  @{}|@{}    c   @{}c   @{}|@{}  c   @{}|@{}}
\hline
\mk{0}&\mk{0}     &\mk{0} &\mk{1}      &\mk{0}    \\
\mk{0}&\mk{0}     &\mk{0} &\mk{0}      &\mk{1}    \\
\hline
\mk{0}&\mk{1}     &\mk{0} &\mk{0}      &\mk{0}    \\
\mk{0}&\mk{0}     &\mk{1} &\mk{0}      &\mk{0}    \\
\hline
{\mk{\lambda}}&\mk{0} &\mk{0}&\mk{0}    &\mk{0}    \\
\hline
\end{array}
&
\begin{array}{@{}l@{}}
\tmk{\scriptstyle{2}}\\\tmk{}
\\\tmk{\scriptstyle{1}}\\\tmk{}
\\\tmk{\scriptstyle{3}}
\end{array}
\\ \phantom{a}
\end{array}
\;\stackrel{(31)}\longmapsto
\;
\begin{array}{@{}c@{}l@{}}
\begin{array}{@{}c@{}c@{}c@{}c@{}c@{}c@{}c@{}}
\smk{_1}&\smk{}&\smk{_2}&\smk{}&\smk{}&\smk{_3}\smk{}
\end{array}
&
\\
\begin{array}{ |@{}c@{}c@{}|@{}    c@{}c @{}|@{}  c@{}c @{}c   @{}|@{}}
\hline
\mk{0}&\mk{0}     &\mk{0} &\mk{0}      &\mk{1}&\mk{0}&\mk{0}    \\
\mk{0}&\mk{0}     &\mk{0} &\mk{0}      &\mk{0}&\mk{1}&\mk{0}    \\
\hline
\mk{0}&\mk{0}     &\mk{0} &\mk{1}      &\mk{0}&\mk{0}&\mk{0}    \\
\mk{0}&\mk{0}     &\mk{0} &\mk{0}      &\mk{0}&\mk{0}&\mk{1}    \\
\hline
\mk{0}&\mk{1}     &\mk{0} &\mk{0}      &\mk{0}&\mk{0}&\mk{0}    \\
\mk{0}&\mk{0}     &\mk{1} &\mk{0}      &\mk{0}&\mk{0}&\mk{0}    \\
\mk{\lambda}&\mk{0} &\mk{0}&\mk{0}     &\mk{0}&\mk{0}&\mk{0}    \\
\hline
\end{array}
&
\begin{array}{@{}l@{}}
\tmk{\scriptstyle{1}}\\\tmk{}
\\\tmk{\scriptstyle{2}}\\\tmk{}\\
\tmk{}
\\\tmk{\scriptstyle{3}}\\\tmk{}
\end{array}
\\ \phantom{a}
\end{array}
\;\stackrel{(12)}\longmapsto
\;
\begin{array}{@{}c@{}l@{}}
\begin{array}{@{}c@{}c@{}c@{}c@{}  c@{}c@{}c  @{}c@{}c@{}}
\smk{}\smk{_1}&\smk{}&\smk{}
&\smk{}&\smk{_3}&\smk{}
&\smk{_2}&\smk{}
\end{array}
&
\\
\begin{array}{ |@{}c@{}c@{}c@{}c@{}|@{}    c@{}c@{}c  @{}|@{} c@{}c @{}|@{}}
\hline
\mk{0}&\mk{0}&\mk{0}&\mk{0}   &\mk{1}&\mk{0}&\mk{0}  &\mk{0}&\mk{0}   \\
\mk{0}&\mk{0}&\mk{0}&\mk{0}   &\mk{0}&\mk{1}&\mk{0}  &\mk{0}&\mk{0}   \\
\mk{0}&\mk{0}&\mk{0}&\mk{0}   &\mk{0}&\mk{0}&\mk{0}  &\mk{1}&\mk{0}   \\
\mk{0}&\mk{0}&\mk{0}&\mk{0}   &\mk{0}&\mk{0}&\mk{0}  &\mk{0}&\mk{1}   \\
\hline
\mk{0}&\mk{0}&\mk{0}&\mk{1}   &\mk{0}&\mk{0}&\mk{0}  &\mk{0}&\mk{0}   \\
\mk{0}&\mk{0}&\mk{0}&\mk{0}   &\mk{0}&\mk{0}&\mk{1}  &\mk{0}&\mk{0}   \\
\hline
\mk{0}&\mk{1}&\mk{0}&\mk{0}   &\mk{0}&\mk{0}&\mk{0}   &\mk{0}&\mk{0}   \\
\mk{0}&\mk{0}&\mk{1}&\mk{0}   &\mk{0}&\mk{0}&\mk{0}   &\mk{0}&\mk{0}   \\
\mk{\lambda}&\mk{0} &\mk{0}&\mk{0}   &\mk{0}&\mk{0}&\mk{0}    &\mk{0}&\mk{0}   \\
\hline
\end{array}
&
\begin{array}{@{}l@{}}
\tmk{}\\\tmk{\scriptstyle{1}}\\\tmk{}\\\tmk{}\\
\tmk{\scriptstyle{2}}\\\tmk{}
\\ \tmk{} \\\tmk{\scriptstyle{3}}\\\tmk{}
\end{array}
\\ \phantom{a}
\end{array}
&&
\end{align*}
\noindent
Let us construct the canonical form 
for the tacnode curve. Besides zeros we also use the  
empty spaces to mark out the blocks $(ij)$, where zeros appear
for some general reasons and the corresponding box contains 
no arrow $j\mo i.$
Note that the order of row and column blocks are chosen 
in such a way that the matrices have block triangular form
(probably with some additional holes). 
\begin{align*}
\begin{array}{@{}c@{}l@{}}
\begin{array}{@{}c@{}}
\smk{_1}
\end{array}
&
\\
\begin{array}{ |@{}c @{}|@{}}
\hline
 {\mk{\lambda }}    \\
\hline
\end{array}
&
\begin{array}{@{}l@{}}
\tmk{\scriptstyle{1}}\\
\end{array}
\\ \phantom{a}
\end{array}
\;\;
\stackrel{(31)(12)}
\longmapsto
\;\;
\begin{array}{@{}c@{}l@{}}
\begin{array}{@{}c@{}c@{}}
\smk{_1}&\smk{_3}
\end{array}
&
\\
\begin{array}{ |@{}c @{}|@{} c @{}|@{}}
\hline
{\mk{\lambda }}      &   \mk{1}    \\
\hline
 {\mk{}}   &     \mk{0}    \\
\hline
\end{array}
&
\begin{array}{@{}l@{}}
\tmk{\scriptstyle{1}}\\\tmk{\scriptstyle{3}}
\end{array}
\\ \phantom{a}
\end{array}
\;\;\stackrel{(13)(32)}\longmapsto
\;\;
\begin{array}{@{}c@{}c@{}}
\begin{array}{@{}c@{}c@{}c}
\smk{_1}&\smk{}&\smk{_3}
\end{array}
&
\\
\begin{array}{ |@{}c @{} c  @{}|@{}    c   @{}|@{}}
\hline
{\mk{\lambda}}&\mk{1}      &     \mk{0}    \\
 {\mk{0}}&{\mk{0}}         &     \mk{1}      \\
\hline
 {\mk{}}&{\mk{}}           &     \mk{0}    \\
\hline
\end{array}
&
\begin{array}{@{}l@{}}
\tmk{\scriptstyle{1}}\\\tmk{\scriptstyle{}}\\\tmk{\scriptstyle{3}}
\end{array}
\\ \phantom{a}
\end{array}
\;\stackrel{(23)}\longmapsto
\;
\begin{array}{@{}c@{}l@{}}
\begin{array}{@{}c@{}c@{}c@{}c}
\smk{_2}&\smk{_1}&\smk{}&\smk{_3}
\end{array}
&
\\
\begin{array}{ |@{}c@{}|@{}c@{}c@{}|@{}c@{}|@{}}
\hline
{\mk{0}}&        \mk{}&\mk{}        &\mk{1}    \\
\hline
\mk{}&    {\mk{\lambda}}&\mk{1}     &\mk{0}    \\
\mk{}&    {\mk{0}}&\mk{0}           &\mk{1}    \\
\hline
\mk{}&    \mk{}&\mk{}               &\mk{0}    \\
\hline
\end{array}
&
\begin{array}{@{}l@{}}
\tmk{\scriptstyle{2}}\\
\tmk{\scriptstyle{1}}\\ \tmk{} \\
\tmk{\scriptstyle{3}}
\end{array}
\\ \phantom{a}
\end{array}
\\
\;\stackrel{(23)}\longmapsto
\;
\begin{array}{@{}c@{}l@{}}
\begin{array}{c@{}c@{}c@{}c@{}c@{}}
\smk{_2}&\smk{}&   \smk{_1}&\smk{}& \smk{_3}
\end{array}
&
\\
\begin{array}{ |@{}c@{}c @{}|@{}c@{}c @{}|@{}c@{}|@{}}
\hline
\mk{0}&\mk{1}   &\mk{} &\mk{}      &\mk{0}    \\
\mk{0}&\mk{0}   &\mk{} &\mk{}      &\mk{1}    \\
\hline
\mk{}&\mk{}     &\mk{\lambda}&\mk{1} &\mk{0}    \\
\mk{}&\mk{}     &\mk{0} &\mk{0}      &\mk{1}    \\
\hline
\mk{}&\mk{}     &\mk{}&\mk{}    &\mk{0}    \\
\hline
\end{array}
&
\begin{array}{@{}l@{}}
\tmk{\scriptstyle{2}}\\\tmk{}
\\\tmk{\scriptstyle{1}}\\\tmk{}
\\\tmk{\scriptstyle{3}}
\end{array}
\\ \phantom{a}
\end{array}
\;\stackrel{(31)}\longmapsto
\;
\begin{array}{@{}c@{}l@{}}
\begin{array}{@{}c@{}c@{}c@{}c@{}c@{}c@{}c@{}}
\smk{_1}&\smk{}   &\smk{_2}&\smk{}   &\smk{}&\smk{_3}\smk{}
\end{array}
&
\\
\begin{array}{ |@{}c@{}c@{}|@{}    c@{}c @{}|@{}  c@{}c @{}c   @{}|@{}}
\hline
\mk{0}&\mk{0}     &\mk{0} &\mk{0}      &\mk{1}&\mk{0}&\mk{0}    \\
\mk{0}&\mk{0}     &\mk{0} &\mk{0}      &\mk{0}&\mk{1}&\mk{0}    \\
\hline
\mk{ }&\mk{ }     &\mk{0} &\mk{1}      &\mk{0}&\mk{0}&\mk{0}    \\
\mk{ }&\mk{ }     &\mk{0} &\mk{0}      &\mk{0}&\mk{0}&\mk{1}    \\
\hline
\mk{ }&\mk{ }     &\mk{ } &\mk{ }      &\mk{\lambda}&\mk{1}&\mk{0}    \\
\mk{ }&\mk{ }     &\mk{ } &\mk{ }      &\mk{0}&\mk{0}&\mk{1}    \\
\mk{ }&\mk{ }     &\mk{ } &\mk{ }     &\mk{0}&\mk{0}&\mk{0}    \\
\hline
\end{array}
&
\begin{array}{@{}l@{}}
\tmk{\scriptstyle{1}}\\\tmk{}
\\\tmk{\scriptstyle{2}}\\\tmk{}\\
\tmk{}
\\\tmk{\scriptstyle{3}}\\\tmk{}
\end{array}
\\ \phantom{a}
\end{array}
\;\stackrel{(12)}\longmapsto
\;
\begin{array}{@{}c@{}l@{}}
\begin{array}{@{}c@{}c@{}c@{}c@{}  c@{}c@{}c  @{}c@{}c@{}}
\smk{}\smk{_1}&\smk{}&\smk{}
&\smk{}&\smk{_3}&\smk{}
&\smk{_2}&\smk{}
\end{array}
&
\\
\begin{array}{ |@{}c@{}c@{}c@{}c@{}|@{}    c@{}c@{}c  @{}|@{} c@{}c @{}|@{}}
\hline
\mk{0}&\mk{0}&\mk{0}&\mk{0}   &\mk{1}&\mk{0}&\mk{0}  &\mk{0}&\mk{0}   \\
\mk{0}&\mk{0}&\mk{0}&\mk{0}   &\mk{0}&\mk{1}&\mk{0}  &\mk{0}&\mk{0}   \\
\mk{0}&\mk{0}&\mk{0}&\mk{1}   &\mk{0}&\mk{0}&\mk{0}  &\mk{1}&\mk{0}   \\
\mk{0}&\mk{0}&\mk{0}&\mk{0}   &\mk{0}&\mk{0}&\mk{1}  &\mk{0}&\mk{1}   \\
\hline
\mk{ }&\mk{ }&\mk{ }&\mk{ }   &\mk{\lambda}&\mk{1}&\mk{0}   &\mk{ }&\mk{ }   \\
\mk{ }&\mk{ }&\mk{ }&\mk{ }   &\mk{0}&\mk{0}&\mk{1}   &\mk{ }&\mk{ }   \\
\mk{ }&\mk{ }&\mk{ }&\mk{ }   &\mk{0}&\mk{0}&\mk{0}   &\mk{ }&\mk{ }   \\
\hline
\mk{ }&\mk{ }&\mk{ }&\mk{ }   &\mk{ }&\mk{ }&\mk{ }  &\mk{0}&\mk{0}   \\
\mk{ }&\mk{ }&\mk{ }&\mk{ }   &\mk{ }&\mk{ }&\mk{ }  &\mk{0}&\mk{0}   \\
\hline
\end{array}
&
\begin{array}{@{}l@{}}
\tmk{}\\\tmk{\scriptstyle{1}}\\\tmk{}\\\tmk{}\\
\tmk{\scriptstyle{3}}\\\tmk{}
\\ \tmk{} \\\tmk{\scriptstyle{2}}\\\tmk{}
\end{array}
\\ \phantom{a}
\end{array}
&
\end{align*}
\end{exam}

\begin{remk}
For a Kodaira fiber II, III and IV 
the parameter $\lambda$ 
of the canonical form of $M(\lambda)$ can be moved to any 
place on the diagonal, as well as it can be distributed 
as ${\frac{\lambda}{r}}$
to all the diagonal entries.
This way the canonical form resembles to the Jordan normal form.
For instance in the last example we get:
\begin{align}
\label{form_diag}
\begin{array}{ |@{}c@{}c@{}c@{}c@{}|@{}    c@{}c@{}c  @{}|@{} c@{}c @{}|@{}}
\hline
\mk{\frac{\lambda}{9}}&\mk{0}&\mk{0}&\mk{0}   &\mk{1}&\mk{0}&\mk{0}  &\mk{0}&\mk{0}   \\
\mk{0}&\mk{\frac{\lambda}{9}}&\mk{0}&\mk{0}   &\mk{0}&\mk{1}&\mk{0}  &\mk{0}&\mk{0}   \\
\mk{0}&\mk{0}&\mk{\frac{\lambda}{9}}&\mk{1}   &\mk{0}&\mk{0}&\mk{0}  &\mk{1}&\mk{0}   \\
\mk{0}&\mk{0}&\mk{0}&\mk{\frac{\lambda}{9}}   &\mk{0}&\mk{0}&\mk{1}  &\mk{0}&\mk{1}   \\
\hline
\mk{ }&\mk{ }&\mk{ }&\mk{ }   &\mk{\frac{\lambda}{9}}&\mk{1}&\mk{0}   &\mk{ }&\mk{ }   \\
\mk{ }&\mk{ }&\mk{ }&\mk{ }   &\mk{0}&\mk{\frac{\lambda}{9}}&\mk{1}   &\mk{ }&\mk{ }   \\
\mk{ }&\mk{ }&\mk{ }&\mk{ }   &\mk{0}&\mk{0}&\mk{\frac{\lambda}{9}}   &\mk{ }&\mk{ }   \\
\hline
\mk{ }&\mk{ }&\mk{ }&\mk{ }   &\mk{ }&\mk{ }&\mk{ }  &\mk{\frac{\lambda}{9}}&\mk{0}   \\
\mk{ }&\mk{ }&\mk{ }&\mk{ }   &\mk{ }&\mk{ }&\mk{ }  &\mk{0}&\mk{\frac{\lambda}{9}}   \\
\hline
\end{array}
\end{align}
\end{remk}